\documentclass[reqno,a4paper]{amsart}
\usepackage[utf8]{inputenc}
\usepackage[T1]{fontenc}
\usepackage[british]{babel}
\usepackage{amssymb,amsthm,amsmath}
\usepackage{calrsfs}
\DeclareSymbolFontAlphabet{\mathcalstd}{symbols}
\usepackage[colorlinks]{hyperref}
\usepackage{mathbbol}
\usepackage{mathabx}
\usepackage[all]{xy}
\usepackage{xcolor}
\usepackage{tikz-cd}
\usepackage{mathtools}
\usepackage{enumitem}
\usepackage{booktabs}
\usepackage{microtype}
\allowdisplaybreaks
\renewcommand{\binom}{\tbinom} 

\theoremstyle{plain}
\newtheorem{theorem}[subsection]{Theorem}
\newtheorem{lemma}[subsection]{Lemma}
\newtheorem{proposition}[subsection]{Proposition}
\newtheorem{corollary}[subsection]{Corollary}
\newtheorem*{maintheorem}{Main Theorem}

\theoremstyle{definition}
\newtheorem{definition}[subsection]{Definition}
\newtheorem{example}[subsection]{Example}
\newtheorem{examples}[subsection]{Examples}
\newtheorem{remark}[subsection]{Remark}

\newtheorem{convention}[subsection]{Convention}
\newtheorem{question}[subsection]{Question}

\newlist{tfae}{enumerate}{1}
\setlist[tfae]{label=(\roman*), ref=(\roman*)}
\newlist{sqaxioms}{enumerate}{1}
\setlist[sqaxioms]{label=(SR\arabic*), ref=(SR\arabic*), leftmargin=*, widest=(SR8)}
\newlist{modaxioms}{enumerate}{1}
\setlist[modaxioms]{label=(M\arabic*), ref=(M\arabic*), leftmargin=*, widest=(M7)}
\setlist[enumerate,1]{label=(\arabic*), ref=(\arabic*)}

\tikzcdset{arrow style=tikz, diagrams={>=stealth}}

\makeatletter
\newcommand*{\comment}[4][]{%
  \@ifundefined{c@#2}{\newcounter{#2}}{}%
  \begingroup
  \ifx\\#1\\\newcommand*{\header}{\csname the#2\endcsname}\else\newcommand*{\header}{(#1\csname the#2\endcsname)}\fi
  \addtocounter{#2}{1}%
  \textnormal{\textcolor{#3}{\header.[#4\@]}}%
  \endgroup
}
\makeatother

\newcommand{\defn}{\textbf}
\newcommand{\noproof}{\hfill\qed}
\newcommand{\DefEq}{\coloneq}

\newcommand{\C}{\ensuremath{\mathcal{C}}}
\newcommand{\D}{\ensuremath{\mathcal{D}}}
\newcommand{\V}{\ensuremath{\mathcal{V}}}
\newcommand{\A}{\ensuremath{\mathcal{A}}}
\newcommand{\R}{\ensuremath{\mathcalstd{R}}}

\newcommand{\Z}{\ensuremath{\mathbb{Z}}}
\newcommand{\RP}{\ensuremath{\mathbb{R}\mathrm{P}}}
\newcommand{\Ab}{\ensuremath{\mathsf{Ab}}}
\newcommand{\Grp}{\ensuremath{\mathsf{Grp}}}
\newcommand{\Gr}{\ensuremath{\mathsf{Gr}}}
\newcommand{\Set}{\ensuremath{\mathsf{Set}}}
\newcommand{\Mod}{\ensuremath{\mathsf{Mod}}}
\newcommand{\Nil}{\ensuremath{\mathsf{Nil}}}
\newcommand{\Func}{\ensuremath{\mathsf{Func}}}
\newcommand{\Free}{\ensuremath{\mathsf{Free}}}
\newcommand{\SqRing}{\ensuremath{\mathsf{SqRing}}}

\newcommand{\CRing}{\ensuremath{\mathsf{CRing}}}

\renewcommand{\P}{\ensuremath{\mathcal{P}}}

\newcommand{\cre}{\ensuremath{\mathrm{cr}}}
\DeclareMathOperator{\ab}{ab}
\DeclareMathOperator{\nil}{nil}
\DeclareMathOperator{\Ker}{Ker}
\DeclareMathOperator{\Coker}{Coker}
\DeclareMathOperator{\coker}{coker}
\DeclareMathOperator{\ev}{ev}

\renewcommand{\Im}{\ensuremath{\mathrm{Im}}}
\newcommand{\U}{\ensuremath{U}}

\newcommand{\cosmash}{\ensuremath{\diamond}}
\newcommand{\cosum}{\ensuremath{\sqcup}}                   
\newcommand{\cm}[2]{\ensuremath{\hat{r}_{#1,#2}}}          
\newcommand{\cmF}[2]{\ensuremath{\hat{r}^{#1}_{#2}}}       
\newcommand{\cc}[1]{\ensuremath{c_{#1}}}                    
\newcommand{\tw}[2]{\ensuremath{\tau_{#1,#2}}}             
\newcommand{\twp}[2]{\ensuremath{\sigma_{#1,#2}}}          
\newcommand{\swc}[1]{\ensuremath{\tau^{\cosmash}_{#1}}}          
\newcommand{\fold}[1]{\ensuremath{\nabla_{#1}}}            
\newcommand{\diag}[1]{\ensuremath{\Delta_{#1}}}            

\newcommand{\sqring}[4]{\ensuremath{(#1\xrightarrow{#3}#2\xrightarrow{#4}#1)}}
\renewcommand{\Re}{\ensuremath{\R_{e}}}
\newcommand{\Ree}{\ensuremath{\R_{ee}}}
\newcommand{\Rbar}{\ensuremath{\overline{R}}}
\newcommand{\brk}[2]{\ensuremath{[#1,#2]}}

\title[Two-step nilpotent varieties and square rings]{Two-step nilpotent varieties as\\ categories of modules over a square ring}
\author{L.~Geenis}
\author{M.~Hartl}
\author{T.~Van der Linden}

\email[L.~Geenis]{logan.geenis@uclouvain.be}
\email[M.~Hartl]{mhartl@gmx.fr}
\email[T.~Van der Linden]{tim.vanderlinden@uclouvain.be}
\email[T.~Van der Linden]{tim.van.der.linden@vub.be}

\address[L.~Geenis, T.~Van der Linden]{Institut de Recherche en Math\'ematique et Physique, Universit\'e catholique de Louvain, che\-min du cyclotron~2 bte~L7.01.02, B--1348 Lou\-vain-la-Neuve, Belgium}
\address[M.~Hartl, T.~Van der Linden]{Mathematics \& Data Science, Vrije Universiteit Brussel, Pleinlaan 2, B--1050 Brussel, Belgium}

\thanks{The third author is a Research Director of the Fonds de la Recherche Scientifique--FNRS}
\date{}

\begin{document}

\begin{abstract}
  We characterise the two-step nilpotent varieties of algebras: the semi-abelian varieties in which all ternary cosmash products vanish. Just as abelian varieties are the categories of modules over a ring, two-step nilpotent varieties are exactly the categories of modules over a square ring in the sense of Baues, Hartl and Pirashvili. We construct the square ring explicitly from the full subcategory on the free algebras on one and two generators, and give complete proofs. As consequences we obtain a description of any two-step nilpotent variety by generators and relations, an intrinsic characterisation in terms of quadratic term operations, and an account of its abelian core and of the lower central series of its objects in terms of the square ring. Along the way, we associate with any commutative ring \(R\) a square ring whose modules, the two-step nilpotent \(R\)-groups, extend Hall's nilpotent \(R\)-groups of class two from \(2\)-binomial rings to arbitrary rings of coefficients.
\end{abstract}

\maketitle

\section*{Introduction}

Any abelian variety of algebras is equivalent to the category of (right) modules over a ring, namely the endomorphism ring of its free algebra on one generator~\cite{Freyd}. At the other end of the spectrum, semi-abelian varieties have been characterised syntactically: they are the varieties admitting a unique constant \(0\), binary operations \(\alpha_{1}\), \dots, \(\alpha_{n}\) and an \((n+1)\)-ary operation \(\theta\) subject to the identities \(\alpha_{i}(x,x)=0\) and \(\theta(\alpha_{1}(x,y),\dots,\alpha_{n}(x,y),y)=x\)~\cite{BoJa}; this description extends to infinitary and many-sorted varieties~\cite{Gran-Rosicky:Monadic}. Characterisations of this kind are useful because they allow us to reason about such varieties in general, providing a hands-on description which is easy to use in practice.

In this paper we prove a result of the same type for a class of varieties situated between the abelian and the semi-abelian ones: the \emph{two-step nilpotent} varieties. A semi-abelian category is two-step nilpotent when all of its ternary cosmash products vanish, or equivalently, when its identity functor is quadratic in the sense of polynomial functor calculus~\cite{Hartl-Vespa, HVdL, DHVdL}. Typical examples are the variety of nilpotent groups of class at most two and, more generally, the two-nilpotent core \(\Nil_{2}(\V)\) of any semi-abelian variety \(\V\). The role of the ring is now played by a \defn{square ring}, a structure introduced by Baues, Hartl and Pirashvili~\cite{BHP} in the context of metastable homotopy theory~\cite{Baues:Metastable}, further developed within the quadratic algebra of square groups~\cite{Baues-Pirashvili, Baues-Pirashvili:Spaces, BJP:SquareGroups}, and later used by Baues and Muro in the study of secondary homotopy operations~\cite{Baues-Muro, Baues-Muro:Smash, Baues-Muro:RingSpectra}. A square ring \(\R\) consists of a nearring \(\Re\), an abelian group \(\Ree\) on which \(\Re\times\Re\) acts on the left and \(\Re\) on the right, quadrilinearly, an involution \(T\) of \(\Ree\) which plays the role of a symmetry, and maps \(H\colon \Re\to\Ree\) and \(P\colon \Ree\to\Re\) subject to a list of identities recalled in Definition~\ref{def:square-ring}; the identities determine \(T\) in terms of \(H\) and \(P\), which is why \(T\) is left out of the notation \(\sqring{\Re}{\Ree}{H}{P}\). Modules over \(\R\) are groups equipped with a unary operation \(m\mapsto m\cdot r\) for each \(r\in\Re\), which is a quadratic map in the sense of~\cite{Passi:GroupRings, Hartl:QuadraticMaps}, and a bilinear operation \((m,n)\mapsto\brk{m}{n}\cdot x\) for each \(x\in\Ree\), see Definition~\ref{def:module}. Throughout this paper, varieties of algebras are understood to be single-sorted and finitary. Our main result is the following.

\begin{maintheorem}
  A variety of algebras \(\V\) is two-step nilpotent if and only if there exists a square ring \(\R\) such that \(\V\) is equivalent to the category \(\Mod_{\R}\) of right \(\R\)-modules.
\end{maintheorem}

The square ring \(\R\) is obtained from the full subcategory of \V\ with objects being the free algebras  of \(\V\) on one and two generators: denoting the former by  \(E\), \(\Re\) is the endomorphism nearring \(\V(E,E)\), \(\Ree\) is the abelian group \(\V(E,E\cosmash E)\) of morphisms into the cosmash product \(E\cosmash E\) (which correspond to certain morphisms \(E\to E+E\)), and the maps \(H\), \(P\) are induced by a chosen comultiplication \(\mu_{E}\colon E\to E+E\), a lifting of the diagonal \(\Delta_E=(1_E,1_E)\colon E\to E\times E\) along the comparison morphism \(\cm{E}{E}\colon E+E\to E\times E\), and by the folding morphism \(\nabla_E=\langle 1_E,1_E\rangle\colon E+E\to E\). The idea that the algebraic theory of a two-step nilpotent variety carries the structure of a square ringoid goes back to the thesis of Defourneau~\cite[Proposition~4.1.24]{Defourneau}. In subsequent work, we hope to provide a different proof of this fact based on a ``nilpotency tower'' associated with semi-abelian categories in terms of linear extensions of categories; compare~\cite{Jibladze-Pirashvili:Linear}, where nilpotent Mal'tsev theories are studied by means of linear extensions, and~\cite{BeBou}. Since we only need the one-object case, and since we wish to keep the paper self-contained, we give a complete proof, constructing \(H\), \(P\), the involution \(T\) and the actions explicitly; see Remark~\ref{rem:defourneau}. Along the way we obtain normal forms for the free objects of \(\V\) and of \(\Mod_{\R}\), a description of \(\V\) by generators and relations (Corollary~\ref{cor:term-equivalence}) and an intrinsic characterisation of two-step nilpotent varieties in terms of quadratic term operations (Corollary~\ref{cor:quadratic-terms}). The result is compatible with the abelian case: the ring which classifies the abelian core \(\Ab(\V)\) is the cokernel of \(P\) while \(\Ree\) together with its symmetry and the three actions of the ring \(\Rbar= \Coker(P)\) on it encodes the cosmash product of \(\V\) and thus the bilinear product \cite{DHVdL} of any variety containing \(\V\) as its two-nilpotent core. These same data \((\Rbar,\Ree)\) also form a two-nilpotent reduced symmetric operad in the monoidal category of abelian groups which encodes the multilinear structure of the graded abelian group consisting of the two potentially non-trivial quotients of the lower central series of \(\R\)-modules and thus of the objects of \(\V\) (Section~\ref{sec:compatibility}).

The notion of a module over a square ring also addresses, for class \(2\), a problem in group theory which is independent of the Main Theorem. Hall's nilpotent \(R\)-groups of class \(2\)~\cite{Hall:Edmonton} are defined only when the ring of coefficients \(R\) is \(2\)-binomial, and the known extensions of this notion to arbitrary commutative rings lack some of its key features~\cite{MR94, IKK}. Following~\cite{P2RM}, we associate with any commutative ring \(R\) a square ring \(R_{\nil}\) which relates to \(R\) as the nil square ring \(\Z_{\nil}\) of~\cite{BHP}, whose modules are the groups of nilpotency class at most two, relates to \(\Z\). We call the modules over \(R_{\nil}\) \emph{two-step nilpotent \(R\)-groups}; they are Hall's nilpotent \(R\)-groups of class \(2\) whenever \(R\) is \(2\)-binomial (Example~\ref{ex:R-groups-any-ring}). The square rings which receive a morphism from \(R_{\nil}\), the \emph{square rings under \(R\)}, are then the quadratic counterpart of rings under \(R\): just as a module over a ring under \(R\) is an \(R\)-module, a module over a square ring under \(R\) is a two-step nilpotent \(R\)-group, by restriction of scalars. Being a square ring under \(R\) turns out to be a property of a square ring equipped with a morphism of nearrings from \(R\) (Proposition~\ref{prop:morphisms-from-R-nil}), and for \(R=\Z\) this recovers the fact that \(\Z_{\nil}\) is an initial object of the category of square rings.

\subsection*{Structure of the paper}
Section~\ref{sec:polynomial} recalls cross-effects, polynomial functors and the resulting notions of nilpotency in a semi-abelian category. Section~\ref{sec:square-rings} introduces square rings and their modules, with complete proofs of the basic facts we need, in particular a description of binary coproducts and of free modules; it also contains the examples, amongst which the square ring \(R_{\nil}\) of a commutative ring \(R\) and the square rings under \(R\). Section~\ref{sec:modules-nilpotent} proves that \(\Mod_{\R}\) is a two-step nilpotent variety for every square ring \(\R\). Section~\ref{sec:square-ring-of-V} constructs the square ring of a two-step nilpotent variety, and Section~\ref{sec:equivalence} proves the equivalence \(\V\simeq\Mod_{\R}\). Section~\ref{sec:syntactic} derives the syntactic characterisations, and Section~\ref{sec:compatibility} discusses compatibility with the abelian core and the bilinear product of a semi-abelian variety, as well as the global structure of the quotients of the lower central series of its two-step nilpotent objects, viewed as modules over the associated square ring.

\section{Polynomial functors and nilpotency}\label{sec:polynomial}

We recall what we need from polynomial functor calculus, following~\cite{Hartl-Vespa, Actions, HVdL, DHVdL}; for polynomial functors in related contexts, see also~\cite{Baues:Metastable, HPV}. Cross-effects of functors between abelian categories were introduced by Eilenberg and Mac~Lane~\cite{Eilenberg-MacLane}; the definition was extended to group-valued functors by Baues and Pirashvili~\cite{Baues-Pirashvili}, and to the present generality by Hartl and Loiseau~\cite{Actions}, who treat binary cross-effects, and by Hartl and Vespa~\cite{Hartl-Vespa}.

\subsection{Conventions and notation}\label{sec:conventions}
Throughout, \(\C\) is a pointed category with finite coproducts (also called \defn{sums}, and written \(X+Y\)) and \(\D\) is a pointed finitely complete category; from Subsection~\ref{sec:folding} on, both are semi-abelian~\cite{Janelidze-Marki-Tholen, WorkingMathematician}. We refer to~\cite{Borceux-Bourn} for the basic theory of semi-abelian categories. A functor \(F\colon \C\to\D\) is \defn{reduced} if \(F(0)=0\); we write \(\Func_{*}(\C,\D)\) for the category of reduced functors.

Given morphisms \(f_{k}\colon X_{k}\to Y\), we write \(\langle f_{1},\dots,f_{n}\rangle\colon X_{1}+\cdots+X_{n}\to Y\) for the morphism induced by the universal property of the sum; given \(g_{k}\colon {X\to Y_{k}}\), we write \((g_{1},\dots,g_{n})\colon X\to Y_{1}\times\cdots\times Y_{n}\) for the morphism induced by the universal property of the product. The sum of morphisms \(f_{k}\colon X_{k}\to Y_{k}\) is written\footnote{We avoid the more common notation \(f_{1}+f_{2}\), because from Section~\ref{sec:square-ring-of-V} on some specific hom-sets of our categories carry a group structure, so that for parallel morphisms \(f\) and \(g\) the symbol \(f+g\) always denotes their sum in this group. The two operations must not be confused: for a finitely generated free algebra \(X\) with coproduct injections \(i_{1}\), \(i_{2}\colon X\to X+X\), we shall find that \(i_{1}+i_{2}\) is the comultiplication \(\mu_{X}\colon {X\to X+X}\) of~\eqref{eq:mu-is-sum}, whereas \(i_{1}\cosum i_{2}\colon X+X\to(X+X)+(X+X)\).}
\[
  f_{1}\cosum\cdots\cosum f_{n}\colon X_{1}+\cdots+X_{n}\to Y_{1}+\cdots+Y_{n}\text{.}
\]
We write \(i_{k}\colon X_{k}\to X_{1}+\cdots+X_{n}\) for the coproduct injections,
\[
  \fold{X}=\langle 1_{X},1_{X}\rangle\colon {X+X\to X}
  \quad\text{and}\quad
  \diag{X}=(1_{X},1_{X})\colon X\to X\times X
\]
for the codiagonal and diagonal, \(\tw{X}{Y}=\langle i_{2},i_{1}\rangle\colon X+Y\to Y+X\) for the switch of the sum, and \(\twp{X}{Y}\colon X\times Y\to Y\times X\) for the switch of the product.

For \(1\leq k\leq n\), let
\[
  r^{(k)}\colon X_{1}+\cdots+X_{n}\to X_{1}+\cdots+\widehat{X_{k}}+\cdots+X_{n}
\]
be the morphism which is the identity on the summands \(X_{j}\) with \(j\neq k\) and zero on \(X_{k}\); here \(\widehat{X_{k}}\) means that the summand \(X_{k}\) is omitted. In the binary case we write \(r_{1}=\langle 1_{X_{1}},0\rangle\colon X_{1}+X_{2}\to X_{1}\) and \(r_{2}=\langle 0,1_{X_{2}}\rangle\colon X_{1}+X_{2}\to X_{2}\), and
\[
  \cm{X_{1}}{X_{2}}\DefEq(r_{1},r_{2})\colon X_{1}+X_{2}\to X_{1}\times X_{2}
\]
for the \defn{comparison morphism}. It is natural in \(X_{1}\) and \(X_{2}\): for \(f\colon X_{1}\to Y_{1}\) and \(g\colon X_{2}\to Y_{2}\),
\begin{equation}\label{eq:r-natural}
  (f\times g)\circ\cm{X_{1}}{X_{2}}=\cm{Y_{1}}{Y_{2}}\circ(f\cosum g)\text{,}
\end{equation}
as one sees by precomposing with the coproduct injections. When \(\C\) is semi-abelian, \(\cm{X_{1}}{X_{2}}\) is a regular epimorphism: \(\C\) is unital, so \((1_{X_{1}},0)\colon {X_{1}\to X_{1}\times X_{2}}\) and \((0,1_{X_{2}})\colon X_{2}\to X_{1}\times X_{2}\) are jointly strongly epimorphic~\cite{Borceux-Bourn}, and strong epimorphisms are regular in a regular category.

\subsection{Cross-effects}\label{sec:cross-effects}

\begin{definition}\label{def:cross-effect}
  Let \(F\colon \C\to\D\) be a reduced functor, \(n\geq 1\), and \(X_{1}\), \dots, \(X_{n}\) objects of \(\C\). Consider the morphism
  \[
    \cmF{F}{n}\DefEq\bigl(F(r^{(1)}),\dots,F(r^{(n)})\bigr)\colon F(X_{1}+\cdots+X_{n})\to\prod_{k=1}^{n}F(X_{1}+\cdots+\widehat{X_{k}}+\cdots+X_{n})\text{,}
  \]
  where for \(n=2\) we order the factors as \(F(X_{1})\times F(X_{2})\), so that \(\cmF{F}{2}=(F(r_{1}),F(r_{2}))\). The \defn{\(n\)-th cross-effect} of \(F\) at \((X_{1},\dots,X_{n})\) is the kernel
  \[
    \cre_{n}F(X_{1},\dots,X_{n})=F(X_{1}|\cdots|X_{n})\DefEq\Ker(\cmF{F}{n})\text{,}
  \]
  with kernel morphism \(\iota^{F}_{X_{1},\dots,X_{n}}\colon F(X_{1}|\cdots|X_{n})\to F(X_{1}+\cdots+X_{n})\). Given morphisms \(f_{k}\colon X_{k}\to Y_{k}\), the morphism \(F(f_{1}\cosum\cdots\cosum f_{n})\) commutes with the \(F(r^{(k)})\) and hence restricts to a morphism
  \[
    F(f_{1}|\cdots|f_{n})\colon F(X_{1}|\cdots|X_{n})\to F(Y_{1}|\cdots|Y_{n})\text{.}
  \]
  This makes \(\cre_{n}F\) a functor \(\C^{\times n}\to\D\).
\end{definition}

Since \(F\) is reduced, \(\cre_{1}F(X)=\Ker(F(X)\to F(0)=0)=F(X)\). The case \(F=1_{\C}\) and \(n=2\) gives the \defn{cosmash product}~\cite{CJ03}
\[
  X\cosmash Y\DefEq 1_{\C}(X|Y)=\Ker(\cm{X}{Y}\colon X+Y\to X\times Y)
\]
with kernel morphism \(\iota_{X,Y}\colon X\cosmash Y\to X+Y\); for \(f\colon X\to X'\) and \(g\colon Y\to Y'\) we write \(f\cosmash g\DefEq 1_{\C}(f|g)\colon X\cosmash Y\to X'\cosmash Y'\) for the restriction of \(f\cosum g\). Likewise \(X\cosmash Y\cosmash Z\DefEq 1_{\C}(X|Y|Z)\) is the \defn{ternary cosmash product}, with kernel morphism \(\iota_{X,Y,Z}\colon X\cosmash Y\cosmash Z\to X+Y+Z\).

\begin{remark}\label{rem:cross-effect-properties}
  The definition above is the one of~\cite[Remark~2.8]{HVdL}. It is equivalent to the inductive definition of~\cite[Definition~1.2]{Hartl-Vespa}: \(\cre_{1}F=F\), \(\cre_{2}F(X_{1},X_{2})=\Ker(\cmF{F}{2})\), and
  \begin{equation}\label{eq:cross-effect-inductive}
    \cre_{n}F(X_{1},\dots,X_{n})\cong\cre_{2}\bigl(\cre_{n-1}F(-,X_{3},\dots,X_{n})\bigr)(X_{1},X_{2})
  \end{equation}
  naturally, for \(n\geq 3\). We note two consequences.
  \begin{enumerate}
    \item \(\cre_{n}F\) is \defn{symmetric}: permuting the \(X_{k}\) permutes the terms of the sum \(X_{1}+\cdots+X_{n}\) and the factors of the product in Definition~\ref{def:cross-effect} accordingly, so \(\cre_{n}F(X_{1},\dots,X_{n})\cong\cre_{n}F(X_{\sigma(1)},\dots,X_{\sigma(n)})\) for every permutation \(\sigma\), naturally in the \(X_{k}\).
    \item \(\cre_{n}F\) is \defn{multireduced}: \(\cre_{n}F(X_{1},\dots,X_{n})=0\) as soon as one \(X_{k}\) is \(0\). For \(n=2\) and \(X_{1}=0\), the morphism \(r_{2}\colon 0+X_{2}\to X_{2}\) is an isomorphism, so that \(\cmF{F}{2}\) is a monomorphism and its kernel vanishes; the case \(X_{2}=0\) is symmetric, and the general case follows from~\eqref{eq:cross-effect-inductive} and symmetry.
  \end{enumerate}
\end{remark}

\begin{remark}\label{rem:r2-regular-epi}
  Suppose that \(\D\) is semi-abelian. Since \(F\) is reduced, \(F(r_{1})\circ F(i_{2})=F(0)=0\), so that
  \[
    \cmF{F}{2}\circ\langle F(i_{1}),F(i_{2})\rangle=\cm{F(X_{1})}{F(X_{2})}\colon F(X_{1})+F(X_{2})\to F(X_{1})\times F(X_{2})\text{.}
  \]
  The right-hand side is a regular epimorphism (Subsection~\ref{sec:conventions}), hence so is \(\cmF{F}{2}\): a morphism through which a regular epimorphism factors is a strong epimorphism, hence regular. Thus
  \[
    \begin{tikzcd}
      0 \arrow[r] & F(X_{1}|X_{2}) \arrow[r, "\iota^{F}_{X_{1},X_{2}}"] & F(X_{1}+X_{2}) \arrow[r, "\cmF{F}{2}"] & F(X_{1})\times F(X_{2}) \arrow[r] & 0
    \end{tikzcd}
  \]
  is a short exact sequence. For \(n\geq 3\), the morphism \(\cmF{F}{n}\) need not be a regular epimorphism.
\end{remark}

\subsection{Polynomial functors}\label{sec:polynomial-functors}

\begin{definition}\label{def:polynomial}
  A reduced functor \(F\colon \C\to\D\) is \defn{polynomial of degree \(\leq n\)} if its \((n+1)\)-th cross-effect vanishes: \(\cre_{n+1}F=0\). We write \(\Func_{\leq n}(\C,\D)\) for the full subcategory of \(\Func_{*}(\C,\D)\) determined by these functors. A functor is \defn{linear} if it is polynomial of degree \(\leq 1\), and \defn{quadratic} if it is polynomial of degree \(\leq 2\).
\end{definition}

If \(F\) is polynomial of degree \(\leq n\), then it is polynomial of degree \(\leq k\) for all \(k\geq n\): by~\eqref{eq:cross-effect-inductive},
\[
  \cre_{n+2}F(X_{1},\dots,X_{n+2})\cong\cre_{2}\bigl(\cre_{n+1}F(-,X_{3},\dots,X_{n+2})\bigr)(X_{1},X_{2})=0\text{.}
\]

\begin{proposition}\label{prop:quadratic-iff-bilinear}
  A reduced functor \(F\colon \C\to\D\) is quadratic if and only if its second cross-effect is \defn{bilinear}: for every object \(X\) of \(\C\), the functors \(F(X|-)\) and \(F(-|X)\colon \C\to\D\) are linear.
\end{proposition}

\begin{proof}
  By~\eqref{eq:cross-effect-inductive}, \(\cre_{3}F(X_{1},X_{2},X_{3})\cong\cre_{2}(F(-|X_{3}))(X_{1},X_{2})\). Hence \(\cre_{3}F=0\) if and only if each \(F(-|X_{3})\) is linear, which by symmetry of \(\cre_{2}F\) is equivalent to each \(F(X_{3}|-)\) being linear.
\end{proof}

The next lemma replaces~\cite[Lemma~3.6]{Hartl-Vespa} in our setting.

\begin{lemma}\label{lem:linear-into-additive}
  Let \(\A\) be an additive category with kernels and \(F\colon \C\to\A\) a reduced linear functor. For all \(X\), \(Y\) in \(\C\), the morphism \(\cmF{F}{2}\colon F(X+Y)\to F(X)\times F(Y)\) is an isomorphism and
  \begin{equation}\label{eq:linear-biproduct}
    1_{F(X+Y)}=F(i_{1})\circ F(r_{1})+F(i_{2})\circ F(r_{2})\text{.}
  \end{equation}
\end{lemma}

\begin{proof}
  The kernel \(F(X|Y)\) of \(\cmF{F}{2}\) vanishes, so \(\cmF{F}{2}\) is a monomorphism, \(\A\) being additive. On the other hand,
  \[
    c\DefEq\cmF{F}{2}\circ\langle F(i_{1}),F(i_{2})\rangle\colon F(X)+F(Y)\to F(X)\times F(Y)
  \]
  is the canonical comparison morphism from the sum to the product (see Remark~\ref{rem:r2-regular-epi}), which in an additive category is an isomorphism. Hence \(\cmF{F}{2}\) is a split epimorphism, so it is an isomorphism, with inverse \(s\DefEq\langle F(i_{1}),F(i_{2})\rangle\circ c^{-1}\). Writing \(j_{1}\), \(j_{2}\) for the coproduct injections of \(F(X)+F(Y)\) and \(p_{1}\), \(p_{2}\) for the projections of \(F(X)\times F(Y)\), we have \(c^{-1}=j_{1}\circ p_{1}+j_{2}\circ p_{2}\), since \(c\circ j_{1}=(1_{F(X)},0)\) and \(c\circ j_{2}=(0,1_{F(Y)})\). Therefore
  \begin{align*}
    1_{F(X+Y)}=s\circ\cmF{F}{2} & =\langle F(i_{1}),F(i_{2})\rangle\circ(j_{1}\circ p_{1}+j_{2}\circ p_{2})\circ(F(r_{1}),F(r_{2})) \\
                                & =F(i_{1})\circ F(r_{1})+F(i_{2})\circ F(r_{2})
  \end{align*}
  by bilinearity of composition.
\end{proof}

By induction on \(n\), the same holds for \(n\)-fold sums: if \(F\) is linear, then we have \(F(X_{1}+\cdots+X_{n})\cong F(X_{1})\times\cdots\times F(X_{n})\) and \(1_{F(X_{1}+\cdots+X_{n})}=\sum_{k=1}^{n}F(i_{k})\circ F(r'_{k})\), where \(r'_{k}\colon X_{1}+\cdots+X_{n}\to X_{k}\) is the retraction which is the identity on \(X_{k}\) and zero on the other summands.

We shall apply this lemma to functors \(F\colon \C\to\D\) with \(\D\) semi-abelian which take their values in a full subcategory \(\A\) of \(\D\) that is closed under kernels and finite products, typically the abelian core \(\Ab(\D)\). The cross-effects of \(F\), computed in \(\D\), then lie in \(\A\) and coincide with the cross-effects of \(F\) viewed as a functor \(\C\to\A\).

\subsection{Folding, Taylorisation and commutators}\label{sec:folding}
From now on, \(\C\) and \(\D\) are semi-abelian categories.

\begin{definition}\label{def:folding}
  Let \(F\colon \C\to\D\) be reduced and \(n\geq 1\). Write \(\Delta^{n}\colon \C\to\C^{\times n}\) for the diagonal functor and \(\fold{X}^{n}=\langle 1_{X},\dots,1_{X}\rangle\colon X+\cdots+X\to X\) for the \(n\)-fold codiagonal. The \defn{folding transformation} \(S^{F}_{n}\colon \cre_{n}F\circ\Delta^{n}\Rightarrow F\) has components
  \[
    (S^{F}_{n})_{X}\DefEq F(\fold{X}^{n})\circ\iota^{F}_{X,\dots,X}\colon F(X|\cdots|X)\to F(X)\text{.}
  \]
  It is natural because \(f\circ\fold{X}^{n}=\fold{Y}^{n}\circ(f\cosum\cdots\cosum f)\) for all \(f\colon X\to Y\).
\end{definition}

For \(F=1_{\C}\) and \(n=2\) we write
\[
  \cc{X}\DefEq(S^{1_{\C}}_{2})_{X}=\fold{X}\circ\iota_{X,X}\colon X\cosmash X\to X\text{,}
\]
and naturality amounts to
\begin{equation}\label{eq:c-natural}
  f\circ\cc{X}=\cc{Y}\circ(f\cosmash f)\qquad\text{for all \(f\colon X\to Y\).}
\end{equation}

\begin{definition}\label{def:taylorisation}
  The \defn{\(n\)-th Taylorisation} of a reduced functor \(F\colon \C\to\D\) is the functor \(T_{n}F\colon \C\to\D\) defined by \(T_{n}F(X)\DefEq\Coker\bigl((S^{F}_{n+1})_{X}\bigr)\). Hartl and Vespa~\cite[Definition~1.9, Proposition~1.10]{Hartl-Vespa} show that, when \(\D\) is the category of groups, \(T_{n}F\) is polynomial of degree \(\leq n\) and \(T_{n}\) is left adjoint to the inclusion \(\Func_{\leq n}(\C,\D)\to\Func_{*}(\C,\D)\). We only need the case \(F=1_{\C}\), where these facts are classical; see Subsection~\ref{sec:nilpotent}.
\end{definition}

\begin{definition}\label{def:commutator}
  Let \(F\colon \C\to\D\) be a reduced functor and let \(x_{k}\colon X_{k}\rightarrowtail X\), for \(1\leq k\leq n\), be subobjects of an object \(X\) of \(\C\). The \defn{\(n\)-fold commutator of \(X_{1}\), \dots, \(X_{n}\) relative to \(F\)} is the regular image
  \[
    [X_{1},\dots,X_{n}]_{F}\DefEq\Im\bigl((S^{F}_{n})_{X}\circ F(x_{1}|\cdots|x_{n})\colon F(X_{1}|\cdots|X_{n})\to F(X)\bigr)\text{,}
  \]
  a subobject of \(F(X)\). We write \(\gamma^{F}_{n}(X)\DefEq[X,\dots,X]_{F}=\Im\bigl((S^{F}_{n})_{X}\bigr)\) for the commutator of \(n\) copies of \(X\) itself. In particular, \([X]_{F}=F(X)\) since \(\cre_{1}F=F\).
\end{definition}

For \(F=1_{\C}\) we drop the index \(F\) and \([X_{1},X_{2}]\) is the \defn{Higgins commutator}~\cite{MM} of the subobjects \(X_{1}\), \(X_{2}\) of \(X\), the image of \(\langle x_{1},x_{2}\rangle\circ\iota_{X_{1},X_{2}}\colon X_{1}\cosmash X_{2}\to X\), and \([X_{1},\dots,X_{n}]\) is the \defn{\(n\)-fold Higgins commutator} of~\cite[Definition~2.14]{HVdL}. Commutators relative to a functor will be studied in the forthcoming~\cite{Hartl-Linearization}.

\begin{proposition}\label{prop:polynomial-iff-commutators}
  A reduced functor \(F\colon \C\to\D\) is polynomial of degree \(\leq n\) if and only if \(\gamma^{F}_{n+1}(X)=0\) for every object \(X\) of \(\C\).
\end{proposition}

\begin{proof}
  If \(\cre_{n+1}F=0\), then \((S^{F}_{n+1})_{X}\) has domain \(0\) and its image vanishes. Conversely, suppose that \(\gamma^{F}_{n+1}(X)=0\), i.e., \((S^{F}_{n+1})_{X}=0\), for all \(X\). Given objects \(X_{1}\), \dots, \(X_{n+1}\), put \(X=X_{1}+\cdots+X_{n+1}\) and consider the coproduct injections \(i_{k}\colon X_{k}\to X\). Since \(\fold{X}^{n+1}\circ(i_{1}\cosum\cdots\cosum i_{n+1})=\langle i_{1},\dots,i_{n+1}\rangle=1_{X}\), we find
  \begin{align*}
    0 & =(S^{F}_{n+1})_{X}\circ F(i_{1}|\cdots|i_{n+1})=F(\fold{X}^{n+1})\circ F(i_{1}\cosum\cdots\cosum i_{n+1})\circ\iota^{F}_{X_{1},\dots,X_{n+1}} \\
      & =\iota^{F}_{X_{1},\dots,X_{n+1}}\text{,}
  \end{align*}
  a monomorphism. Hence \(F(X_{1}|\cdots|X_{n+1})=0\).
\end{proof}

\subsection{Nilpotent objects and categories}\label{sec:nilpotent}

\begin{definition}\label{def:nilpotent}
  An object \(X\) of \(\C\) is \defn{\(n\)-step nilpotent} if the \((n+1)\)-fold Higgins commutator \(\gamma_{n+1}(X)\DefEq[X,\dots,X]=\Im(\fold{X}^{n+1}\circ\iota_{X,\dots,X})\) vanishes. For \(n=1\) these are the \defn{abelian} objects of \(\C\). We write \(\Nil_{n}(\C)\) for the full subcategory of \(\C\) determined by the \(n\)-step nilpotent objects, and \(\Ab(\C)=\Nil_{1}(\C)\) for the \defn{abelian core} of \(\C\).

  The category \(\C\) itself is \defn{\(n\)-step nilpotent} if its identity functor \(1_{\C}\) is polynomial of degree \(\leq n\), i.e., if all \((n+1)\)-fold cosmash products in \(\C\) vanish. By Proposition~\ref{prop:polynomial-iff-commutators}, this happens if and only if every object of \(\C\) is \(n\)-step nilpotent.
\end{definition}

This illustrates Hartl's slogan that \emph{polynomiality subsumes nilpotency}. The subcategory \(\Nil_{n}(\C)\) is a Birkhoff subcategory~\cite{Janelidze-Kelly} of \(\C\)~\cite{HVdL}: it is reflective and closed under subobjects and regular quotients, and its reflector \(\nil_{n}\colon \C\to\Nil_{n}(\C)\) sends \(X\) to the quotient \(\nil_{n}(X)=X/\gamma_{n+1}(X)\), which is the cokernel of \((S^{1_{\C}}_{n+1})_{X}\)---so \(\nil_{n}=T_{n}1_{\C}\) is the \(n\)-th Taylorisation of the identity functor. For \(n=1\) this is the \defn{abelianisation} functor \(\ab\colon \C\to\Ab(\C)\), whose unit we write \(\ab_{X}\colon X\to\ab(X)\); by construction, \(\ab_{X}\) is the cokernel \(\coker(\cc{X})\) of \(\cc{X}\colon X\cosmash X\to X\), so that \(\ab_{X}\circ\cc{X}=0\). We thus obtain the \defn{nilpotency tower}
\[
  \Ab(\C)=\Nil_{1}(\C)\subseteq\Nil_{2}(\C)\subseteq\cdots\subseteq\Nil_{n}(\C)\subseteq\cdots\subseteq\C\text{.}
\]
Nilpotency in this categorical sense goes back to Huq~\cite{Huq}; for the semi-abelian and the exact Mal'tsev settings, see~\cite{Bourn-Gran, BeBou}, where the nilpotency tower is studied by means of central extensions and functor calculus. For the commutator-theoretic approach to nilpotency in universal algebra, see~\cite{Sm, Freese-McKenzie}, and~\cite{Higgins} for the higher commutators, which play a role comparable to that of the higher cosmash products. The abelian core \(\Ab(\C)\) is an abelian category, in which finite sums are biproducts, so that it is \(1\)-step nilpotent; the \defn{two-nilpotent core} \(\Nil_{2}(\C)\) is a two-step nilpotent semi-abelian category~\cite[Subsection~2.3]{DHVdL}.

\begin{example}\label{ex:nilpotent-groups}
  The variety \(\Nil_{2}(\Grp)\) of groups of nilpotency class at most \(2\) is a two-step nilpotent category. This can also be seen from the normal form of its binary coproducts~\cite[Proposition~1.2]{Jibladze-Pirashvili}: for \(X\), \(Y\) in \(\Nil_{2}(\Grp)\),
  \[
    X+Y\cong X\times Y\times(\ab(X)\otimes\ab(Y))
  \]
  as sets, where \(\otimes\) is the tensor product of abelian groups. The kernel of \(\cm{X}{Y}\) is the summand \(\ab(X)\otimes\ab(Y)\), and the elements of \(X+Y+Z\) which vanish under all three retractions \(r^{(k)}\) are products of triple commutators \([[x,y],z]\), which are trivial. This coproduct formula is the prototype of the description of coproducts of modules over a square ring in Proposition~\ref{prop:coproduct}.
\end{example}

\begin{proposition}\label{prop:polynomial-values-nilpotent}
  If \(F\colon \C\to\D\) is a reduced polynomial functor of degree \(\leq n\), then \(F\) takes its values in \(\Nil_{n}(\D)\).
\end{proposition}

\begin{proof}
  Let \(X\) be an object of \(\C\) and write \(Q\DefEq F(X)\cosmash\cdots\cosmash F(X)\) for the \((n+1)\)-fold cosmash product in \(\D\), so that \(\gamma_{n+1}(F(X))\) is the image of
  \[
    \fold{F(X)}^{n+1}\circ\iota_{F(X),\dots,F(X)}\colon Q\to F(X)\text{.}
  \]
  Consider \(j\DefEq\langle F(i_{1}),\dots,F(i_{n+1})\rangle\colon F(X)+\cdots+F(X)\to F(X+\cdots+X)\). By naturality of the \(r^{(k)}\) and since \(F\) is reduced, \(\cmF{F}{n+1}\circ j\) factors through \(\hat{r}^{1_{\D}}_{n+1}\), so \(j\) restricts to a morphism \(j'\colon Q\to F(X|\cdots|X)\) with \(\iota^{F}_{X,\dots,X}\circ j'=j\circ\iota_{F(X),\dots,F(X)}\). Moreover, \(F(\fold{X}^{n+1})\circ j=\langle F(\fold{X}^{n+1}\circ i_{k})\rangle_{k}=\langle 1_{F(X)},\dots,1_{F(X)}\rangle=\fold{F(X)}^{n+1}\). Hence
  \[
    \fold{F(X)}^{n+1}\circ\iota_{F(X),\dots,F(X)}=F(\fold{X}^{n+1})\circ j\circ\iota_{F(X),\dots,F(X)}=(S^{F}_{n+1})_{X}\circ j'\text{,}
  \]
  whose image is contained in \(\gamma^{F}_{n+1}(X)\), which vanishes by Proposition~\ref{prop:polynomial-iff-commutators}. So \(\gamma_{n+1}(F(X))=0\).
\end{proof}

\begin{proposition}\label{prop:cosmash-abelian}
  If \(\C\) is a two-step nilpotent category, then for all objects \(X\), \(Y\) of \(\C\) the cosmash product \(X\cosmash Y\) is an abelian object, and the functors \(X\cosmash(-)\) and \((-)\cosmash Y\colon \C\to\Ab(\C)\) are linear.
\end{proposition}

\begin{proof}
  By symmetry of the cosmash product, it suffices to consider the first case. The identity functor is quadratic, so \(X\cosmash(-)=1_{\C}(X|-)\) is linear by Proposition~\ref{prop:quadratic-iff-bilinear}, and takes its values in \(\Nil_{1}(\C)=\Ab(\C)\) by Proposition~\ref{prop:polynomial-values-nilpotent}.
\end{proof}

The next proposition is where we use the right exactness results of~\cite{DHVdL}. Recall that a functor which preserves coequalisers of reflexive graphs preserves regular epimorphisms, since every regular epimorphism is the coequaliser of its kernel pair.

\begin{proposition}\label{prop:F-nil}
  Let \(F\colon \C\to\D\) be a reduced polynomial functor of degree \(\leq n\) which preserves coequalisers of reflexive graphs. Then for every object \(X\) of \(\C\), the morphism \(F((\nil_{n})_{X})\colon F(X)\to F(\nil_{n}(X))\) is an isomorphism. In other words, \(F\cong F\circ\nil_{n}\).
\end{proposition}

\begin{proof}
  Write \(k\colon K\DefEq\gamma_{n+1}(X)\rightarrowtail X\), so that \((\nil_{n})_{X}\colon X\to\nil_{n}(X)\) is the cokernel of \(k\). By~\cite[Proposition~9.8]{DHVdL}, \(F\) sends this cokernel to a right exact sequence
  \[
    F(K|X)\rtimes F(K)\xrightarrow{\;\langle (S^{F}_{2})_{X}\circ F(k|1_{X}),\,F(k)\rangle\;}F(X)\xrightarrow{\;F((\nil_{n})_{X})\;}F(\nil_{n}(X))\to 0\text{,}
  \]
  so it suffices to show that both \(F(k)\) and \((S^{F}_{2})_{X}\circ F(k|1_{X})\) are zero. Let \(Q\DefEq X\cosmash\cdots\cosmash X\) be the \((n+1)\)-fold cosmash product and \(q\colon Q\to K\) the regular epimorphism with \(k\circ q=\fold{X}^{n+1}\circ\iota_{X,\dots,X}\).

  Since \(F\) preserves regular epimorphisms,
  \[
    \Im(F(k))=\Im(F(k)\circ F(q))=\Im\bigl(F(\fold{X}^{n+1})\circ F(\iota_{X,\dots,X})\bigr)\text{.}
  \]
  Now \(F(r^{(j)})\circ F(\iota_{X,\dots,X})=F(r^{(j)}\circ\iota_{X,\dots,X})=0\) for all \(j\), so \(F(\iota_{X,\dots,X})\) factors through \(\iota^{F}_{X,\dots,X}\); hence \(\Im(F(k))\subseteq\Im((S^{F}_{n+1})_{X})=\gamma^{F}_{n+1}(X)=0\) by Proposition~\ref{prop:polynomial-iff-commutators}.

  By~\cite[Theorem~9.5]{DHVdL}, the functor \(F(-|X)\) also preserves coequalisers of reflexive graphs, so \(F(q|1_{X})\colon F(Q|X)\to F(K|X)\) is a regular epimorphism and it suffices to show that \(u\DefEq(S^{F}_{2})_{X}\circ F(k\circ q|1_{X})\colon F(Q|X)\to F(X)\) vanishes. We have
  \[
    u=F(\fold{X})\circ F(k\circ q\cosum 1_{X})\circ\iota^{F}_{Q,X}=F(\fold{X}^{n+2})\circ F(\iota_{X,\dots,X}\cosum 1_{X})\circ\iota^{F}_{Q,X}\text{,}
  \]
  because \(\fold{X}\circ(k\circ q\cosum 1_{X})=\langle\fold{X}^{n+1}\circ\iota_{X,\dots,X},1_{X}\rangle=\fold{X}^{n+2}\circ(\iota_{X,\dots,X}\cosum 1_{X})\). The morphism \(F(\iota_{X,\dots,X}\cosum 1_{X})\circ\iota^{F}_{Q,X}\) factors through \(\iota^{F}_{X,\dots,X}\), the kernel of \(\cmF{F}{n+2}\): for \(j\leq n+1\) we have \(r^{(j)}\circ(\iota_{X,\dots,X}\cosum 1_{X})=(r^{(j)}\circ\iota_{X,\dots,X})\cosum 1_{X}=0\cosum 1_{X}=i_{2}\circ r_{2}\), and \(F(r_{2})\circ\iota^{F}_{Q,X}=0\); for \(j=n+2\) we have \(r^{(n+2)}\circ(\iota_{X,\dots,X}\cosum 1_{X})=\iota_{X,\dots,X}\circ r_{1}\) and \(F(r_{1})\circ\iota^{F}_{Q,X}=0\). Hence \(\Im(u)\subseteq\gamma^{F}_{n+2}(X)\), which vanishes since \(F\) is also polynomial of degree \(\leq n+1\).
\end{proof}

\begin{corollary}\label{cor:cosmash-abelianisation}
  Let \(\C\) be a two-step nilpotent category. For all objects \(X\), \(Y\) of \(\C\), the morphisms
  \[
    1_{X}\cosmash\ab_{Y}\colon X\cosmash Y\to X\cosmash\ab(Y)\qquad\text{and}\qquad\ab_{X}\cosmash 1_{Y}\colon X\cosmash Y\to\ab(X)\cosmash Y
  \]
  are isomorphisms. Consequently \(1_{X}\cosmash\cc{Y}=0\colon X\cosmash(Y\cosmash Y)\to X\cosmash Y\) and \(\cc{X}\cosmash 1_{Y}=0\).
\end{corollary}

\begin{proof}
  The functor \(X\cosmash(-)\) is linear by Proposition~\ref{prop:cosmash-abelian} and preserves coequalisers of reflexive graphs by~\cite[Corollary~9.6]{DHVdL}, so Proposition~\ref{prop:F-nil} with \(n=1\) applies; the second morphism is treated symmetrically. Since \((1_{X}\cosmash\ab_{Y})\circ(1_{X}\cosmash\cc{Y})=1_{X}\cosmash(\ab_{Y}\circ\cc{Y})=0\) and \(1_{X}\cosmash\ab_{Y}\) is a monomorphism, \(1_{X}\cosmash\cc{Y}=0\).
\end{proof}

\section{Square rings and their modules}\label{sec:square-rings}

Square rings and their modules were introduced by Baues, Hartl and Pirashvili~\cite{BHP} in the context of metastable homotopy theory; they were later used by Baues and Muro in the study of secondary homotopy operations~\cite{Baues-Muro}, and quadratic maps between modules over a square ring are studied in~\cite{Gaudier-Hartl:Morphismes}; see also the unpublished notes~\cite{Pirashvili-Square-Rings, Hartl-WIP}. The motivating example is the endomorphism square ring of the suspended projective plane \(\Sigma\RP^{2}\), which is a square ring by a result of Barratt~\cite{Barratt-TrackGroups} (see~\cite{Baues-Iwase} for square rings associated with elements of homotopy groups of spheres in general); by~\cite[Theorem~8.1]{BHP}, the category of free modules over it is equivalent to the homotopy category of Moore spaces \(M(V,2)\), where \(V\) is a \(\Z/2\Z\)-vector space. Although this example may seem remote from our subject, its categorical structure is what underlies the Main Theorem. In this section we recall the definitions and prove, in full, the elementary facts we need.

\subsection{Square rings}\label{sec:square-rings-def}

\begin{convention}\label{conv:nearring}
  A \defn{left nearring} is a set \(R\) with a group structure \((R,+,0)\), not necessarily abelian, and a monoid structure \((R,\cdot,1)\) such that multiplication is left distributive over addition: \(r\cdot(s+t)=r\cdot s+r\cdot t\) for all \(r\), \(s\), \(t\in R\). It follows that \(r\cdot 0=0\) and \(r\cdot(-s)=-r\cdot s\). In~\cite{BHP} these structures are called \emph{right} nearrings, because the sum being distributed appears on the right; we prefer to name the side on which the multiplication distributes. A \defn{morphism of nearrings} is a map which preserves \(+\), \(0\), \(\cdot\) and \(1\).
\end{convention}

Let \(G\) and \(G'\) be groups. Following~\cite{BHP}, a map \(f\colon G\to G'\) is \defn{linear} if \(f(g+g')=f(g)+f(g')\) for all \(g\), \(g'\in G\), i.e., if it is a group homomorphism; a map of several variables is \defn{multilinear} if it is linear in each variable separately.

\begin{definition}[{\cite[Section~7]{BHP}}]\label{def:square-ring}
  A (left) \defn{square ring} \(\R\) consists of the following data:
  \begin{itemize}
    \item a left nearring \(\Re\);
    \item an abelian group \(\Ree\);
    \item a quadrilinear map
          \begin{equation}\label{eq:quadrilinear}
            \Re\times\Re\times\Ree\times\Re\to\Ree\colon (r,s,x,t)\mapsto(r,s)\cdot x\cdot t
          \end{equation}
          which makes \(\Ree\) a left module over the multiplicative monoid \(\Re\times\Re\) and a right module over the multiplicative monoid \(\Re\), the two actions commuting with each other;
    \item a group homomorphism \(T\colon \Ree\to\Ree\) with \(T\circ T=1_{\Ree}\);
    \item a map \(H\colon \Re\to\Ree\), not necessarily a homomorphism, and a group homomorphism \(P\colon \Ree\to\Re\);
  \end{itemize}
  such that, for all \(r\), \(s\), \(t\in\Re\) and \(x\), \(y\in\Ree\), writing \(2\DefEq 1+1\in\Re\):
  \begin{sqaxioms}
    \item\label{SR:P-binatural} \(P((r,r)\cdot x\cdot s)=r\cdot P(x)\cdot s\);
    \item\label{SR:T} \(T(x)=H(P(x))-x\);
    \item\label{SR:PT} \(P(T(x))=P(x)\);
    \item\label{SR:P-kills} \((P(x),r)\cdot y=0=(r,P(x))\cdot y\) and \(y\cdot P(x)=0\);
    \item\label{SR:H-sum} \(H(r+s)=H(r)+H(s)+(s,r)\cdot H(2)\);
    \item\label{SR:H-product} \(H(r\cdot s)=(r,r)\cdot H(s)+H(r)\cdot s\);
    \item\label{SR:right-distributivity} \((r+s)\cdot t=r\cdot t+s\cdot t+P((r,s)\cdot H(t))\);
    \item\label{SR:T-action} \(T((r,s)\cdot x\cdot t)=(s,r)\cdot T(x)\cdot t\).
  \end{sqaxioms}
  We write \(\R=\sqring{\Re}{\Ree}{H}{P}\), or simply \(\R\). Right square rings are defined dually; we only consider left square rings.
\end{definition}

Some immediate consequences: \(H(1)=0\) by~\ref{SR:H-product} with \(r=s=1\); \(0\cdot t=0\) by~\ref{SR:right-distributivity} with \(r=s=0\) and quadrilinearity; \((2,2)\cdot x=4x\) by quadrilinearity; and \(P\circ H\circ P=P+P\), since by~\ref{SR:PT} and~\ref{SR:T} we have \(P(x)=P(T(x))=P(H(P(x)))-P(x)\). In~\cite{BHP} the latter identity is part of the definition. By~\ref{SR:T}, the involution \(T\) is determined by \(H\) and \(P\); it is nevertheless part of the data, as in~\cite{BHP}, because the axioms are most naturally stated in terms of it.

\begin{remark}
  Square rings are the one-object case of \emph{square ringoids}, the quadratic analogues of ringoids (preadditive categories). For the proof of our Main Theorem, the one-object case suffices; we refer to~\cite{BHP} for the general theory of square ringoids.
\end{remark}

\begin{examples}\label{ex:square-rings}
  \begin{enumerate}
    \item\label{ex:ring-as-square-ring} Any ring \(R\) is a square ring with \(\Re=R\), \(\Ree=0\) and \(H=P=T=0\). Square rings thus generalise rings.
    \item\label{ex:Z-nil} The \defn{nil square ring} \(\Z_{\nil}=\sqring{\Z}{\Z}{H}{0}\) has \(\Re=\Ree=\Z\), the action \((r,s)\cdot x\cdot t=rsxt\), \(P=0\), \(T=-1\) and \(H(r)=\binom{r}{2}=r(r-1)/2\). It is the endomorphism square ring of the free group of nilpotency class \(\leq 2\) on one generator, which is \(\Z\)~\cite[Section~7]{BHP}.
  \end{enumerate}
\end{examples}

\begin{proposition}\label{prop:coker-P}
  Let \(\R\) be a square ring. The image of \(P\colon \Ree\to\Re\) is a central subgroup of \((\Re,+)\), stable under left and right multiplication by elements of \(\Re\). The quotient \(\Rbar\DefEq\Coker(P)=\Re/\Im(P)\) is a ring, with multiplication induced by that of \(\Re\). We write \(\overline{r}\) for the class of \(r\in\Re\) in \(\Rbar\).
\end{proposition}

\begin{proof}
  Let \(a\in\Re\) and \(x\in\Ree\). By left distributivity, \((a+P(x))\cdot 2=a+P(x)+a+P(x)\), while by~\ref{SR:right-distributivity} and~\ref{SR:P-kills}, \((a+P(x))\cdot 2=a\cdot 2+P(x)\cdot 2+P((a,P(x))\cdot H(2))=a+a+P(x)+P(x)\). Cancelling \(a\) on the left and \(P(x)\) on the right gives \(P(x)+a=a+P(x)\), so \(\Im(P)\) is central, hence normal. By~\ref{SR:P-binatural}, \(P(x)\cdot t=P(x\cdot t)\) and \(r\cdot P(x)=P((r,r)\cdot x)\), so \(\Im(P)\) is stable under multiplication on either side. It follows that the multiplication descends to \(\Rbar\): indeed \((r+P(x))\cdot t=r\cdot t+P(x\cdot t)+P((r,P(x))\cdot H(t))=r\cdot t+P(x\cdot t)\) by~\ref{SR:right-distributivity} and~\ref{SR:P-kills}, and \(r\cdot(t+P(x))=r\cdot t+r\cdot P(x)\). The quotient \(\Rbar\) is a monoid under multiplication, and left distributivity is inherited from \(\Re\). Right distributivity holds in \(\Rbar\) by~\ref{SR:right-distributivity}, since the correction term lies in \(\Im(P)\). Finally, addition in \(\Rbar\) is commutative: right distributivity with \(t=2\) gives \(\overline{r}+\overline{s}+\overline{r}+\overline{s}=(\overline{r}+\overline{s})\cdot\overline{2}=\overline{r}\cdot\overline{2}+\overline{s}\cdot\overline{2}=\overline{r}+\overline{r}+\overline{s}+\overline{s}\), and cancelling yields \(\overline{s}+\overline{r}=\overline{r}+\overline{s}\).
\end{proof}

\begin{lemma}\label{lem:actions-factor}
  Let \(\R\) be a square ring. The actions of \(\Re\times\Re\) and of \(\Re\) on \(\Ree\) factor through \(\Rbar\times\Rbar\) and \(\Rbar\): the maps \(r\mapsto(r,s)\cdot x\), \(s\mapsto(r,s)\cdot x\) and \(t\mapsto x\cdot t\) are group homomorphisms \(\Re\to\Ree\) which vanish on \(\Im(P)\). Consequently, \(\Ree\) is an \((\Rbar\otimes\Rbar)\)-\(\Rbar\)-bimodule via \((\overline{r}\otimes\overline{s})\cdot x\cdot\overline{t}\DefEq(r,s)\cdot x\cdot t\), where \(\Rbar\otimes\Rbar\) is the tensor product of abelian groups with its induced ring structure.
\end{lemma}

\begin{proof}
  The maps are additive by quadrilinearity and vanish on \(\Im(P)\) by~\ref{SR:P-kills}. The bimodule structure follows since the actions are monoid actions which commute with each other.
\end{proof}

The involution \(T\) can be incorporated into this picture by means of a skew group ring. For any ring \(R\), the swap \(\sigma\colon R\otimes R\to R\otimes R\colon r\otimes s\mapsto s\otimes r\) is an involutive ring automorphism, so the symmetric group \(\mathfrak{S}_2=\{1,\tau\}\) acts on \(R\otimes R\) by ring automorphisms. The resulting \defn{skew group ring} \((R\otimes R)\rtimes\mathfrak{S}_2\), also called the \defn{wreath product} of \(R\) with \(\mathfrak{S}_2\) and then written \(R\wr\mathfrak{S}_2\), is the ring
\[
  R\wr\mathfrak{S}_2=(R\otimes R)\rtimes\mathfrak{S}_2\DefEq(R\otimes R)\oplus(R\otimes R)\tau
\]
in which \(\tau^2=1\) and \(\tau a=\sigma(a)\tau\) for all \(a\in R\otimes R\). Explicitly, the product of two elements is
\[
  (a_1+a_2\tau)(a_1'+a_2'\tau)=a_1a_1'+a_2\sigma(a_2')+\bigl(a_1a_2'+a_2\sigma(a_1')\bigr)\tau\text{,}
\]
so that, on elementary tensors,
\begin{multline*}
  (r_1\otimes s_1+(r_2\otimes s_2)\tau)(r_1'\otimes s_1'+(r_2'\otimes s_2')\tau)\\*
  =r_1r_1'\otimes s_1s_1'+r_2s_2'\otimes s_2r_2'+(r_1r_2'\otimes s_1s_2'+r_2s_1'\otimes s_2r_1')\tau
\end{multline*}
for \(r_i\), \(r_i'\), \(s_i\), \(s_i'\in R\). A left \((R\otimes R)\rtimes\mathfrak{S}_2\)-module is thus the same thing as a left \((R\otimes R)\)-module \(M\) equipped with an additive involution \(\tau\colon M\to M\) such that \(\tau(a\cdot x)=\sigma(a)\cdot\tau(x)\) for all \(a\in R\otimes R\) and \(x\in M\).

\begin{lemma}\label{lem:wreath}
  Let \(\R\) be a square ring. The bimodule structure of Lemma~\ref{lem:actions-factor} extends to an \(((\Rbar\otimes\Rbar)\rtimes\mathfrak{S}_2)\)-\(\Rbar\)-bimodule structure on \(\Ree\), in which \(\tau\) acts as \(T\):
  \[
    (\overline{r}_1\otimes\overline{s}_1+(\overline{r}_2\otimes\overline{s}_2)\tau)\cdot x\cdot\overline{t}\DefEq(r_1,s_1)\cdot x\cdot t+(r_2,s_2)\cdot T(x)\cdot t\text{.}
  \]
\end{lemma}

\begin{proof}
  By Lemma~\ref{lem:actions-factor}, \(\Ree\) is an \((\Rbar\otimes\Rbar)\)-\(\Rbar\)-bimodule, and \(T\) is an additive involution of \(\Ree\). By the description of left \(((\Rbar\otimes\Rbar)\rtimes\mathfrak{S}_2)\)-modules given above, it remains to check that \(T((\overline{r}\otimes\overline{s})\cdot x\cdot\overline{t})=(\overline{s}\otimes\overline{r})\cdot T(x)\cdot\overline{t}\) for all \(r\), \(s\), \(t\in\Re\) and \(x\in\Ree\): with \(t=1\) this says that \(T\) is compatible with the left action, and with \(r=s=1\) that it commutes with the right action. Both are instances of~\ref{SR:T-action}.
\end{proof}

\begin{definition}\label{def:square-ring-morphism}
  A \defn{morphism of square rings} \(f=(f_{e},f_{ee})\colon \R\to\R'\) consists of a nearring morphism \(f_{e}\colon \Re\to\R'_{e}\) and a group homomorphism \(f_{ee}\colon \Ree\to\R'_{ee}\) such that \(H'\circ f_{e}=f_{ee}\circ H\), \(P'\circ f_{ee}=f_{e}\circ P\) and \(f_{ee}((r,s)\cdot x\cdot t)=(f_{e}(r),f_{e}(s))\cdot f_{ee}(x)\cdot f_{e}(t)\) for all \(r\), \(s\), \(t\in\Re\) and \(x\in\Ree\). Then also \(T'\circ f_{ee}=f_{ee}\circ T\) by~\ref{SR:T}. Square rings and their morphisms form a category \(\SqRing\), composition being componentwise. The square ring \(\Z_{\nil}\) is an \defn{initial object} of \(\SqRing\), see~\cite{BHP}.
\end{definition}

\subsection{Modules over a square ring}\label{sec:modules}

\begin{definition}[{\cite[Definition~7.7]{BHP}}]\label{def:module}
  Let \(\R\) be a square ring. A (right) \defn{\(\R\)-module} is a group \(M\), written additively but not necessarily abelian, equipped with two operations
  \[
    M\times\Re\to M\colon (m,r)\mapsto m\cdot r\qquad\text{and}\qquad M\times M\times\Ree\to M\colon (m,n,x)\mapsto\brk{m}{n}\cdot x
  \]
  such that, for all \(m\), \(n\in M\), \(r\), \(s\in\Re\) and \(x\in\Ree\):
  \begin{modaxioms}
    \item\label{M:action} \(m\cdot 1=m\), \((m\cdot r)\cdot s=m\cdot(r\cdot s)\) and \(m\cdot(r+s)=m\cdot r+m\cdot s\);
    \item\label{M:sum} \((m+n)\cdot r=m\cdot r+n\cdot r+\brk{m}{n}\cdot H(r)\);
    \item\label{M:P} \(m\cdot P(x)=\brk{m}{m}\cdot x\);
    \item\label{M:T} \(\brk{m}{n}\cdot T(x)=\brk{n}{m}\cdot x\);
    \item\label{M:compatible} \(\brk{m\cdot r}{n\cdot s}\cdot x=\brk{m}{n}\cdot((r,s)\cdot x)\) and \((\brk{m}{n}\cdot x)\cdot r=\brk{m}{n}\cdot(x\cdot r)\);
    \item\label{M:trilinear} \(\brk{m}{n}\cdot x\) is linear in each of \(m\), \(n\) and \(x\);
    \item\label{M:nil} \(\brk{m}{n}\cdot x=0\) whenever \(m\in[M]\), where \([M]\) denotes the subgroup of \(M\) generated by the elements \(\brk{h}{z}\cdot y\) for \(h\), \(z\in M\) and \(y\in\Ree\).
  \end{modaxioms}
  By~\ref{M:T}, also \(\brk{m}{n}\cdot x=0\) whenever \(n\in[M]\). A \defn{morphism of \(\R\)-modules} \(\alpha\colon M\to N\) is a group homomorphism such that
  \[
    \alpha(m\cdot r)=\alpha(m)\cdot r
    \quad\text{and}\quad
    \alpha(\brk{m}{m'}\cdot x)=\brk{\alpha(m)}{\alpha(m')}\cdot x
  \]
  for all \(m\), \(m'\in M\), \(r\in\Re\) and \(x\in\Ree\). Since these conditions are stable under composition, \(\R\)-modules form a category \(\Mod_{\R}\).
\end{definition}

The intuition is as follows. An element \(r\in\Re\) encodes a unary operation \(m\mapsto m\cdot r\) on \(M\), which by~\ref{M:sum} is \emph{quadratic}: its deviation from being additive is the bilinear operation \((m,n)\mapsto\brk{m}{n}\cdot H(r)\), so \(H\) assigns to a unary operation its cross-effect. An element \(x\in\Ree\) encodes a bilinear operation \((m,n)\mapsto\brk{m}{n}\cdot x\), and by~\ref{M:P} the map \(P\) assigns to it the associated squaring operation \(m\mapsto\brk{m}{m}\cdot x\). The involution \(T\) encodes the symmetry of bilinear operations, and~\ref{M:nil} says that brackets of brackets vanish.

Just as for rings, the first natural example of an \(\R\)-module is its underlying nearring \(\R_e\) itself:

\begin{example}\label{ex:R_e}
  For any square ring \(\R\), the nearring \(\Re\) is an \(\R\)-module with \(r\cdot s\) the nearring multiplication and \(\brk{r}{s}\cdot x\DefEq P((r,s)\cdot x)\). Indeed, \ref{M:action} is left distributivity and the monoid axioms; \ref{M:sum} is~\ref{SR:right-distributivity}; \ref{M:P} is~\ref{SR:P-binatural} with \(s=1\); \ref{M:T} follows from~\ref{SR:T-action} and~\ref{SR:PT}: \(P((r,s)\cdot T(x))=P(T((s,r)\cdot x))=P((s,r)\cdot x)\); \ref{M:compatible} follows from the monoid action and~\ref{SR:P-binatural}: \(P((r,s)\cdot x)\cdot t=P((r,s)\cdot x\cdot t)\); \ref{M:trilinear} holds since \(P\) is additive and the action is quadrilinear; and \ref{M:nil} holds since \([\Re]\subseteq\Im(P)\) and \(\brk{P(y)}{s}\cdot x=P((P(y),s)\cdot x)=0\) by~\ref{SR:P-kills}. Note that \([\Re]=\Im(P)\), because \(P(x)=P((1,1)\cdot x)=\brk{1}{1}\cdot x\), so that \(\overline{\Re}=\Rbar\).

  This \(\R\)-module, denoted simply by \(\Re\) in the sequel, actually is a \defn{free \(\R\)-module of rank one} whose unique basis element is the unit \(1\), see Proposition~\ref{prop:free-rank-one}.
\end{example}

We describe modules over further important examples of squares rings in Examples~\ref{ex:modules}, after stating two fundamental properties of modules over square rings:

\begin{lemma}\label{lem:module-basic}
  Let \(M\) be an \(\R\)-module. Then, for all \(m\), \(n\in M\), \(r\in\Re\) and \(x\in\Ree\):
  \begin{enumerate}
    \item\label{it:zero} \(0\cdot r=0\) and \(\brk{0}{n}\cdot x=0=\brk{m}{0}\cdot x\);
    \item\label{it:central} the elements \(\brk{m}{n}\cdot x\) are central in \(M\), so that \([M]\) is a central subgroup;
    \item\label{it:commutator} \(n+m-n-m=\brk{m}{n}\cdot H(2)\); hence \([M]\) contains the derived subgroup of \(M\), and \(M\) is nilpotent of class \(\leq 2\);
    \item\label{it:quotient} the quotient \(\overline{M}\DefEq M/[M]\) is an abelian group, and it is a right \(\Rbar\)-module via \(\overline{m}\cdot\overline{r}\DefEq\overline{m\cdot r}\), where \(\overline{m}\) denotes the class of \(m\).
  \end{enumerate}
\end{lemma}

\begin{proof}
  \ref{it:zero} By~\ref{M:trilinear}, \(\brk{0}{n}\cdot x=\brk{0+0}{n}\cdot x=\brk{0}{n}\cdot x+\brk{0}{n}\cdot x\), so \(\brk{0}{n}\cdot x=0\), and similarly \(\brk{m}{0}\cdot x=0\). Then by~\ref{M:sum}, \(0\cdot r=(0+0)\cdot r=0\cdot r+0\cdot r+\brk{0}{0}\cdot H(r)=0\cdot r+0\cdot r\), so \(0\cdot r=0\).

  \ref{it:central} Applying~\ref{M:sum} with \(r=2\) and then~\ref{M:action}, we find for all \(a\), \(b\in M\) that
  \begin{align*}
    \brk{a}{b}\cdot H(2) & =-b\cdot 2-a\cdot 2+(a+b)\cdot 2          \\
                         & =-b-b-a-a+a+b+a+b=-b+(-b-a+b+a)+b\text{.}
  \end{align*}
  Now let \(c=\brk{m}{n}\cdot x\) and \(b\in M\). By~\ref{M:nil}, \(\brk{c}{b}\cdot H(2)=0\), so the above gives \(-b-c+b+c=0\), i.e., \(c+b=b+c\).

  \ref{it:commutator} Since \(c=\brk{m}{n}\cdot H(2)\) is central, \ref{M:sum} with \(r=2\) reads \(m+n+m+n=m+m+n+n+c=m+c+m+n+n\), whence \(n+m=c+m+n\) and \(c=n+m-n-m\). Thus the derived subgroup of \(M\), generated by the commutators \(n+m-n-m\), lies in the central subgroup \([M]\), and \(M\) is nilpotent of class \(\leq 2\).

  \ref{it:quotient} Since \([M]\) is central and contains all commutators, \(\overline{M}\) is an abelian group. The action is well defined: for \(c\in[M]\) we have \((m+c)\cdot r=m\cdot r+c\cdot r+\brk{m}{c}\cdot H(r)=m\cdot r+c\cdot r\) by~\ref{M:nil}, and \(c\cdot r\in[M]\) since \((\brk{h}{z}\cdot y)\cdot r=\brk{h}{z}\cdot(y\cdot r)\) by~\ref{M:compatible}; it is additive in \(m\) by~\ref{M:sum}, additive and multiplicative in \(r\) by~\ref{M:action}, and it factors through \(\Rbar\) because \(m\cdot P(x)=\brk{m}{m}\cdot x\in[M]\) by~\ref{M:P}.
\end{proof}

\begin{proposition}\label{prop:restriction-of-scalars}
  Let \(f=(f_{e},f_{ee})\colon \R\to\R'\) be a morphism of square rings. Then a functor \(f^{*}\colon \Mod_{\R'}\to\Mod_{\R}\) is defined by \(f^{*}(M)=M\) as a group, with \(m\cdot r=m\cdot f_{e}(r)\) and \(\brk{m}{m'}\cdot x=\brk{m}{m'}\cdot f_{ee}(x)\) for \(m\), \(m'\in M\), \(r\in\Re\) and \(x\in\Ree\).\noproof
\end{proposition}

As for ordinary rings, this functor is called the \defn{restriction of scalars} along~\(f\); a left adjoint \emph{extension of scalars} along \(f\) is constructed in subsequent work.

Now we recall some of the numerous examples in~\cite{BHP} (thereby correcting a number of confusing typos there), and an additional one of group-theoretic interest described in less detail in~\cite{P2RM}, where it went largely unnoticed by the group theory community, likely because of the module-theoretic context of that paper.

\begin{examples}\label{ex:modules}
  \begin{enumerate}
    \item\label{ex:module-over-ring} Let \(R\) be a ring, viewed as a square ring as in Examples~\ref{ex:square-rings}\ref{ex:ring-as-square-ring}. Then the bracket operation vanishes, \ref{M:sum} says that each \(m\mapsto m\cdot r\) is additive, and Lemma~\ref{lem:module-basic}\ref{it:commutator} says that \(M\) is abelian. Modules over the square ring \(R\) are thus the usual right \(R\)-modules.
    \item\label{ex:module-Z-nil} A module over the nil square ring \(\Z_{\nil}\) of Examples~\ref{ex:square-rings}\ref{ex:Z-nil} is the same thing as a group of nilpotency class \(\leq 2\), with \(m\cdot a\) the \(a\)-fold sum of \(m\) and \(\brk{m}{n}\cdot x\) the \(x\)-fold sum of the commutator \(n+m-n-m\): see~\cite[Example~7.9]{BHP}. This gives an isomorphism of categories \(\Nil_{2}(\Grp)\cong\Mod_{\Z_{\nil}}\), a first instance of the Main Theorem.
    \item\label{ex:quotients-of-Z-nil} In~\cite[Subsection~8.1]{BHP} the finite quotient square rings of \(\Z_{\nil}\) are determined, that is, the square rings \(\R\) such that both components of the unique morphism of square rings \((\eta_{e},\eta_{ee})\colon \Z_{\nil}\to\R\) are surjective with finite images. They are of the form (after correction of typos in~\cite{BHP})
          \[
            \Z_{\nil}^{r,s}=\sqring{\Z/r\Z}{\Z/s\Z}{H}{P}
          \]
          with integers \(r\), \(s\geq 1\) such that \(s\) divides \(r\) if \(r\) is odd and \(2s\) divides \(r\) if \(r\) is even; moreover, \(H(\overline{n})=\overline{\binom{n}{2}}\) and \(P=0\). The modules over \(\Z_{\nil}^{r,s}\) are the groups of nilpotency class \(\leq 2\) and exponent dividing \(r\) whose commutators have order dividing \(s\).
    \item\label{ex:p-restricted} In particular, the modules over \(\Z_{\nil}^{p^{2},p}\) for a prime \(p\) are the \defn{\(p\)-restricted two-step nilpotent groups}, that is, the groups for which the third term of the \(p\)-restricted lower central series (i.e., of the fastest descending central series whose quotients are vector spaces over \(\Z/p\Z\)) is trivial.
    \item\label{ex:moore-spaces} The homotopy category of Moore spaces in degree \(2\) (that is, of simply connected CW-complexes whose reduced homology is concentrated in degree \(2\)) is a quadratic category~\cite{BHP} which plays an important role in low-dimensional homotopy theory~\cite{Baues:HomotopyType}.

    \item\label{ex:R-groups-binomial}
          The square ring \(\Z_{\nil}\) can be generalised by replacing \(\Z\) with any \defn{\(2\)-binomial ring} of coefficients, that is, a commutative ring \(R\) in which for every \(r\in R\) the element \(r(r-1)\) is uniquely \(2\)-divisible, so that the binomial coefficient \(\binom{r}{2}\DefEq r(r-1)/2\) is well defined; taking \(r=0\) we see that such a ring has no \(2\)-torsion. For such a ring \(R\) define the square ring
          \[
            R_{\nil}^{2}=\sqring{R}{R}{H_{R}^{2}}{0}
          \]
          by \((r,s)\cdot x\cdot t=rsxt\), \(T(x)=-x\) and \(H_{R}^{2}(r)=\binom{r}{2}\). The only axioms which are not immediate are~\ref{SR:H-sum} and~\ref{SR:H-product}, which amount to the identities \(\binom{r+s}{2}=\binom{r}{2}+\binom{s}{2}+rs\) and \(\binom{rs}{2}=r^{2}\binom{s}{2}+\binom{r}{2}s\). Then a module over \(R_{\nil}^{2}\) is the same thing as a \defn{nilpotent \(R\)-group of class \(2\)} as introduced by Philip Hall in~\cite{Hall:Edmonton}; we prove this in Example~\ref{ex:R-groups-any-ring} below, after recalling Hall's definition.
  \end{enumerate}
\end{examples}

For the convenience of the reader, and for comparison with the more general notion of \(R\)-group introduced in Example~\ref{ex:R-groups-any-ring}, let us recall Hall's notion explicitly, adapted to our (in particular additive) notation.

\begin{definition}\label{def:hall-R-group}
  Let \(R\) be a \(2\)-binomial ring. A \defn{nilpotent \(R\)-group of class \(2\)} is a group \(M\) of nilpotency class \(\leq 2\), endowed with a map \(M\times R\to M\colon (m,r)\mapsto m\cdot r\) satisfying the following properties:
  \begin{enumerate}
    \item\label{it:hall-action} \(m\cdot 1=m\), \(m\cdot(r+s)=m\cdot r+m\cdot s\) and \(m\cdot(rs)=(m\cdot r)\cdot s\) for all \(m\in M\) and \(r\), \(s\in R\);
    \item\label{it:hall-conjugation} \(n+m\cdot r-n=(n+m-n)\cdot r\) for all \(m\), \(n\in M\) and \(r\in R\);
    \item\label{it:hall-petresco} the \defn{Hall--Petresco formula}
          \[
            m_{1}\cdot r+\cdots+m_{k}\cdot r=(m_{1}+\cdots+m_{k})\cdot r+\bigl(\sum_{1\leq i<j\leq k}m_{i}+m_{j}-m_{i}-m_{j}\bigr)\cdot\binom{r}{2}
          \]
          holds for all \(k\geq 2\), all \((m_{1},\dots,m_{k})\in M^{k}\) and all \(r\in R\).
  \end{enumerate}
\end{definition}

Note that the case \(k=2\) of condition~\ref{it:hall-petresco} implies all the remaining cases by induction on \(k\), since the commutator \(m+n-m-n\) is biadditive and central in a group of nilpotency class \(\leq 2\).

Our last example, borrowed from~\cite{P2RM}, concerns the problem of extending the notion of a nilpotent \(R\)-group of class \(n\), introduced by Hall for \(n\)-binomial rings~\(R\) \cite{Hall:Edmonton} (see also~\cite{Warfield:Nilpotent}), to arbitrary commutative rings of coefficients. The class of \(n\)-binomial rings is not closed under important operations, quotients in particular, which prevents the development of a meaningful extension of classical commutative algebra to nilpotent \(R\)-groups as suggested by Hall. The exponential groups of Myasnikov and Remeslennikov~\cite{MR94}, defined over any ring and going back to work of Lyndon~\cite{Lyndon1960} and Baumslag~\cite{Baumslag:UniqueRoots}, lack a key feature of Hall's notion: they are not stable under central quotients, and for arbitrary \(R\) no stable substitute is known~\cite[Theorem~11.14 and Problem~1]{IKK}; see also~\cite{Amaglobeli-Miasnikov} for recent work on two-step nilpotent exponential groups. For \(n=2\), the notion of a module over a square ring allows us to solve this problem, as indicated in~\cite{P2RM} in a somewhat sketchy way; here we develop it more explicitly, as follows.

\begin{example}\label{ex:R-groups-any-ring}
  Let \(R\) be any commutative ring and \(I_{2}\) the ideal of \(R\) generated by the elements \(r(r-1)\). It is an \((R\wr\mathfrak{S}_{2})\)-\(R\)-bimodule, \(R\otimes R\) and \(R\) acting by multiplication in \(R\) and the involution \(\tau\) acting as \(-1\), or more explicitly,
  \[
    (r_{1}\otimes s_{1}+(r_{2}\otimes s_{2})\tau)\cdot x\cdot t=r_{1}s_{1}xt-r_{2}s_{2}xt\text{.}
  \]
  Then a square ring \(R_{\nil}=\sqring{R}{I_{2}}{H_{R}}{0}\) is defined by \(H_{R}(r)=r(r-1)\): as in Examples~\ref{ex:modules}\ref{ex:R-groups-binomial}, the only axioms which are not immediate are~\ref{SR:H-sum} and~\ref{SR:H-product}, which now amount to the identities \((r+s)(r+s-1)=r(r-1)+s(s-1)+2rs\) and \(rs(rs-1)=r^{2}s(s-1)+r(r-1)s\). We may thus define a \defn{two-step nilpotent \(R\)-group} to be a module over \(R_{\nil}\). (For \(R=\Z\) this square ring is \(\sqring{\Z}{2\Z}{H_{\Z}}{0}\), which is isomorphic but not equal to the nil square ring \(\Z_{\nil}\) of Examples~\ref{ex:square-rings}\ref{ex:Z-nil}; see below.)

  In order to compare this notion with Hall's, let \(M\) be a two-step nilpotent \(R\)-group. Note first that, by~\ref{M:sum} and~\ref{M:nil}, the map \(c\mapsto c\cdot r\) is additive on the central subgroup \([M]\), for every \(r\in R\). The action of \(R\) on \(M\) satisfies property~\ref{it:hall-action} of Definition~\ref{def:hall-R-group} by~\ref{M:action}. As to property~\ref{it:hall-conjugation}, let \(m\), \(n\in M\) and \(r\in R\), and write \(c\DefEq\brk{m}{n}\cdot H_{R}(2)=n+m-n-m\) (Lemma~\ref{lem:module-basic}\ref{it:commutator}), so that \(n+m-n=c+m\). Then \((c+m)\cdot r=c\cdot r+m\cdot r+\brk{c}{m}\cdot H_{R}(r)=c\cdot r+m\cdot r\) by~\ref{M:sum} and~\ref{M:nil}, while
  \begin{align*}
    c\cdot r & =\brk{m}{n}\cdot(H_{R}(2)\,r)         &  & \text{by~\ref{M:compatible}}                               \\
             & =\brk{m}{n}\cdot((r,1)\cdot H_{R}(2)) &  & \text{by commutativity of \(R\)}                           \\
             & =\brk{m\cdot r}{n}\cdot H_{R}(2)      &  & \text{by~\ref{M:compatible}}                               \\
             & =n+m\cdot r-n-m\cdot r                &  & \text{by Lemma~\ref{lem:module-basic}\ref{it:commutator},}
  \end{align*}
  so that \((n+m-n)\cdot r=n+m\cdot r-n-m\cdot r+m\cdot r=n+m\cdot r-n\). As to property~\ref{it:hall-petresco}, \ref{M:sum}, \ref{M:trilinear} and induction on \(k\) show that
  \[
    m_{1}\cdot r+\cdots+m_{k}\cdot r=(m_{1}+\cdots+m_{k})\cdot r-\sum_{1\leq i<j\leq k}\brk{m_{i}}{m_{j}}\cdot r(r-1)\text{.}
  \]
  Now if \(R\) is \(2\)-binomial, then \(r(r-1)=2\binom{r}{2}=H_{R}(2)\binom{r}{2}\), so that, by~\ref{M:compatible}, additivity of the action on \([M]\) and Lemma~\ref{lem:module-basic}\ref{it:commutator},
  \begin{align*}
    \sum_{1\leq i<j\leq k}\brk{m_{i}}{m_{j}}\cdot r(r-1)
     & =\sum_{1\leq i<j\leq k}\bigl(\brk{m_{i}}{m_{j}}\cdot H_{R}(2)\bigr)\cdot\binom{r}{2} \\
     & =\bigl(\sum_{1\leq i<j\leq k}\brk{m_{i}}{m_{j}}\cdot H_{R}(2)\bigr)\cdot\binom{r}{2} \\
     & =\bigl(\sum_{1\leq i<j\leq k}m_{j}+m_{i}-m_{j}-m_{i}\bigr)\cdot\binom{r}{2}          \\
     & =-\bigl(\sum_{1\leq i<j\leq k}m_{i}+m_{j}-m_{i}-m_{j}\bigr)\cdot\binom{r}{2}\text{,}
  \end{align*}
  and the Hall--Petresco formula holds. Thus, when \(R\) is \(2\)-binomial, a two-step nilpotent \(R\)-group in our sense is a nilpotent \(R\)-group of class \(2\) in the sense of Hall.

  From a more conceptual viewpoint, note that for a \(2\)-binomial ring \(R\) we have \(I_{2}=2R\), since \(2=H_{R}(2)\in I_{2}\) and \(r(r-1)=2\binom{r}{2}\). Hence the \(R\)-linear map \(\mu_{2}\colon R\to I_{2}\colon r\mapsto 2r\) is surjective, and it is injective since \(R\) has no \(2\)-torsion. The pair \((1_{R},\mu_{2})\colon R_{\nil}^{2}\to R_{\nil}\) is a morphism of square rings, because
  \[
    \mu_{2}\binom{r}{2}=r(r-1)=H_{R}(r)
  \]
  and \(\mu_{2}(rsxt)=(r,s)\cdot\mu_{2}(x)\cdot t\); so it is an isomorphism of square rings, and the restriction of scalars \((1_{R},\mu_{2})^{*}\colon \Mod_{R_{\nil}}\to\Mod_{R_{\nil}^{2}}\) of Proposition~\ref{prop:restriction-of-scalars} is an isomorphism of categories, with inverse \(((1_{R},\mu_{2})^{-1})^{*}\). It sends a two-step nilpotent \(R\)-group \(M\) to the \(R_{\nil}^{2}\)-module with the same action of \(R\) and with bracket \(\brk{m}{n}\cdot y=\brk{m}{n}\cdot 2y=(\brk{m}{n}\cdot H_{R}(2))\cdot y=(n+m-n-m)\cdot y\) for \(y\in R\), by~\ref{M:compatible} and Lemma~\ref{lem:module-basic}\ref{it:commutator}.

  Conversely, let \(R\) be \(2\)-binomial and let \(M\) be a nilpotent \(R\)-group of class \(2\) in the sense of Hall. Then \(M\) becomes an \(R_{\nil}^{2}\)-module when we define \(\brk{m}{n}\cdot y\DefEq(n+m-n-m)\cdot y\) for \(m\), \(n\in M\) and \(y\in R\). Indeed, the case \(k=2\) of~\ref{it:hall-petresco} shows that \((c+m)\cdot r=c\cdot r+m\cdot r\) whenever \(c\) is central in \(M\), so that the action of \(R\) is additive on the centre of \(M\); in particular \(0\cdot r=0\) and \((-c)\cdot r=-(c\cdot r)\) for central \(c\). Now \ref{M:action} is~\ref{it:hall-action}; \ref{M:sum} is the case \(k=2\) of~\ref{it:hall-petresco}, since \(-(m+n-m-n)=n+m-n-m\); \ref{M:P} holds because \(P=0\), \(m\cdot 0=0\) and \(m+m-m-m=0\); \ref{M:T} holds because \(T=-1\) and \(-(n+m-n-m)=m+n-m-n\); \ref{M:trilinear} holds because the commutator is biadditive in a group of nilpotency class \(\leq 2\), and the action is additive in the scalar and on the centre; and \ref{M:nil} holds because \([M]\) consists of central elements, by~\ref{it:hall-conjugation}. For~\ref{M:compatible}, the second identity is the associativity in~\ref{it:hall-action}, while for the first one, \ref{it:hall-conjugation} and additivity on the centre give
  \begin{align*}
    n+m\cdot r-n-m\cdot r & =(n+m-n)\cdot r-m\cdot r                                  \\
                          & =\bigl((n+m-n-m)+m\bigr)\cdot r-m\cdot r=(n+m-n-m)\cdot r
  \end{align*}
  for all \(m\), \(n\in M\) and \(r\in R\). Applying this twice, and using that \(n\cdot s+m-n\cdot s-m=-(m+n\cdot s-m-n\cdot s)=-((m+n-m-n)\cdot s)=(n+m-n-m)\cdot s\), we find
  \begin{align*}
    n\cdot s+m\cdot r-n\cdot s-m\cdot r & =(n\cdot s+m-n\cdot s-m)\cdot r                                 \\
                                        & =\bigl((n+m-n-m)\cdot s\bigr)\cdot r=(n+m-n-m)\cdot(sr)\text{,}
  \end{align*}
  which is \(\brk{m\cdot r}{n\cdot s}\cdot y=\brk{m}{n}\cdot(rsy)\) for \(y=1\), and then for all \(y\in R\) by~\ref{it:hall-action}. The two constructions are inverse to each other, because in any \(R_{\nil}^{2}\)-module the bracket is determined by the action: \(\brk{m}{n}\cdot y=(\brk{m}{n}\cdot H_{R}^{2}(2))\cdot y=(n+m-n-m)\cdot y\) by~\ref{M:compatible} and Lemma~\ref{lem:module-basic}\ref{it:commutator}, as \(H_{R}^{2}(2)=1\). Since a morphism of \(R_{\nil}^{2}\)-modules is then simply a group homomorphism which preserves the action of \(R\), this proves the claim made in Examples~\ref{ex:modules}\ref{ex:R-groups-binomial}: the category \(\Mod_{R_{\nil}^{2}}\), and hence also \(\Mod_{R_{\nil}}\), is isomorphic to the category of nilpotent \(R\)-groups of class \(2\) in the sense of Hall. So our definition of a two-step nilpotent \(R\)-group might eventually inspire a way to define \(n\)-step nilpotent \(R\)-groups, and thus nilpotent \(R\)-groups of any class \(n\), for an arbitrary commutative ring of coefficients \(R\).
\end{example}

Just as \(\Z_{\nil}\) is an initial object of \(\SqRing\), the square ring \(R_{\nil}\) is initial amongst the square rings which receive a morphism from it; we now describe those square rings, which we think of as square rings under \(R\), the square ring \(R_{\nil}\) playing the role of \(R\) itself. Note first that the construction of \(R_{\nil}\) is functorial: a homomorphism of commutative rings \(\alpha\colon R\to S\) sends each generator \(r(r-1)\) of the ideal \(I_{2}\) of \(R\) to the generator \(\alpha(r)(\alpha(r)-1)\) of the corresponding ideal of \(S\), so that \(\alpha_{\nil}\DefEq(\alpha,\alpha|_{I_{2}})\colon R_{\nil}\to S_{\nil}\) is a morphism of square rings, as \(H_{S}\circ\alpha=\alpha\circ H_{R}\), and we obtain a functor \((-)_{\nil}\colon \CRing\to\SqRing\) from the category of commutative rings.

\begin{definition}\label{def:square-ring-under-R}
  Let \(R\) be a commutative ring. A \defn{square ring under \(R\)} is a square ring \(\R\) equipped with a morphism of square rings \(\eta=(\eta_{e},\eta_{ee})\colon R_{\nil}\to\R\), and a \defn{morphism \((\R,\eta)\to(\R',\eta')\) of square rings under \(R\)} is a morphism of square rings \(f\colon \R\to\R'\) such that \(f\circ\eta=\eta'\). In other words, the category \(\SqRing_{R}\DefEq R_{\nil}/\SqRing\) of square rings under \(R\) is the coslice category of \(\SqRing\) under \(R_{\nil}\), and \(R_{\nil}\), equipped with its identity, is an initial object of \(\SqRing_{R}\).
\end{definition}

Square rings under \(R\) are the quadratic counterpart of rings under \(R\): just as every module over a ring under \(R\) is an \(R\)-module, via Proposition~\ref{prop:restriction-of-scalars}, restriction of scalars along~\(\eta\) makes every module over a square ring under \(R\) a two-step nilpotent \(R\)-group, and every morphism between such modules a morphism of two-step nilpotent \(R\)-groups.

To be a square ring under \(R\) is a property of a square ring \(\R\) together with a morphism of nearrings \(R\to\Re\), as we now show. For such a morphism \(\eta\colon R\to\Re\), write \(\Ree^{\eta}\) for the set of those \(x\in\Ree\) which satisfy \((\eta(r),1)\cdot x=(1,\eta(r))\cdot x=x\cdot\eta(r)\) for all \(r\in R\). It is a subgroup of \(\Ree\) by quadrilinearity, and \((\eta(s),1)\cdot x\in\Ree^{\eta}\) whenever \(x\in\Ree^{\eta}\) and \(s\in R\), because the three actions commute with one another and \(\eta\) preserves products; hence \(\Ree^{\eta}\) is an \(R\)-module via \(r\cdot x\DefEq(\eta(r),1)\cdot x=(1,\eta(r))\cdot x=x\cdot\eta(r)\).

\begin{proposition}\label{prop:morphisms-from-R-nil}
  Let \(R\) be a commutative ring, \(\R=\sqring{\Re}{\Ree}{H}{P}\) a square ring and \(\eta\colon R\to\Re\) a morphism of nearrings. The following are equivalent:
  \begin{tfae}
    \item\label{it:eta-first-component} \(\eta\) is the first component of a morphism of square rings \(R_{\nil}\to\R\);
    \item\label{it:H-eta-balanced} \(H(\eta(r))\in\Ree^{\eta}\) for all \(r\in R\).
  \end{tfae}
  If these conditions hold, then the morphism in~\ref{it:eta-first-component} is unique. In particular, they hold when \(\Ree^{\eta}=\Ree\), i.e., when the three actions of \(R\) on \(\Ree\) induced by \(\eta\) coincide.
\end{proposition}

\begin{proof}
  \ref{it:H-eta-balanced}~\(\Rightarrow\)~\ref{it:eta-first-component} Recall from~\cite[Theorem~3.6]{P2RM} that the map \(H_{R}\colon R\to I_{2}\) is a universal \defn{quadratic derivation}, i.e., it is initial among all maps \(d\colon R\to N\) to \(R\)-modules \(N\) which satisfy the identities \(d(r+s)=d(r)+d(s)+rs\,d(2)\) and \(d(rs)=r^{2}d(s)+s\,d(r)\). By~\ref{it:H-eta-balanced}, the map \(d\DefEq H\circ\eta\) takes values in the \(R\)-module \(\Ree^{\eta}\), and there it satisfies these identities: by~\ref{SR:H-sum}, \(d(r+s)=d(r)+d(s)+(\eta(s),\eta(r))\cdot d(2)=d(r)+d(s)+rs\,d(2)\) because \(d(2)\in\Ree^{\eta}\), and by~\ref{SR:H-product}, \(d(rs)=(\eta(r),\eta(r))\cdot d(s)+d(r)\cdot\eta(s)=r^{2}d(s)+s\,d(r)\) because \(d(r)\), \(d(s)\in\Ree^{\eta}\). Hence there exists a unique homomorphism of \(R\)-modules \(\tilde{d}\colon I_{2}\to\Ree^{\eta}\) such that \(\tilde{d}\circ H_{R}=H\circ\eta\). We check that the pair \((\eta,\tilde{d})\colon R_{\nil}\to\R\) is a morphism of square rings. Let \(r\), \(s\), \(t\in R\) and \(x\in I_{2}\). Since \(\tilde{d}(x)\in\Ree^{\eta}\), we have \(\tilde{d}((r,s)\cdot x\cdot t)=\tilde{d}(rstx)=rst\,\tilde{d}(x)=(\eta(r),\eta(s))\cdot\tilde{d}(x)\cdot\eta(t)\). Moreover, by~\ref{SR:right-distributivity},
  \begin{align*}
    \eta((r+s)t) & =(\eta(r)+\eta(s))\,\eta(t)                                                   \\
                 & =\eta(r)\eta(t)+\eta(s)\eta(t)+P\bigl((\eta(r),\eta(s))\cdot H(\eta(t))\bigr) \\
                 & =\eta((r+s)t)+P\bigl(rs\,\tilde{d}(H_{R}(t))\bigr)                            \\
                 & =\eta((r+s)t)+P\bigl(\tilde{d}(rst(t-1))\bigr)\text{,}
  \end{align*}
  whence \(P\circ\tilde{d}\) vanishes on the elements \(rst(t-1)\), which generate \(I_{2}\) as an abelian group; since \(P\circ\tilde{d}\) is additive, \(P\circ\tilde{d}=0=\eta\circ 0\).

  \ref{it:eta-first-component}~\(\Rightarrow\)~\ref{it:H-eta-balanced} Let \(f=(\eta,f_{ee})\colon R_{\nil}\to\R\) be a morphism of square rings. For \(a\), \(u\in R\) we have \(a\,H_{R}(u)=(a,1)\cdot H_{R}(u)=(1,a)\cdot H_{R}(u)=H_{R}(u)\cdot a\) in \(I_{2}\). Applying \(f_{ee}\) and using \(f_{ee}\circ H_{R}=H\circ\eta\), we find \((\eta(a),1)\cdot H(\eta(u))=(1,\eta(a))\cdot H(\eta(u))=H(\eta(u))\cdot\eta(a)\), which is~\ref{it:H-eta-balanced}. The same computation shows that \(f_{ee}(a\,u(u-1))=(\eta(a),1)\cdot H(\eta(u))\); since the elements \(a\,u(u-1)\) generate \(I_{2}\) as an abelian group and \(f_{ee}\) is additive, \(f_{ee}\) is determined by \(\eta\), which proves uniqueness.
\end{proof}

\begin{remark}\label{rem:square-ring-under-R}
  By Proposition~\ref{prop:morphisms-from-R-nil}, a square ring under \(R\) is determined by the square ring \(\R\) and the morphism of nearrings \(\eta_{e}\colon R\to\Re\). We therefore write it as \((R\xrightarrow{\eta}\Re\xrightarrow{H}\Ree\xrightarrow{P}\Re)\), where \(\eta=\eta_{e}\); a morphism of square rings under \(R\) is then a morphism of square rings \(f\) such that \(f_{e}\circ\eta=\eta'\). Condition~\ref{it:H-eta-balanced} is not automatic, even when \(\Re\) is a commutative ring and \(\eta\) is an isomorphism. Let \(R=\Z[t]\), let \(\phi\colon R\to R\colon p\mapsto p(0)\) be the evaluation at \(0\), and consider \(\R=\sqring{R}{R}{H}{0}\) with \((r,s)\cdot x\cdot u\DefEq\phi(r)\phi(s)xu\), \(T=-1\) and \(H(p)\DefEq\phi(p)^{2}-p\). Then \(H(2)=2\), and \ref{SR:H-sum}, \ref{SR:H-product} and \ref{SR:T-action} hold because \(\phi\) is a ring homomorphism and \(R\) is commutative; the remaining axioms are immediate, as \(P=0\). For \(\eta=1_{R}\), condition~\ref{it:H-eta-balanced} fails: \(H(t)=-t\), whereas \((t,1)\cdot(-t)=\phi(t)(-t)=0\neq-t^{2}=(-t)\cdot t\). So no morphism of square rings \(R_{\nil}\to\R\) has the identity of \(R\) as its first component.
\end{remark}

For \(R=\Z\), every square ring \(\R\) admits exactly one morphism of nearrings \(\eta_{\R}\colon \Z\to\Re\), which sends \(1\) to the unit of \(\Re\), and condition~\ref{it:H-eta-balanced} of Proposition~\ref{prop:morphisms-from-R-nil} holds for it by quadrilinearity: \((n,1)\cdot x=nx=(1,n)\cdot x=x\cdot n\) for all \(n\in\Z\) and \(x\in\Ree\). Hence every square ring is a square ring under \(\Z\) in exactly one way, the square ring \((\Z\xrightarrow{1_{\Z}}\Z\xrightarrow{H_{\Z}}2\Z\xrightarrow{0}\Z)\) is an initial object of \(\SqRing\), and the forgetful functor \(\SqRing_{\Z}\to\SqRing\) is an isomorphism of categories. Via the isomorphism
\[
  (1_{\Z},\mu_{2})\colon (\Z\xrightarrow{1_{\Z}}\Z\xrightarrow{H_{\Z}^{2}}\Z\xrightarrow{0}\Z)\to(\Z\xrightarrow{1_{\Z}}\Z\xrightarrow{H_{\Z}}2\Z\xrightarrow{0}\Z)
\]
of Example~\ref{ex:R-groups-any-ring}, this recovers the fact, due to~\cite{BHP}, that the nil square ring \(\Z_{\nil}\) is an initial object of \(\SqRing\). Thus \(R_{\nil}\) may be considered as a generalisation of \(\Z_{\nil}\) in which \(\Z\) is replaced by an arbitrary commutative ring of coefficients \(R\).

\subsection{Coproducts and free modules}\label{sec:coproducts}
Let \(\R\) be a square ring. For \(\R\)-modules \(M\) and \(N\), Lemma~\ref{lem:module-basic}\ref{it:quotient} and Lemma~\ref{lem:actions-factor} allow us to define the abelian group
\[
  A(M,N)\DefEq\overline{M}\otimes\overline{N}\otimes_{\Rbar\otimes\Rbar}\Ree\text{,}
\]
where \(\overline{M}\otimes\overline{N}\) is the tensor product of abelian groups, a right \((\Rbar\otimes\Rbar)\)-module via \((\overline{m}\otimes\overline{n})\cdot(\overline{r}\otimes\overline{s})=\overline{m}\cdot\overline{r}\otimes\overline{n}\cdot\overline{s}\). Thus \(A(M,N)\) is generated by elements \(\overline{m}\otimes\overline{n}\otimes x\), additive in each variable, subject to \(\overline{m\cdot r}\otimes\overline{n\cdot s}\otimes x=\overline{m}\otimes\overline{n}\otimes(r,s)\cdot x\); it is a right \(\Rbar\)-module via the right action on \(\Ree\). The following description of binary coproducts in \(\Mod_{\R}\) generalises both the construction for \(\Nil_{2}(\Grp)\) given in~\cite{Jibladze-Pirashvili} and its analogue for the two-nilpotent core of the category of algebras over a reduced symmetric operad in the monoidal category of modules over a commutative ring given in~\cite{DHVdL}.

\begin{proposition}\label{prop:coproduct}
  Let \(M\) and \(N\) be \(\R\)-modules. The set \(M+N\DefEq M\times N\times A(M,N)\) is an \(\R\)-module with the operations
  \begin{align*}
    (m,n,u)+(m',n',u')               & \DefEq(m+m',\,n+n',\,u+u'+\overline{m'}\otimes\overline{n}\otimes H(2))\text{,}                                                                  \\
    (m,n,u)\cdot r                   & \DefEq(m\cdot r,\,n\cdot r,\,u\cdot\overline{r}+\overline{m}\otimes\overline{n}\otimes H(r))\text{,}                                             \\
    \brk{(m,n,u)}{(m',n',u')}\cdot x & \DefEq(\brk{m}{m'}\cdot x,\,\brk{n}{n'}\cdot x,\,\overline{m}\otimes\overline{n'}\otimes x+\overline{m'}\otimes\overline{n}\otimes T(x))\text{,}
  \end{align*}
  and it is the coproduct of \(M\) and \(N\) in \(\Mod_{\R}\), with coproduct injections \(i_{M}(m)=(m,0,0)\) and \(i_{N}(n)=(0,n,0)\). Its neutral element is \((0,0,0)\), the inverse of \((m,n,u)\) is \((-m,-n,-u+\overline{m}\otimes\overline{n}\otimes H(2))\), and \([M+N]=[M]\times[N]\times A(M,N)\), so that \(\overline{M+N}=\overline{M}\times\overline{N}\).
\end{proposition}

\begin{proof}
  Write \(\delta(a,a')\DefEq\overline{m'}\otimes\overline{n}\otimes H(2)\) for \(a=(m,n,u)\) and \(a'=(m',n',u')\). This term is biadditive in \((\overline{m'},\overline{n})\), and it satisfies \(\delta(a,a')=0\) as soon as \(m'\in[M]\) or \(n\in[N]\). The associativity of \(+\) amounts to the cocycle identity \(\delta(a,a')+\delta(a+a',a'')=\delta(a',a'')+\delta(a,a'+a'')\), which holds by biadditivity; \((0,0,0)\) is neutral since \(\delta\) vanishes when \(m'=0\) or \(n=0\); and the given element is a two-sided inverse of \((m,n,u)\) because \(\overline{-m}=-\overline{m}\). So \(M+N\) is a group.

  We verify the module axioms. Throughout we use \(H(1)=0\), \((2,2)\cdot x=4x\), \(\overline{P(x)}=0\), the balancing relation of \(A(M,N)\), and that \(\overline{c}=0\) for \(c\in[M]\) or \(c\in[N]\). The first two components of each identity hold because they hold in \(M\) and in \(N\); we only discuss the third component, writing \(a=(m,n,u)\) and \(a'=(m',n',u')\).
  \begin{itemize}
    \item \ref{M:action}: \(a\cdot 1=a\) since \(H(1)=0\). For \((a\cdot r)\cdot s=a\cdot(r\cdot s)\) we need \(\overline{m}\otimes\overline{n}\otimes H(r)\cdot s+\overline{m\cdot r}\otimes\overline{n\cdot r}\otimes H(s)=\overline{m}\otimes\overline{n}\otimes H(r\cdot s)\), which is~\ref{SR:H-product}. For \(a\cdot(r+s)=a\cdot r+a\cdot s\) we need \(\overline{m}\otimes\overline{n}\otimes H(r+s)=\overline{m}\otimes\overline{n}\otimes H(r)+\overline{m}\otimes\overline{n}\otimes H(s)+\overline{m\cdot s}\otimes\overline{n\cdot r}\otimes H(2)\), which is~\ref{SR:H-sum}.
    \item \ref{M:sum}: expanding both sides of \((a+a')\cdot r=a\cdot r+a'\cdot r+\brk{a}{a'}\cdot H(r)\), all terms agree except that the left-hand side contains \(\overline{m'}\otimes\overline{n}\otimes(H(2)\cdot r+H(r))\) where the right-hand side contains \(\overline{m'}\otimes\overline{n}\otimes((r,r)\cdot H(2)+T(H(r)))\). These agree because
          \begin{align*}
            H(2)\cdot r & =H(2\cdot r)-(2,2)\cdot H(r)=H(r+r+P(H(r)))-4H(r) \\
                        & =H(r)+H(r)+(r,r)\cdot H(2)+H(P(H(r)))-4H(r)       \\
                        & =(r,r)\cdot H(2)+T(H(r))+H(r)-2H(r)               \\
                        & =(r,r)\cdot H(2)+T(H(r))-H(r)
          \end{align*}
          by~\ref{SR:H-product}, \ref{SR:right-distributivity}, \ref{SR:H-sum}, \ref{SR:P-kills} and \ref{SR:T}.
    \item \ref{M:P}: \(a\cdot P(x)\) has third component
          \[
            u\cdot\overline{P(x)}+\overline{m}\otimes\overline{n}\otimes H(P(x))=\overline{m}\otimes\overline{n}\otimes(x+T(x))\text{,}
          \]
          which is the third component of \(\brk{a}{a}\cdot x\).
    \item \ref{M:T}, \ref{M:compatible} and \ref{M:trilinear} follow from \(T\circ T=1_{\Ree}\), \ref{SR:T-action} (which also gives \(T(x\cdot r)=T(x)\cdot r\)), quadrilinearity and the balancing relation; the terms \(\delta\) which appear when adding brackets vanish since brackets have first two components in \([M]\) and \([N]\).
    \item \ref{M:nil}: the brackets generate \([M]\times[N]\times A(M,N)\), which is therefore \([M+N]\), and for \(a\) in this subgroup \(\overline{m}=0=\overline{n}\), so \(\brk{a}{a'}\cdot x=0\).
  \end{itemize}
  The coproduct injections are module morphisms because all correction terms vanish when \(n=0\) and \(m=0\), respectively.

  Let \(f\colon M\to Z\) and \(g\colon N\to Z\) be module morphisms. By~\ref{M:trilinear}, \ref{M:nil} and~\ref{M:compatible}, the assignment \(\overline{m}\otimes\overline{n}\otimes x\mapsto\brk{f(m)}{g(n)}\cdot x\) is well defined and extends to a group homomorphism \(\beta\colon A(M,N)\to[Z]\), with \(\beta(u\cdot\overline{r})=\beta(u)\cdot r\). Define \(h(m,n,u)\DefEq f(m)+g(n)+\beta(u)\). Then \(h\) is a group homomorphism: since \(g(n)+f(m')=f(m')+g(n)+(g(n)+f(m')-g(n)-f(m'))\) and the commutator equals \(\brk{f(m')}{g(n)}\cdot H(2)=\beta(\delta(a,a'))\) by Lemma~\ref{lem:module-basic}\ref{it:commutator}, we get \(h(a)+h(a')=h(a+a')\). Using that \(\beta(u)\in[Z]\) is central and~\ref{M:sum}, \ref{M:nil}, \ref{M:trilinear}, \ref{M:T},
  \begin{align*}
    h(a)\cdot r              & =(f(m)+g(n))\cdot r+\beta(u)\cdot r                                                               \\
                             & =f(m)\cdot r+g(n)\cdot r+\brk{f(m)}{g(n)}\cdot H(r)+\beta(u\cdot\overline{r})=h(a\cdot r)\text{,} \\
    \brk{h(a)}{h(a')}\cdot x & =\brk{f(m)+g(n)}{f(m')+g(n')}\cdot x                                                              \\
                             & =\brk{f(m)}{f(m')}\cdot x+\brk{g(n)}{g(n')}\cdot x                                                \\
                             & \qquad+\brk{f(m)}{g(n')}\cdot x+\brk{f(m')}{g(n)}\cdot T(x)\text{,}
  \end{align*}
  which is \(h(\brk{a}{a'}\cdot x)\). Clearly \(h\circ i_{M}=f\) and \(h\circ i_{N}=g\). Conversely, any module morphism \(h'\) with these properties agrees with \(h\): we have \((m,n,u)=i_{M}(m)+i_{N}(n)+(0,0,u)\), the map \(u\mapsto(0,0,u)\) is additive, and \((0,0,\overline{m'}\otimes\overline{n'}\otimes x)=\brk{i_{M}(m')}{i_{N}(n')}\cdot x\), so \(h'(0,0,\overline{m'}\otimes\overline{n'}\otimes x)=\brk{f(m')}{g(n')}\cdot x=\beta(\overline{m'}\otimes\overline{n'}\otimes x)\).
\end{proof}

\begin{proposition}\label{prop:free-rank-one}
  For any \(\R\)-module \(M\), the map \(\Mod_{\R}(\Re,M)\to M\colon {f\mapsto f(1)}\) is a bijection, with inverse sending \(m\in M\) to \(f_{m}\colon \Re\to M\colon r\mapsto m\cdot r\). Hence \(\Re\) is the free \(\R\)-module on one generator.
\end{proposition}

\begin{proof}
  The map \(f_{m}\) is a module morphism: it is additive and preserves the action by~\ref{M:action}, and it preserves brackets by~\ref{M:P} and~\ref{M:compatible}:
  \[
    f_{m}(\brk{r}{s}\cdot x)=m\cdot P((r,s)\cdot x)=\brk{m}{m}\cdot((r,s)\cdot x)=\brk{m\cdot r}{m\cdot s}\cdot x\text{.}
  \]
  We have \(f_{m}(1)=m\), and \(f(r)=f(1\cdot r)=f(1)\cdot r=f_{f(1)}(r)\) for any module morphism \(f\).
\end{proof}

\begin{corollary}\label{cor:free-modules}
  For \(n\geq 0\), let \(\Re^{(n)}\DefEq\Re+\cdots+\Re\) be the coproduct of \(n\) copies of \(\Re\) in \(\Mod_{\R}\) (so \(\Re^{(0)}=0\)), the free \(\R\)-module on the generators \(e_{1}\), \dots, \(e_{n}\), where \(e_{k}\) is the image of \(1\) under the \(k\)-th coproduct injection. Every element of \(\Re^{(n)}\) can be written in a unique way as
  \begin{equation}\label{eq:normal-form-module}
    \sum_{k=1}^{n}e_{k}\cdot r_{k}+\sum_{1\leq k<l\leq n}\brk{e_{k}}{e_{l}}\cdot x_{kl}
  \end{equation}
  with \(r_{k}\in\Re\) and \(x_{kl}\in\Ree\), where the second sum is central and the first is taken in the order \(k=1\), \dots, \(n\). In other words, the map \(\Re^{n}\times\Ree^{\binom{n}{2}}\to\Re^{(n)}\) sending \(((r_{k})_{k},(x_{kl})_{k<l})\) to~\eqref{eq:normal-form-module} is a bijection. Moreover, \(\overline{\Re^{(n)}}\cong\Rbar^{n}\).
\end{corollary}

\begin{proof}
  By induction on \(n\). The case \(n=0\) is trivial, and the case \(n=1\) is Proposition~\ref{prop:free-rank-one} together with Example~\ref{ex:R_e}. Suppose the claims hold for \(n\), and apply Proposition~\ref{prop:coproduct} to the coproduct \(\Re^{(n)}+\Re=\Re^{(n+1)}\). Since \(\overline{\Re^{(n)}}\cong\Rbar^{n}\) and \(\overline{\Re}=\Rbar\), the abelian group \(A(\Re^{(n)},\Re)\) is
  \[
    \Rbar^{n}\otimes\Rbar\otimes_{\Rbar\otimes\Rbar}\Ree\cong\Ree^{n}\text{,}
  \]
  an isomorphism being given by \(\sum_{k}\overline{e_{k}}\otimes\overline{1}\otimes x_{k}\mapsto(x_{k})_{k}\). The element \((a,r,(x_{k})_{k})\) of \(\Re^{(n)}\times\Re\times\Ree^{n}\) corresponds to \(i_{1}(a)+i_{2}(r)+\sum_{k}\brk{e_{k}}{e_{n+1}}\cdot x_{k}\) in \(\Re^{(n+1)}\), by the description of the elements \((0,0,u)\) in the proof of Proposition~\ref{prop:coproduct}; substituting the normal form of \(a\) gives~\eqref{eq:normal-form-module} for \(n+1\). The last claim follows from \(\overline{M+N}=\overline{M}\times\overline{N}\).
\end{proof}

\begin{remark}\label{rem:freeness}
  Freeness of \(\Re^{(n)}\) means that for every \(\R\)-module \(M\), the map
  \[
    \Mod_{\R}(\Re^{(n)},M)\to M^{n}\colon v\mapsto(v(e_{1}),\dots,v(e_{n}))
  \]
  is a bijection. Given \(\mathbf{m}=(m_{1},\dots,m_{n})\in M^{n}\), we write \(\ev_{\mathbf{m}}\colon \Re^{(n)}\to M\) for the module morphism with \(\ev_{\mathbf{m}}(e_{k})=m_{k}\); on the normal form~\eqref{eq:normal-form-module} it is given by
  \[
    \ev_{\mathbf{m}}\bigl(\sum_{k}e_{k}\cdot r_{k}+\sum_{k<l}\brk{e_{k}}{e_{l}}\cdot x_{kl}\bigr)=\sum_{k}m_{k}\cdot r_{k}+\sum_{k<l}\brk{m_{k}}{m_{l}}\cdot x_{kl}\text{.}
  \]
\end{remark}

\section{Categories of modules are two-step nilpotent varieties}\label{sec:modules-nilpotent}

Throughout this section, \(\R\) is a square ring. We prove one implication of the Main Theorem.

\begin{proposition}\label{prop:mod-semi-abelian}
  The category \(\Mod_{\R}\) is a variety of \(\Omega\)-groups in the sense of Higgins~\cite{Higgins}; in particular, it is semi-abelian.
\end{proposition}

\begin{proof}
  An \(\R\)-module is a group \(M\) equipped with a unary operation \(\rho_{r}(m)\DefEq m\cdot r\) for each \(r\in\Re\) and a binary operation \(\phi_{x}(m,n)\DefEq\brk{m}{n}\cdot x\) for each \(x\in\Ree\), and the axioms of Definition~\ref{def:module} are identities in these operations: \ref{M:nil} amounts, via~\ref{M:trilinear}, to the identities \(\phi_{x}(\phi_{y}(h,z),n)=0\). Module morphisms are precisely the maps preserving the operations. So \(\Mod_{\R}\) is a variety of algebras, and since \(\rho_{r}(0)=0\) and \(\phi_{x}(0,0)=0\) by Lemma~\ref{lem:module-basic}\ref{it:zero}, it is a variety of \(\Omega\)-groups. Every such variety is semi-abelian~\cite[Example~2.6(2)]{Janelidze-Marki-Tholen}.
\end{proof}

Being a variety, \(\Mod_{\R}\) has free objects; the free module on one generator is \(\Re\) by Proposition~\ref{prop:free-rank-one}, and the free module on \(n\) generators is \(\Re^{(n)}\) of Corollary~\ref{cor:free-modules}.

\begin{proposition}\label{prop:cosmash-mod}
  For \(\R\)-modules \(M\) and \(N\), the cosmash product \(M\cosmash N\) is the abelian group \(A(M,N)=\overline{M}\otimes\overline{N}\otimes_{\Rbar\otimes\Rbar}\Ree\), with zero bracket and with \(\Re\) acting through \(\Rbar\), and the kernel morphism \(\iota_{M,N}\colon M\cosmash N\to M+N\) is \(u\mapsto(0,0,u)\). The morphism \(\cc{M}\colon M\cosmash M\to M\) is given on generators by \(\overline{m}\otimes\overline{m'}\otimes x\mapsto\brk{m}{m'}\cdot x\). Consequently, the Higgins commutator \([M,M]\) is the subgroup \([M]\), the abelianisation of \(M\) is \(\ab(M)=\overline{M}\), the abelian objects of \(\Mod_{\R}\) are the modules with zero bracket, and \(\Ab(\Mod_{\R})\cong\Mod_{\Rbar}\).
\end{proposition}

\begin{proof}
  In the description of Proposition~\ref{prop:coproduct}, the retractions \(r_{1}=\langle 1_{M},0\rangle\) and \(r_{2}=\langle 0,1_{N}\rangle\) are the projections \((m,n,u)\mapsto m\) and \((m,n,u)\mapsto n\), since the map \(\beta\) in the proof of that proposition vanishes when \(f=0\) or \(g=0\). Hence \(\Ker(\cm{M}{N})=\{(0,0,u)\mid u\in A(M,N)\}\), which by the formulas of Proposition~\ref{prop:coproduct} carries the structure described. The codiagonal \(\fold{M}=\langle 1_{M},1_{M}\rangle\) sends \((0,0,\overline{m}\otimes\overline{m'}\otimes x)\) to \(\beta(\overline{m}\otimes\overline{m'}\otimes x)=\brk{m}{m'}\cdot x\), which gives the formula for \(\cc{M}\); its image is \([M]\). Thus \(\ab(M)=\Coker(\cc{M})=M/[M]=\overline{M}\), and \(M\) is abelian if and only if \([M]=0\), i.e., if and only if all brackets vanish. For such \(M\), \ref{M:sum} says that the action is additive in \(m\), \ref{M:P} says that it factors through \(\Rbar\), and \(M\) is an abelian group by Lemma~\ref{lem:module-basic}\ref{it:commutator}; conversely, any right \(\Rbar\)-module with zero bracket satisfies the axioms. Module morphisms between such objects are the \(\Rbar\)-linear maps.
\end{proof}

The short exact sequence \(0\to A(M,N)\to M+N\to M\times N\to 0\) which underlies Proposition~\ref{prop:cosmash-mod}, and thus the description of the cosmash product of two modules, appears as Lemma~3.6.1 in Pirashvili's report~\cite{Pirashvili-Square-Rings}, with a different proof: there it is deduced from the case of free modules by means of simplicial resolutions and derived functors.

\begin{theorem}\label{thm:mod-nilpotent}
  For any square ring \(\R\), the category \(\Mod_{\R}\) is a two-step nilpotent variety.
\end{theorem}

\begin{proof}
  By Proposition~\ref{prop:mod-semi-abelian} and Proposition~\ref{prop:quadratic-iff-bilinear}, it suffices to show that for each \(\R\)-module \(X\), the functor \(X\cosmash(-)\colon \Mod_{\R}\to\Mod_{\R}\) is linear. By Proposition~\ref{prop:cosmash-mod}, this functor is the composite
  \[
    \Mod_{\R}\xrightarrow{\;\ab\;}\Mod_{\Rbar}\xrightarrow{\;\overline{X}\otimes(-)\otimes_{\Rbar\otimes\Rbar}\Ree\;}\Mod_{\Rbar}=\Ab(\Mod_{\R})\hookrightarrow\Mod_{\R}\text{,}
  \]
  where on morphisms \(X\cosmash f\) is \(\overline{x}\otimes\overline{m}\otimes y\mapsto\overline{x}\otimes\overline{f(m)}\otimes y\). The abelianisation functor is a left adjoint, so it sends the coproduct \(M+N\) to the coproduct \(\overline{M}\oplus\overline{N}\) in the abelian category \(\Mod_{\Rbar}\), which is a biproduct; the tensor product functor is additive, so it preserves biproducts; and the inclusion of \(\Ab(\Mod_{\R})\) preserves products. Hence \(\cmF{X\cosmash(-)}{2}\colon X\cosmash(M+N)\to(X\cosmash M)\times(X\cosmash N)\) is an isomorphism, its kernel \(\cre_{2}(X\cosmash(-))(M,N)\) vanishes, and \(X\cosmash(-)\) is linear.
\end{proof}

\section{The square ring of a two-step nilpotent variety}\label{sec:square-ring-of-V}

Throughout this section, \(\V\) is a two-step nilpotent variety of algebras and \(E\) is the free algebra of \(\V\) on one generator \(e\). We write \(\U\colon \V\to\Set\) for the forgetful functor. Recall that \(E\) is a regular projective object, that the map
\[
  \V(E,X)\to\U(X)\colon \alpha\mapsto\alpha(e)
\]
is a bijection for every \(X\), and that \(E^{(n)}\DefEq E+\cdots+E\), the coproduct of \(n\) copies of \(E\), is the free algebra on \(n\) generators (with \(E^{(0)}=0\)). We write \(\Free(\V)\) for the full subcategory of \(\V\) determined by the objects \(E^{(n)}\), \(n\geq 0\); this is the algebraic theory of \(\V\), and morphisms \(E\to E^{(n)}\) are the \(n\)-ary terms of \(\V\). Since \(\V\) is semi-abelian, it is pointed, so \(0\) is its only constant. Our aim is to prove that
\[
  \R(\V)\DefEq\bigl(\V(E,E)\xrightarrow{H}\V(E,E\cosmash E)\xrightarrow{P}\V(E,E)\bigr)
\]
is a square ring, for suitable \(H\) and \(P\) and a suitable group structure on the hom-sets. This will be Theorem~\ref{thm:square-ring-of-V}.

\begin{remark}\label{rem:defourneau}
  In his thesis~\cite[Proposition~4.1.24]{Defourneau}, Defourneau shows that the algebraic theory \(\Free(\V)\) of a two-step nilpotent variety carries the structure of a square ring\emph{oid}. We only need the one-object case, and some of the results quoted in that proof are incomplete, while others are used in a generality in which they are not stated, so we give a self-contained proof.
\end{remark}

\subsection{Group structures on hom-sets}\label{sec:hom-groups}
Recall~\cite{Borceux-Bourn} that in a semi-abelian category an object \(Z\) is abelian if and only if it carries an internal group structure, which is then unique and abelian. Its addition \(m_{Z}\colon Z\times Z\to Z\) is the unique morphism with \(m_{Z}\circ\cm{Z}{Z}=\fold{Z}\), which exists because \(\fold{Z}\circ\iota_{Z,Z}=\cc{Z}=0\) and \(\cm{Z}{Z}=\coker(\iota_{Z,Z})\). For any object \(X\), the set \(\V(X,Z)\) is then an abelian group under \(f\boxplus g\DefEq m_{Z}\circ(f,g)\), with neutral element \(0\); the inverse is obtained by composing with the inverse of~\(Z\). Precomposition with any morphism \(X'\to X\) and postcomposition with any morphism \(Z\to Z'\) between abelian objects are group homomorphisms for \(\boxplus\), and \(\Ab(\V)\) is an additive category with this structure.

By Proposition~\ref{prop:cosmash-abelian}, all cosmash products \(Y\cosmash Z\) are abelian objects in \(\V\). In particular \(\V(E,E\cosmash E)\) is an abelian group under \(\boxplus\). Since \(E\) need not be abelian, this construction does not, in general, provide a group structure on \(\V(E,E)\); we use the freeness of \(E\) instead.

\begin{lemma}\label{lem:comonoid}
  There exists a morphism \(\mu_{E}\colon E\to E+E\) such that \(\cm{E}{E}\circ\mu_{E}=\diag{E}\), i.e., \(r_{1}\circ\mu_{E}=1_{E}=r_{2}\circ\mu_{E}\). Any such \(\mu_{E}\) is coassociative: \((\mu_{E}\cosum 1_{E})\circ\mu_{E}=(1_{E}\cosum\mu_{E})\circ\mu_{E}\). Hence \((E,\mu_{E},0)\) is a comonoid in \(\V\).
\end{lemma}

\begin{proof}
  Since \(\cm{E}{E}\) is a regular epimorphism (Subsection~\ref{sec:conventions}) and \(E\) is regular projective, \(\diag{E}\) factors through \(\cm{E}{E}\). For coassociativity, since \(E\cosmash E\cosmash E=0\), the morphism \(\hat{r}^{1_{\V}}_{3}=(r^{(1)},r^{(2)},r^{(3)})\colon E+E+E\to\prod_{k}(E+E)\) (where \(r^{(k)}\) sends the \(k\)-th term of the coproduct to zero) is a monomorphism, so it suffices to check that both sides agree after composition with each \(r^{(k)}\). Using \(r_{1}\circ\mu_{E}=1_{E}=r_{2}\circ\mu_{E}\), we find \(r^{(1)}\circ(\mu_{E}\cosum 1_{E})\circ\mu_{E}=((r_{2}\circ\mu_{E})\cosum 1_{E})\circ\mu_{E}=1_{E+E}\circ\mu_{E}=\mu_{E}\) and \(r^{(1)}\circ(1_{E}\cosum\mu_{E})\circ\mu_{E}=\langle 0,\mu_{E}\rangle\circ\mu_{E}=\mu_{E}\circ r_{2}\circ\mu_{E}=\mu_{E}\), so both composites with \(r^{(1)}\) equal \(\mu_{E}\); the cases \(k=2\), \(3\) are similar. Counitality \(\langle 1_{E},0\rangle\circ\mu_{E}=1_{E}=\langle 0,1_{E}\rangle\circ\mu_{E}\) is the defining property of \(\mu_{E}\).
\end{proof}

We fix, once and for all, a morphism \(\mu_{E}\) as in Lemma~\ref{lem:comonoid}. Note that \(\mu_{E}\) need not be unique: it corresponds to a binary term \(t(x,y)=\mu_{E}(e)\) of \(\V\) with \(t(x,0)=x=t(0,x)\), which is associative; see Remark~\ref{rem:choice-of-mu}. For any object \(Y\) of \(\V\), the comonoid structure induces a monoid structure on \(\V(E,Y)\): for \(f\), \(g\in\V(E,Y)\),
\begin{equation}\label{eq:def-sum}
  f+g\DefEq\langle f,g\rangle\circ\mu_{E}\text{,}
\end{equation}
with neutral element the zero morphism. Indeed,
\[
  (f+g)+h=\langle f,g,h\rangle\circ(\mu_{E}\cosum 1_{E})\circ\mu_{E}=\langle f,g,h\rangle\circ(1_{E}\cosum\mu_{E})\circ\mu_{E}=f+(g+h)
\]
by coassociativity, and \(0+f=\langle 0,f\rangle\circ\mu_{E}=f\circ r_{2}\circ\mu_{E}=f\), similarly \(f+0=f\). This structure is natural in \(Y\): for \(u\colon Y\to Y'\),
\begin{equation}\label{eq:post-additive}
  u\circ(f+g)=u\circ\langle f,g\rangle\circ\mu_{E}=\langle u\circ f,u\circ g\rangle\circ\mu_{E}=u\circ f+u\circ g\text{,}
\end{equation}
so postcomposition with \(u\) is a monoid homomorphism. Moreover, for an object \(X\) of \(\V\) and \(m\in\V(E,X)\), \(n\in\V(E,Y)\), we have \(\langle m,n\rangle=\fold{X}\circ(m\cosum n)\) when \(X=Y\), and \(m\cosum n\) restricts to \(m\cosmash n\colon E\cosmash E\to X\cosmash Y\):
\begin{equation}\label{eq:restriction}
  (m\cosum n)\circ\iota_{E,E}=\iota_{X,Y}\circ(m\cosmash n)\text{,}\qquad (m\circ r)\cosmash(n\circ s)=(m\cosmash n)\circ(r\cosmash s)\text{.}
\end{equation}

\begin{lemma}\label{lem:abelian-target}
  Let \(Z\) be an abelian object of \(\V\). Then the monoid structure~\eqref{eq:def-sum} on \(\V(E,Z)\) coincides with the abelian group structure \(\boxplus\). In particular, \(\V(E,Z)\) is an abelian group, and precomposition \(\V(E,Z)\to\V(E,Z)\colon a\mapsto a\circ t\) with any \(t\in\V(E,E)\) is a group homomorphism.
\end{lemma}

\begin{proof}
  For \(f\), \(g\in\V(E,Z)\), using~\eqref{eq:r-natural} and \(\cm{E}{E}\circ\mu_{E}=\diag{E}\),
  \begin{align*}
    f+g & =\fold{Z}\circ(f\cosum g)\circ\mu_{E} =m_{Z}\circ\cm{Z}{Z}\circ(f\cosum g)\circ\mu_{E}=m_{Z}\circ(f\times g)\circ\cm{E}{E}\circ\mu_{E} \\
        & =m_{Z}\circ(f\times g)\circ\diag{E}=m_{Z}\circ(f,g)=f\boxplus g\text{.}
  \end{align*}
  The last claim holds because precomposition is a homomorphism for \(\boxplus\).
\end{proof}

\begin{lemma}\label{lem:monoid-extension}
  Let \(\pi\colon G\to Q\) be a surjective homomorphism of monoids, written additively but not assumed commutative. If \(Q\) is a group and the submonoid \(K\DefEq\pi^{-1}(0)\) of \(G\)---the kernel of \(\pi\) in the category of monoids---is a group, then \(G\) is a group.
\end{lemma}

\begin{proof}
  Let \(g\in G\) and choose \(s\in G\) with \(\pi(s)=-\pi(g)\). Then \(\pi(g+s)=0\), so \(g+s\in K\) and there is \(k\in K\) with \(g+s+k=0\): thus \(g\) has a right inverse. Symmetrically, \(g\) has a left inverse. In any monoid, a left inverse \(a\) and a right inverse \(b\) of \(g\) coincide, since \(a=a+(g+b)=(a+g)+b=b\).
\end{proof}

\begin{proposition}\label{prop:hom-group}
  For every object \(Y\) of \(\V\), the monoid \((\V(E,Y),+,0)\) is a group. The inverse of \(f\in\V(E,Y)\) is \(-f=f\circ v_{E}\), where \(v_{E}\DefEq-1_{E}\) is the inverse of \(1_{E}\) in \(\V(E,E)\); in particular \((E,\mu_{E},v_{E},0)\) is a cogroup in \(\V\), and the functor \(\V(E,-)\colon \V\to\Set\) factors through the category of groups.
\end{proposition}

\begin{proof}
  Let \(K\subseteq\V(E,Y)\) be the image of the monoid homomorphism
  \[
    (\cc{Y})_{*}\colon \V(E,Y\cosmash Y)\to\V(E,Y)\colon \xi\mapsto\cc{Y}\circ\xi\text{.}
  \]
  Its domain is a group by Lemma~\ref{lem:abelian-target}, so \(K\) is a group. Consider the monoid homomorphism \((\ab_{Y})_{*}\colon \V(E,Y)\to\V(E,\ab(Y))\). Its codomain is a group by Lemma~\ref{lem:abelian-target}, and it is surjective because \(\ab_{Y}\) is a regular epimorphism and \(E\) is regular projective. We claim that \((\ab_{Y})_{*}^{-1}(0)=K\). Clearly \(\ab_{Y}\circ\cc{Y}=0\). Conversely, if \(\ab_{Y}\circ f=0\), then \(f\) factors through the kernel \([Y,Y]\rightarrowtail Y\) of \(\ab_{Y}=\coker(\cc{Y})\), which is the image of \(\cc{Y}\); since \(E\) is projective, \(f\) lifts along the regular epimorphism \(Y\cosmash Y\to[Y,Y]\), so \(f\in K\). By Lemma~\ref{lem:monoid-extension}, \(\V(E,Y)\) is a group. Finally, \(f+f\circ v_{E}=f\circ(1_{E}+v_{E})=f\circ 0=0\) by~\eqref{eq:post-additive}.
\end{proof}

Since the inverse is natural, \(\V(E,-)\) takes values in the category of groups, and the cogroup structure of \(E\) extends to all objects of \(\Free(\V)\): for \(X=E^{(n)}\) with coproduct injections \(i_{k}\colon E\to X\), let \(\mu_{X}\colon X\to X+X\) be defined by \(\mu_{X}\circ i_{k}=(i_{k}\cosum i_{k})\circ\mu_{E}\), and put \(f+g\DefEq\langle f,g\rangle\circ\mu_{X}\) for \(f\), \(g\in\V(X,Y)\). Then \((f+g)\circ i_{k}=f\circ i_{k}+g\circ i_{k}\), so that under the bijection \(\V(X,Y)\cong\V(E,Y)^{n}\colon v\mapsto(v\circ i_{k})_{k}\), this is the componentwise group structure. In particular \(\V(X,Y)\) is a group, \eqref{eq:post-additive} holds for \(f\), \(g\in\V(X,Y)\), and \(r_{1}\circ\mu_{X}=1_{X}=r_{2}\circ\mu_{X}\), so that \(\cm{X}{X}\circ\mu_{X}=\diag{X}\). Lemma~\ref{lem:abelian-target} extends right away to \(\V(X,Z)\) for \(X\) in \(\Free(\V)\) and \(Z\) abelian. Two further identities will be used constantly: for \(X\) in \(\Free(\V)\), with \(i_{1}\), \(i_{2}\colon X\to X+X\) the coproduct injections,
\begin{equation}\label{eq:mu-is-sum}
  \mu_{X}=\langle i_{1},i_{2}\rangle\circ\mu_{X}=i_{1}+i_{2}\text{,}
\end{equation}
and for any objects \(Y\), \(Z\) of \(\V\), with \(i_{1}\colon Y\to Y+Z\) and \(i_{2}\colon Z\to Y+Z\),
\begin{equation}\label{eq:c-of-coproduct-injections}
  \cc{Y+Z}\circ(i_{1}\cosmash i_{2})=\fold{Y+Z}\circ\iota_{Y+Z,Y+Z}\circ(i_{1}\cosmash i_{2})=\fold{Y+Z}\circ(i_{1}\cosum i_{2})\circ\iota_{Y,Z}=\iota_{Y,Z}\text{.}
\end{equation}

\begin{corollary}\label{cor:nearring}
  The set \(\V(E,E)\), with the addition~\eqref{eq:def-sum} and composition as multiplication, is a left nearring with unit \(1_{E}\).
\end{corollary}

\begin{proof}
  It is a group by Proposition~\ref{prop:hom-group} and a monoid under composition; left distributivity \(g\circ(f_{1}+f_{2})=g\circ f_{1}+g\circ f_{2}\) is~\eqref{eq:post-additive}.
\end{proof}

We stress that precomposition \(f\mapsto f\circ t\) is in general \emph{not} additive; the failure of right distributivity is measured by the map \(H\) constructed in Subsection~\ref{sec:HPT}.

\subsection{Bilinearity of the cosmash product and centrality}\label{sec:bilinearity}

\begin{proposition}\label{prop:cosmash-bilinear}
  Let \(X\), \(X'\) be objects of \(\Free(\V)\) and \(Y\), \(Y'\) objects of \(\V\). The map
  \[
    \V(X,Y)\times\V(X',Y')\to\V(X\cosmash X',Y\cosmash Y')\colon (f,g)\mapsto f\cosmash g
  \]
  is biadditive, the group structure on the right being \(\boxplus\).
\end{proposition}

\begin{proof}
  By Proposition~\ref{prop:cosmash-abelian} and symmetry, \(F\DefEq(-)\cosmash X'\colon \V\to\Ab(\V)\) is linear, so by Lemma~\ref{lem:linear-into-additive} we have \(1_{(X+X)\cosmash X'}=(i_{1}\cosmash 1_{X'})\circ(r_{1}\cosmash 1_{X'})+(i_{2}\cosmash 1_{X'})\circ(r_{2}\cosmash 1_{X'})\) in the additive category \(\Ab(\V)\). Hence, for \(f_{1}\), \(f_{2}\in\V(X,Y)\) and \(g\in\V(X',Y')\), using~\eqref{eq:restriction}, bilinearity of composition in \(\Ab(\V)\) and \(r_{k}\circ\mu_{X}=1_{X}\),
  \begin{align*}
    (f_{1}+f_{2})\cosmash g & =(\langle f_{1},f_{2}\rangle\cosmash g)\circ(\mu_{X}\cosmash 1_{X'})                                                                     \\
                            & =\sum_{k=1}^{2}(\langle f_{1},f_{2}\rangle\cosmash g)\circ(i_{k}\cosmash 1_{X'})\circ(r_{k}\cosmash 1_{X'})\circ(\mu_{X}\cosmash 1_{X'}) \\
                            & =\sum_{k=1}^{2}(\langle f_{1},f_{2}\rangle\circ i_{k}\circ r_{k}\circ\mu_{X})\cosmash g=f_{1}\cosmash g+f_{2}\cosmash g\text{.}
  \end{align*}
  Additivity in \(g\) is proved symmetrically.
\end{proof}

\begin{proposition}\label{prop:central}
  Let \(X\) be an object of \(\Free(\V)\), \(Y\) an object of \(\V\), \(f\in\V(X,Y)\) and \(\xi\in\V(X,Y\cosmash Y)\). Then \(f+\cc{Y}\circ\xi=\cc{Y}\circ\xi+f\). In other words, the homomorphism \((\cc{Y})_{*}\colon \V(X,Y\cosmash Y)\to\V(X,Y)\) takes values in the centre of \(\V(X,Y)\).
\end{proposition}

\begin{proof}
  By~\eqref{eq:restriction} and Corollary~\ref{cor:cosmash-abelianisation},
  \[
    \langle f,\cc{Y}\circ\xi\rangle\circ\iota_{X,X}=\fold{Y}\circ\iota_{Y,Y}\circ\bigl(f\cosmash(\cc{Y}\circ\xi)\bigr)=\cc{Y}\circ(1_{Y}\cosmash\cc{Y})\circ(f\cosmash\xi)=0\text{.}
  \]
  Since \(\cm{X}{X}=\coker(\iota_{X,X})\), there is a unique \(\phi_{1}\colon X\times X\to Y\) with \(\langle f,\cc{Y}\circ\xi\rangle=\phi_{1}\circ\cm{X}{X}\), and likewise a unique \(\phi_{2}\) with \(\langle\cc{Y}\circ\xi,f\rangle=\phi_{2}\circ\cm{X}{X}\). Now \(\twp{X}{X}\circ\cm{X}{X}=\cm{X}{X}\circ\tw{X}{X}\), so \(\phi_{1}\circ\twp{X}{X}\circ\cm{X}{X}=\langle f,\cc{Y}\circ\xi\rangle\circ\tw{X}{X}=\langle\cc{Y}\circ\xi,f\rangle=\phi_{2}\circ\cm{X}{X}\), whence \(\phi_{1}\circ\twp{X}{X}=\phi_{2}\) because \(\cm{X}{X}\) is an epimorphism. Using \(\cm{X}{X}\circ\mu_{X}=\diag{X}\) and \(\twp{X}{X}\circ\diag{X}=\diag{X}\),
  \[
    f+\cc{Y}\circ\xi=\phi_{1}\circ\cm{X}{X}\circ\mu_{X}=\phi_{1}\circ\twp{X}{X}\circ\diag{X}=\phi_{2}\circ\diag{X}=\phi_{2}\circ\cm{X}{X}\circ\mu_{X}=\cc{Y}\circ\xi+f\text{.}\qedhere
  \]
\end{proof}

\begin{corollary}\label{cor:iota-central}
  For \(X\) in \(\Free(\V)\), objects \(Y\), \(Z\) of \(\V\) and \(\xi\in\V(X,Y\cosmash Z)\), the morphism \(\iota_{Y,Z}\circ\xi\) is central in \(\V(X,Y+Z)\).
\end{corollary}

\begin{proof}
  By~\eqref{eq:c-of-coproduct-injections}, \(\iota_{Y,Z}\circ\xi=\cc{Y+Z}\circ(i_{1}\cosmash i_{2})\circ\xi\), which is central by Proposition~\ref{prop:central}.
\end{proof}

For an object \(Y\) of \(\V\), let \(\swc{Y}\colon Y\cosmash Y\to Y\cosmash Y\) be the restriction of the switch \(\tw{Y}{Y}\) to the kernel of \(\cm{Y}{Y}\); it exists because \(\twp{Y}{Y}\circ\cm{Y}{Y}=\cm{Y}{Y}\circ\tw{Y}{Y}\), and it is an involution. Since \(\tw{Y}{Y}\circ(f\cosum g)=(g\cosum f)\circ\tw{X}{X}\) and \(\fold{Y+Y}\circ(i_{2}\cosum i_{1})=\langle i_{2},i_{1}\rangle=\tw{Y}{Y}\), we have
\begin{equation}\label{eq:switch}
  \swc{Y}\circ(f\cosmash g)=(g\cosmash f)\circ\swc{X}\qquad\text{and}\qquad\cc{Y+Y}\circ(i_{2}\cosmash i_{1})=\iota_{Y,Y}\circ\swc{Y}
\end{equation}
for all \(f\), \(g\colon X\to Y\), the second identity being the analogue of~\eqref{eq:c-of-coproduct-injections} for the coproduct injections in the opposite order. (Note that \(i_{2}\cosmash i_{1}\) need not coincide with \((i_{1}\cosmash i_{2})\circ\swc{Y}\); they do agree after composition with \(\cc{Y+Y}\).)

\subsection{The maps \texorpdfstring{\(H\)}{H}, \texorpdfstring{\(P\)}{P} and \texorpdfstring{\(T\)}{T}, and the actions}\label{sec:HPT}

\begin{definition}\label{def:H}
  Let \(X\), \(Y\) be objects of \(\Free(\V)\) and \(\alpha\in\V(X,Y)\). Consider the morphism \(\mu_{Y}\circ\alpha-i_{2}\circ\alpha-i_{1}\circ\alpha\in\V(X,Y+Y)\). By~\eqref{eq:post-additive} and \(r_{k}\circ\mu_{Y}=1_{Y}\), its composites with \(r_{1}\) and \(r_{2}\) vanish, so it factors uniquely through \(\iota_{Y,Y}=\ker(\cm{Y}{Y})\). We define \(H_{X,Y}(\alpha)\in\V(X,Y\cosmash Y)\) to be the unique morphism such that
  \begin{equation}\label{eq:def-H}
    \iota_{Y,Y}\circ H_{X,Y}(\alpha)=(i_{1}+i_{2})\circ\alpha-i_{2}\circ\alpha-i_{1}\circ\alpha\text{.}
  \end{equation}
  We write \(H_{X}\DefEq H_{X,X}\) and \(H\DefEq H_{E,E}\colon \V(E,E)\to\V(E,E\cosmash E)\).
\end{definition}

By Corollary~\ref{cor:iota-central}, \(\iota_{Y,Y}\circ H_{X,Y}(\alpha)\) is central, so~\eqref{eq:def-H} may be rewritten as
\begin{equation}\label{eq:def-H-bis}
  (i_{1}+i_{2})\circ\alpha=i_{1}\circ\alpha+i_{2}\circ\alpha+\iota_{Y,Y}\circ H_{X,Y}(\alpha)\text{.}
\end{equation}
Thus \(H_{X,Y}(\alpha)\) measures the failure of \(\alpha\) to be ``additive''. The map \(H_{X,Y}\) is not a group homomorphism in general: its deviation is computed in Proposition~\ref{prop:H-sum}. We write \(2_{X}\DefEq 1_{X}+1_{X}\in\V(X,X)\) and \([a,b]\DefEq a+b-a-b\) for the commutator of two elements of a group.

\begin{lemma}\label{lem:H2}
  For \(X\) in \(\Free(\V)\), with \(i_{1}\), \(i_{2}\colon X\to X+X\):
  \[
    \iota_{X,X}\circ H_{X}(2_{X})=[i_{2},i_{1}]\qquad\text{and}\qquad\swc{X}\circ H_{X}(2_{X})=-H_{X}(2_{X})\text{.}
  \]
\end{lemma}

\begin{proof}
  By~\eqref{eq:def-H} and~\eqref{eq:post-additive},
  \[
    \iota_{X,X}\circ H_{X}(2_{X})=(i_{1}+i_{2})+(i_{1}+i_{2})-i_{2}-i_{2}-i_{1}-i_{1}=i_{1}+[i_{2},i_{1}]-i_{1}\text{,}
  \]
  and since the left-hand side is central (Corollary~\ref{cor:iota-central}), conjugating by \(i_{1}\) gives the first identity. For the second, \(\iota_{X,X}\circ\swc{X}\circ H_{X}(2_{X})=\tw{X}{X}\circ[i_{2},i_{1}]=[i_{1},i_{2}]=-[i_{2},i_{1}]=\iota_{X,X}\circ(-H_{X}(2_{X}))\) by~\eqref{eq:post-additive}, and \(\iota_{X,X}\) is a monomorphism.
\end{proof}

\begin{proposition}\label{prop:right-distributivity}
  Let \(X\), \(Y\) be objects of \(\Free(\V)\), \(Z\) an object of \(\V\), \(f\in\V(X,Y)\) and \(g_{1}\), \(g_{2}\in\V(Y,Z)\). Then
  \[
    (g_{1}+g_{2})\circ f=g_{1}\circ f+g_{2}\circ f+\cc{Z}\circ(g_{1}\cosmash g_{2})\circ H_{X,Y}(f)\text{.}
  \]
\end{proposition}

\begin{proof}
  By~\eqref{eq:mu-is-sum}, \eqref{eq:def-H-bis} and~\eqref{eq:post-additive}, \((g_{1}+g_{2})\circ f=\langle g_{1},g_{2}\rangle\circ(i_{1}+i_{2})\circ f=g_{1}\circ f+g_{2}\circ f+\langle g_{1},g_{2}\rangle\circ\iota_{Y,Y}\circ H_{X,Y}(f)\), and \(\langle g_{1},g_{2}\rangle\circ\iota_{Y,Y}=\fold{Z}\circ(g_{1}\cosum g_{2})\circ\iota_{Y,Y}=\fold{Z}\circ\iota_{Z,Z}\circ(g_{1}\cosmash g_{2})=\cc{Z}\circ(g_{1}\cosmash g_{2})\) by~\eqref{eq:restriction}.
\end{proof}

\begin{proposition}\label{prop:commutator}
  Let \(X\) be an object of \(\Free(\V)\), \(Y\) an object of \(\V\) and \(f\), \(g\in\V(X,Y)\). Then
  \[
    g+f=f+g+\cc{Y}\circ(f\cosmash g)\circ H_{X}(2_{X})\text{.}
  \]
\end{proposition}

\begin{proof}
  By~\eqref{eq:mu-is-sum}, \(g+f=\langle g,f\rangle\circ\mu_{X}=\langle f,g\rangle\circ\tw{X}{X}\circ(i_{1}+i_{2})=\langle f,g\rangle\circ(i_{2}+i_{1})\), and \(i_{2}+i_{1}=[i_{2},i_{1}]+i_{1}+i_{2}=i_{1}+i_{2}+\iota_{X,X}\circ H_{X}(2_{X})\) by Lemma~\ref{lem:H2} and centrality. Now apply~\eqref{eq:post-additive} and \(\langle f,g\rangle\circ\iota_{X,X}=\cc{Y}\circ(f\cosmash g)\).
\end{proof}

\begin{proposition}\label{prop:H-sum}
  Let \(X\), \(Y\) be objects of \(\Free(\V)\) and \(f\), \(g\in\V(X,Y)\). Then
  \[
    H_{X,Y}(f+g)=H_{X,Y}(f)+H_{X,Y}(g)+(g\cosmash f)\circ H_{X}(2_{X})\text{.}
  \]
\end{proposition}

\begin{proof}
  Write \(\iota=\iota_{Y,Y}\) and \(H=H_{X,Y}\). By~\eqref{eq:def-H} and~\eqref{eq:post-additive},
  \[
    \iota\circ H(f+g)=(i_{1}+i_{2})\circ f+(i_{1}+i_{2})\circ g-i_{2}\circ g-(i_{1}\circ g+i_{2}\circ f)-i_{1}\circ f\text{.}
  \]
  By Proposition~\ref{prop:commutator}, \eqref{eq:restriction}, \eqref{eq:switch} and Lemma~\ref{lem:H2},
  \begin{align*}
    i_{1}\circ g+i_{2}\circ f & =i_{2}\circ f+i_{1}\circ g+\cc{Y+Y}\circ(i_{2}\cosmash i_{1})\circ(f\cosmash g)\circ H_{X}(2_{X}) \\
                              & =i_{2}\circ f+i_{1}\circ g+\iota\circ\swc{Y}\circ(f\cosmash g)\circ H_{X}(2_{X})                  \\
                              & =i_{2}\circ f+i_{1}\circ g+\iota\circ(g\cosmash f)\circ\swc{X}\circ H_{X}(2_{X})                  \\
                              & =i_{2}\circ f+i_{1}\circ g-\iota\circ(g\cosmash f)\circ H_{X}(2_{X})\text{.}
  \end{align*}
  Substituting, and moving the central terms \(\iota\circ(g\cosmash f)\circ H_{X}(2_{X})\), \(\iota\circ H(g)\) (Corollary~\ref{cor:iota-central}) to the right,
  \begin{align*}
    \iota\circ H(f+g) & =(i_{1}+i_{2})\circ f+\bigl((i_{1}+i_{2})\circ g-i_{2}\circ g-i_{1}\circ g\bigr)-i_{2}\circ f-i_{1}\circ f \\
                      & \qquad+\iota\circ(g\cosmash f)\circ H_{X}(2_{X})                                                           \\
                      & =\iota\circ H(f)+\iota\circ H(g)+\iota\circ(g\cosmash f)\circ H_{X}(2_{X})\text{,}
  \end{align*}
  and the claim follows from~\eqref{eq:post-additive}, \(\iota\) being a monomorphism.
\end{proof}

\begin{proposition}\label{prop:H-composite}
  Let \(X\), \(Y\), \(Z\) be objects of \(\Free(\V)\), \(f\in\V(X,Y)\) and \(g\in\V(Y,Z)\). Then
  \[
    H_{X,Z}(g\circ f)=(g\cosmash g)\circ H_{X,Y}(f)+H_{Y,Z}(g)\circ f\text{.}
  \]
\end{proposition}

\begin{proof}
  Write \(\iota=\iota_{Z,Z}\). By~\eqref{eq:def-H-bis}, \((i_{1}+i_{2})\circ g\circ f=(i_{1}\circ g+i_{2}\circ g+\iota\circ H_{Y,Z}(g))\circ f\). We expand this using Proposition~\ref{prop:right-distributivity} twice. First,
  \[
    (i_{1}\circ g+i_{2}\circ g)\circ f=i_{1}\circ g\circ f+i_{2}\circ g\circ f+\cc{Z+Z}\circ(i_{1}\cosmash i_{2})\circ(g\cosmash g)\circ H_{X,Y}(f)\text{,}
  \]
  and \(\cc{Z+Z}\circ(i_{1}\cosmash i_{2})=\iota\) by~\eqref{eq:c-of-coproduct-injections}. Second, the correction term in
  \[
    \bigl((i_{1}\circ g+i_{2}\circ g)+\iota\circ H_{Y,Z}(g)\bigr)\circ f
  \]
  is \(\cc{Z+Z}\circ\bigl((i_{1}\circ g+i_{2}\circ g)\cosmash(\iota\circ H_{Y,Z}(g))\bigr)\circ H_{X,Y}(f)\), which vanishes because \(\iota\circ H_{Y,Z}(g)=\cc{Z+Z}\circ(i_{1}\cosmash i_{2})\circ H_{Y,Z}(g)\) and \(1_{Z+Z}\cosmash\cc{Z+Z}=0\) by Corollary~\ref{cor:cosmash-abelianisation}. Hence
  \[
    (i_{1}+i_{2})\circ g\circ f=i_{1}\circ g\circ f+i_{2}\circ g\circ f+\iota\circ(g\cosmash g)\circ H_{X,Y}(f)+\iota\circ H_{Y,Z}(g)\circ f\text{,}
  \]
  and comparing with~\eqref{eq:def-H-bis} for \(\alpha=g\circ f\) gives
  \[
    \iota\circ H_{X,Z}(g\circ f)=\iota\circ\bigl((g\cosmash g)\circ H_{X,Y}(f)+H_{Y,Z}(g)\circ f\bigr)
  \]
  by~\eqref{eq:post-additive}.
\end{proof}

\begin{proposition}\label{prop:T-formula}
  Let \(X\), \(Y\) be objects of \(\Free(\V)\) and \(\xi\in\V(X,Y\cosmash Y)\). Then
  \[
    H_{X,Y}(\cc{Y}\circ\xi)=\xi+\swc{Y}\circ\xi\text{.}
  \]
\end{proposition}

\begin{proof}
  Write \(\iota=\iota_{Y,Y}\) and \(c=\cc{Y+Y}\). By~\eqref{eq:def-H} and the naturality~\eqref{eq:c-natural} of \(\cc{}\),
  \[
    \iota\circ H_{X,Y}(\cc{Y}\circ\xi)=c\circ\bigl((i_{1}+i_{2})\cosmash(i_{1}+i_{2})\bigr)\circ\xi-c\circ(i_{2}\cosmash i_{2})\circ\xi-c\circ(i_{1}\cosmash i_{1})\circ\xi\text{.}
  \]
  By Proposition~\ref{prop:cosmash-bilinear}, \((i_{1}+i_{2})\cosmash(i_{1}+i_{2})\) is the sum, in the abelian group
  \[
    \V(Y\cosmash Y,(Y+Y)\cosmash(Y+Y))\text{,}
  \]
  of the four terms \(i_{k}\cosmash i_{l}\). Precomposition with \(\xi\) and postcomposition with \(c\) preserve sums (Lemma~\ref{lem:abelian-target} and~\eqref{eq:post-additive}), and all the terms \(c\circ(i_{k}\cosmash i_{l})\circ\xi\) are central by Proposition~\ref{prop:central}. Hence, by~\eqref{eq:c-of-coproduct-injections} and~\eqref{eq:switch},
  \[
    \iota\circ H_{X,Y}(\cc{Y}\circ\xi)=c\circ(i_{1}\cosmash i_{2})\circ\xi+c\circ(i_{2}\cosmash i_{1})\circ\xi=\iota\circ\xi+\iota\circ\swc{Y}\circ\xi=\iota\circ(\xi+\swc{Y}\circ\xi)\text{,}
  \]
  and \(\iota\) is a monomorphism.
\end{proof}

\begin{definition}\label{def:PT-actions}
  We define
  \begin{align*}
    P & \DefEq(\cc{E})_{*}\colon \V(E,E\cosmash E)\to\V(E,E)\colon x\mapsto\cc{E}\circ x=\fold{E}\circ\iota_{E,E}\circ x\text{,} \\
    T & \DefEq(\swc{E})_{*}\colon \V(E,E\cosmash E)\to\V(E,E\cosmash E)\colon x\mapsto\swc{E}\circ x\text{,}
  \end{align*}
  and, for \(r\), \(s\), \(t\in\V(E,E)\) and \(x\in\V(E,E\cosmash E)\),
  \[
    (r,s)\cdot x\cdot t\DefEq(r\cosmash s)\circ x\circ t\in\V(E,E\cosmash E)\text{.}
  \]
\end{definition}

\begin{lemma}\label{lem:PT-actions}
  The maps \(P\) and \(T\) are group homomorphisms, \(T\circ T=1_{\Ree}\), and \((r,s,x,t)\mapsto(r,s)\cdot x\cdot t\) is a quadrilinear map making \(\V(E,E\cosmash E)\) a left module over the multiplicative monoid \(\V(E,E)\times\V(E,E)\) and a right module over the monoid \(\V(E,E)\), the two actions commuting.
\end{lemma}

\begin{proof}
  \(P\) and \(T\) are homomorphisms by~\eqref{eq:post-additive}, and \(\swc{E}\) is an involution. The action axioms hold by functoriality of \(\cosmash\) (\(1_{E}\cosmash 1_{E}=1_{E\cosmash E}\) and \((r\cosmash s)\circ(r'\cosmash s')=(r\circ r')\cosmash(s\circ s')\) by~\eqref{eq:restriction}) and associativity of composition. Additivity in \(r\) and in \(s\) follows from Proposition~\ref{prop:cosmash-bilinear}, since pre- and postcomposition preserve \(\boxplus\) between abelian objects and \(\boxplus=+\) on \(\V(E,E\cosmash E)\) by Lemma~\ref{lem:abelian-target}; additivity in \(x\) holds by~\eqref{eq:post-additive} and Lemma~\ref{lem:abelian-target}; additivity in \(t\) holds by~\eqref{eq:post-additive}.
\end{proof}

\subsection{The square ring axioms}\label{sec:square-ring-axioms}

\begin{theorem}\label{thm:square-ring-of-V}
  Let \(\V\) be a two-step nilpotent variety, \(E\) its free algebra on one generator with a fixed comultiplication \(\mu_{E}\) as in Lemma~\ref{lem:comonoid}. With the nearring structure of Corollary~\ref{cor:nearring} on \(\Re\DefEq\V(E,E)\), the abelian group \(\Ree\DefEq\V(E,E\cosmash E)\), the map \(H\) of Definition~\ref{def:H} and the maps \(P\), \(T\) and the actions of Definition~\ref{def:PT-actions},
  \[
    \R(\V)=\sqring{\V(E,E)}{\V(E,E\cosmash E)}{H}{P}
  \]
  is a square ring.
\end{theorem}

\begin{proof}
  The data are as required by Lemma~\ref{lem:PT-actions}. Let \(r\), \(s\), \(t\in\Re\) and \(x\), \(y\in\Ree\). We verify the axioms of Definition~\ref{def:square-ring}.
  \begin{itemize}
    \item \ref{SR:P-binatural}: \(P((r,r)\cdot x\cdot s)=\cc{E}\circ(r\cosmash r)\circ x\circ s=r\circ\cc{E}\circ x\circ s=r\cdot P(x)\cdot s\) by~\eqref{eq:c-natural}.
    \item \ref{SR:T}: \(H(P(x))=x+T(x)\) is Proposition~\ref{prop:T-formula} with \(X=Y=E\).
    \item \ref{SR:PT}: \(P(T(x))=\fold{E}\circ\iota_{E,E}\circ\swc{E}\circ x=\fold{E}\circ\tw{E}{E}\circ\iota_{E,E}\circ x=\fold{E}\circ\iota_{E,E}\circ x=P(x)\), since \(\fold{E}\circ\tw{E}{E}=\fold{E}\).
    \item \ref{SR:P-kills}: \((P(x),r)\cdot y=(\cc{E}\cosmash 1_{E})\circ(x\cosmash r)\circ y=0\) by~\eqref{eq:restriction} and Corollary~\ref{cor:cosmash-abelianisation}, and likewise \((r,P(x))\cdot y=0\). Since \(E\cosmash E\) is abelian, \(y\colon E\to E\cosmash E\) factors through \(\ab_{E}\), and \(\ab_{E}\circ\cc{E}=0\), so \(y\cdot P(x)=y\circ\cc{E}\circ x=0\).
    \item \ref{SR:H-sum} is Proposition~\ref{prop:H-sum}, \ref{SR:H-product} is Proposition~\ref{prop:H-composite}, and \ref{SR:right-distributivity} is Proposition~\ref{prop:right-distributivity}, all with \(X=Y=Z=E\).
    \item \ref{SR:T-action}: \(T((r,s)\cdot x\cdot t)=\swc{E}\circ(r\cosmash s)\circ x\circ t=(s\cosmash r)\circ\swc{E}\circ x\circ t=(s,r)\cdot T(x)\cdot t\) by~\eqref{eq:switch}.\qedhere
  \end{itemize}
\end{proof}

\begin{remark}\label{rem:choice-of-mu}
  The square ring \(\R(\V)\) depends on the choice of \(\mu_{E}\), i.e., on the choice of a binary term \(t(x,y)\) of \(\V\) satisfying \(t(x,0)=x=t(0,x)\) (which forces \(t\) to be associative and to admit inverses, by Lemma~\ref{lem:comonoid} and Proposition~\ref{prop:hom-group}). Two such choices \(\mu_{E}\), \(\mu'_{E}\) yield square rings \(\R(\V)\), \(\R'(\V)\) whose module categories are both equivalent to \(\V\) by the Main Theorem, and hence to each other. The relation between \(\R(\V)\) and \(\R'(\V)\) is thus an instance of what could be called a ``quadratic Morita equivalence''.
\end{remark}

\section{The equivalence}\label{sec:equivalence}

We keep the setting and notation of Section~\ref{sec:square-ring-of-V}: \(\V\) is a two-step nilpotent variety, \(E\) its free algebra on the generator \(e\), \(\mu_{E}\) is fixed, and \(\R=\R(\V)\) is the square ring of Theorem~\ref{thm:square-ring-of-V}, so that \(\Re=\V(E,E)\) and \(\Ree=\V(E,E\cosmash E)\).

\subsection{The functor \texorpdfstring{\(\Phi\colon \V\to\Mod_{\R}\)}{Phi}}\label{sec:Phi}

\begin{proposition}\label{prop:hom-module}
  For every object \(X\) of \(\V\), the group \(\V(E,X)\) of Proposition~\ref{prop:hom-group} is an \(\R\)-module under
  \[
    m\cdot r\DefEq m\circ r\qquad\text{and}\qquad\brk{m}{n}\cdot x\DefEq\langle m,n\rangle\circ\iota_{E,E}\circ x
  \]
  for \(m\), \(n\in\V(E,X)\), \(r\in\Re\) and \(x\in\Ree\). For \(f\colon X\to Y\) in \(\V\), postcomposition \(f_{*}\colon \V(E,X)\to\V(E,Y)\) is a morphism of \(\R\)-modules. This defines a faithful functor \(\Phi\DefEq\V(E,-)\colon \V\to\Mod_{\R}\).
\end{proposition}

\begin{proof}
  Note first that \(\langle m,n\rangle\circ\iota_{E,E}=\fold{X}\circ(m\cosum n)\circ\iota_{E,E}=\cc{X}\circ(m\cosmash n)\) by~\eqref{eq:restriction}. We check the axioms of Definition~\ref{def:module}.
  \begin{itemize}
    \item \ref{M:action}: \(m\circ 1_{E}=m\), composition is associative, and \(m\circ(r+s)=m\circ r+m\circ s\) by~\eqref{eq:post-additive}.
    \item \ref{M:sum}: by Proposition~\ref{prop:right-distributivity}, \((m+n)\circ r=m\circ r+n\circ r+\cc{X}\circ(m\cosmash n)\circ H(r)=m\cdot r+n\cdot r+\brk{m}{n}\cdot H(r)\).
    \item \ref{M:P}: \(m\circ P(x)=m\circ\fold{E}\circ\iota_{E,E}\circ x=\langle m,m\rangle\circ\iota_{E,E}\circ x=\brk{m}{m}\cdot x\).
    \item \ref{M:T}: \(\brk{m}{n}\cdot T(x)=\langle m,n\rangle\circ\iota_{E,E}\circ\swc{E}\circ x=\langle m,n\rangle\circ\tw{E}{E}\circ\iota_{E,E}\circ x=\langle n,m\rangle\circ\iota_{E,E}\circ x=\brk{n}{m}\cdot x\).
    \item \ref{M:compatible}: \(\brk{m\circ r}{n\circ s}\cdot x=\langle m,n\rangle\circ(r\cosum s)\circ\iota_{E,E}\circ x=\langle m,n\rangle\circ\iota_{E,E}\circ(r\cosmash s)\circ x=\brk{m}{n}\cdot((r,s)\cdot x)\) by~\eqref{eq:restriction}, and \((\brk{m}{n}\cdot x)\cdot r=\brk{m}{n}\cdot(x\cdot r)\) by associativity of composition.
    \item \ref{M:trilinear}: additivity in \(x\) is~\eqref{eq:post-additive}. For additivity in \(m\), by Proposition~\ref{prop:cosmash-bilinear} we have \((m+m')\cosmash n=(m\cosmash n)\boxplus (m'\cosmash n)\); precomposition with \(x\) and postcomposition with \(\cc{X}\) preserve sums by Lemma~\ref{lem:abelian-target} and~\eqref{eq:post-additive}, whence \(\brk{m+m'}{n}\cdot x=\cc{X}\circ(m\cosmash n)\circ x+\cc{X}\circ(m'\cosmash n)\circ x\) follows. Additivity in \(n\) is symmetric.
    \item \ref{M:nil}: if \(m=\brk{h}{z}\cdot y=\cc{X}\circ(h\cosmash z)\circ y\), then \(m\cosmash n=(\cc{X}\cosmash 1_{X})\circ((h\cosmash z)\circ y\cosmash n)=0\) by Corollary~\ref{cor:cosmash-abelianisation}, so \(\brk{m}{n}\cdot x=0\). By ~\ref{M:trilinear}, the same vanishing holds for every \(m\) in the subgroup generated by these brackets.
  \end{itemize}
  Postcomposition with \(f\) is additive by~\eqref{eq:post-additive}, and \(f\circ(m\circ r)=(f\circ m)\circ r\), \(f\circ\langle m,n\rangle\circ\iota_{E,E}\circ x=\langle f\circ m,f\circ n\rangle\circ\iota_{E,E}\circ x\). Functoriality is clear, and \(\Phi\) is faithful because \(E\) is a generator.
\end{proof}

\subsection{Normal forms}\label{sec:normal-forms}
The following elementary decomposition replaces~\cite[Corollary~2.10]{BHP} in our context.

\begin{lemma}\label{lem:decomposition}
  For all objects \(X\), \(Y\) of \(\V\), the map
  \[
    \V(E,X)\times\V(E,Y)\times\V(E,X\cosmash Y)\to\V(E,X+Y)\colon (a,b,z)\mapsto i_{1}\circ a+i_{2}\circ b+\iota_{X,Y}\circ z
  \]
  is a bijection, and the last summand is central in \(\V(E,X+Y)\).
\end{lemma}

\begin{proof}
  Given \(\xi\in\V(E,X+Y)\), put \(a\DefEq r_{1}\circ\xi\) and \(b\DefEq r_{2}\circ\xi\). By~\eqref{eq:post-additive}, the composites of \(\xi-i_{2}\circ b-i_{1}\circ a\) with \(r_{1}\) and \(r_{2}\) vanish, so \(\xi-i_{2}\circ b-i_{1}\circ a=\iota_{X,Y}\circ z\) for a unique \(z\), and \(\xi=i_{1}\circ a+i_{2}\circ b+\iota_{X,Y}\circ z\) because \(\iota_{X,Y}\circ z\) is central (Corollary~\ref{cor:iota-central}). Conversely, if \(\xi=i_{1}\circ a+i_{2}\circ b+\iota_{X,Y}\circ z\), then \(r_{1}\circ\xi=a\), \(r_{2}\circ\xi=b\), and \(z\) is determined since \(\iota_{X,Y}\) is a monomorphism.
\end{proof}

\begin{proposition}\label{prop:normal-form-V}
  For \(n\geq 0\), with \(i_{k}\colon E\to E^{(n)}\) the coproduct injections, every morphism \(f\colon E\to E^{(n)}\) can be written in a unique way as
  \begin{equation}\label{eq:normal-form-V}
    f=\sum_{k=1}^{n}i_{k}\circ r_{k}+\sum_{1\leq k<l\leq n}\langle i_{k},i_{l}\rangle\circ\iota_{E,E}\circ x_{kl}
  \end{equation}
  with \(r_{k}\in\Re\) and \(x_{kl}\in\Ree\), the second sum being central. In other words, \(\V(E,E^{(n)})\cong\Re^{n}\times\Ree^{\binom{n}{2}}\) as sets.
\end{proposition}

\begin{proof}
  By induction on \(n\); the cases \(n=0\) and \(n=1\) are trivial. Suppose the claim holds for \(n\) and write \(E^{(n+1)}=E^{(n)}+E\), with \(j_{1}\colon E^{(n)}\to E^{(n+1)}\) and \(j_{2}=i_{n+1}\colon E\to E^{(n+1)}\). By Lemma~\ref{lem:decomposition}, every \(f\colon E\to E^{(n+1)}\) is uniquely \(f=j_{1}\circ a+i_{n+1}\circ r_{n+1}+\iota_{E^{(n)},E}\circ z\) with \(a\in\V(E,E^{(n)})\), \(r_{n+1}\in\Re\) and \(z\in\V(E,E^{(n)}\cosmash E)\). Since \((-)\cosmash E\) is linear (Proposition~\ref{prop:cosmash-abelian}), Lemma~\ref{lem:linear-into-additive} gives \(1_{E^{(n)}\cosmash E}=\sum_{k=1}^{n}(i_{k}\cosmash 1_{E})\circ(r'_{k}\cosmash 1_{E})\), so \(z=\sum_{k}(i_{k}\cosmash 1_{E})\circ z_{k}\) with \(z_{k}\DefEq(r'_{k}\cosmash 1_{E})\circ z\in\Ree\) uniquely determined by \(z\) (Lemma~\ref{lem:abelian-target}). Hence, by~\eqref{eq:restriction} and~\eqref{eq:post-additive}, and since \((i_{k}\cosum 1_{E})=\langle j_{1}\circ i_{k},i_{n+1}\rangle=\langle i_{k},i_{n+1}\rangle\) as morphisms \(E+E\to E^{(n+1)}\),
  \[
    \iota_{E^{(n)},E}\circ z=\sum_{k=1}^{n}(i_{k}\cosum 1_{E})\circ\iota_{E,E}\circ z_{k}=\sum_{k=1}^{n}\langle i_{k},i_{n+1}\rangle\circ\iota_{E,E}\circ z_{k}\text{.}
  \]
  Applying \(j_{1}\) to the normal form of \(a\) (using~\eqref{eq:post-additive} and \(j_{1}\circ\langle i_{k},i_{l}\rangle=\langle i_{k},i_{l}\rangle\)) and putting \(x_{k,n+1}\DefEq z_{k}\), we obtain~\eqref{eq:normal-form-V} for \(f\); uniqueness follows from the uniqueness in each step.
\end{proof}

\begin{proposition}\label{prop:Theta}
  For \(n\geq 0\), there is a unique isomorphism of \(\R\)-modules
  \[
    \theta_{n}\colon \Re^{(n)}\to\Phi(E^{(n)})=\V(E,E^{(n)})
  \]
  with \(\theta_{n}(e_{k})=i_{k}\) for all \(k\); it sends the element~\eqref{eq:normal-form-module} of \(\Re^{(n)}\) to the morphism~\eqref{eq:normal-form-V} with the same coefficients \(r_{k}\) and \(x_{kl}\). For a module morphism \(v\colon \Re^{(m)}\to\Re^{(n)}\), let
  \[
    \Theta(v)\DefEq\langle\theta_{n}(v(e_{1})),\dots,\theta_{n}(v(e_{m}))\rangle\colon E^{(m)}\to E^{(n)}\text{,}
  \]
  the unique morphism with \(\Theta(v)\circ i_{k}=\theta_{n}(v(e_{k}))\). Then
  \begin{equation}\label{eq:Theta-theta}
    \Theta(v)\circ\theta_{m}(a)=\theta_{n}(v(a))\qquad\text{for all \(a\in\Re^{(m)}\),}
  \end{equation}
  and \(\Theta\colon \Free(\Mod_{\R})\to\Free(\V)\colon \Re^{(n)}\mapsto E^{(n)}\) is an equivalence of categories which preserves finite coproducts, where \(\Free(\Mod_{\R})\) is the full subcategory of \(\Mod_{\R}\) determined by the free modules \(\Re^{(n)}\), \(n\geq 0\).
\end{proposition}

\begin{proof}
  By Remark~\ref{rem:freeness}, there is a unique module morphism \(\theta_{n}\) with \(\theta_{n}(e_{k})=i_{k}\), and it sends \(e_{k}\cdot r_{k}\) to \(i_{k}\circ r_{k}\) and \(\brk{e_{k}}{e_{l}}\cdot x_{kl}\) to \(\brk{i_{k}}{i_{l}}\cdot x_{kl}=\langle i_{k},i_{l}\rangle\circ\iota_{E,E}\circ x_{kl}\). Hence it maps~\eqref{eq:normal-form-module} to~\eqref{eq:normal-form-V}, coordinate by coordinate, and it is a bijection by Corollary~\ref{cor:free-modules} and Proposition~\ref{prop:normal-form-V}.

  For~\eqref{eq:Theta-theta}, both \(a\mapsto\theta_{n}(v(a))\) and \(a\mapsto\Theta(v)\circ\theta_{m}(a)=\Phi(\Theta(v))(\theta_{m}(a))\) are module morphisms \(\Re^{(m)}\to\V(E,E^{(n)})\) which agree on the generators \(e_{k}\); so they agree. It follows that \(\Theta\) is a functor: \(\Theta(1_{\Re^{(n)}})\circ i_{k}=\theta_{n}(e_{k})=i_{k}\), and \(\Theta(v\circ w)\circ i_{k}=\theta_{n}(v(w(e_{k})))=\Theta(v)\circ\theta_{m}(w(e_{k}))=\Theta(v)\circ\Theta(w)\circ i_{k}\), the \(i_{k}\) being jointly epimorphic. It is fully faithful, since the composite bijection
  \[
    \Mod_{\R}(\Re^{(m)},\Re^{(n)})\cong(\Re^{(n)})^{m}\xrightarrow{\;\theta_{n}^{m}\;}\V(E,E^{(n)})^{m}\cong\V(E^{(m)},E^{(n)})
  \]
  is \(v\mapsto\Theta(v)\). It is essentially surjective, every object of \(\Free(\V)\) being of the form \(E^{(n)}=\Theta(\Re^{(n)})\), hence an equivalence, and it sends the coproduct injections
  \[
    \Re^{(m)}\to\Re^{(m)}+\Re^{(n)}=\Re^{(m+n)}
  \]
  (which map \(e_{k}\) to \(e_{k}\)) to the coproduct injections \(E^{(m)}\to E^{(m+n)}\), so it preserves finite coproducts.
\end{proof}

\begin{remark}\label{rem:quadratic-category}
  The category \(\Free(\V)\) is a \defn{quadratic category} in the sense of~\cite[Definition~2.1]{BHP}. Its hom-sets are groups (Proposition~\ref{prop:hom-group}); composition is linear in the second variable by~\eqref{eq:post-additive} and quadratic in the first by Proposition~\ref{prop:right-distributivity}, its deviation being bilinear by Proposition~\ref{prop:cosmash-bilinear} and central by Proposition~\ref{prop:central}; the object \(E^{(0)}=0\) is a zero object; and the coproducts \(E^{(m)}+E^{(n)}\) are \defn{quadratic biproducts}, which is to say that their coproduct injections and retractions are linear and satisfy \(i_{1}\circ r_{1}+i_{2}\circ r_{2}=1_{E^{(m)}+E^{(n)}}\). The last point holds because the group structure on \(\V(E^{(k)},Y)\) is componentwise (Subsection~\ref{sec:hom-groups}), so that both sides agree after composition with the coproduct injections \(E^{(k)}\to E^{(m)}+E^{(n)}\), which are jointly epimorphic. Proposition~\ref{prop:Theta} then says that \(\Free(\V)\) is equivalent, as a quadratic category, to the category of free \(\R(\V)\)-modules of finite rank.
\end{remark}

\subsection{The functor \texorpdfstring{\(\Psi\colon \Mod_{\R}\to\V\)}{Psi}}\label{sec:Psi}
Let \(M\) be an \(\R\)-module. We define an algebra \(\Psi M\) of \(\V\) whose underlying set is that of the group \(M\). Given an \(n\)-ary term \(q\colon {E\to E^{(n)}}\) of \(\V\), put \(a_{q}\DefEq\theta_{n}^{-1}(q)\in\Re^{(n)}\) and, with the notation of Remark~\ref{rem:freeness}, define
\begin{equation}\label{eq:interpretation}
  q^{\Psi M}\colon M^{n}\to M\colon \mathbf{m}=(m_{1},\dots,m_{n})\mapsto\ev_{\mathbf{m}}(a_{q})\text{.}
\end{equation}
For \(n=0\), the only term is \(0\colon E\to 0\), interpreted as the constant \(0\in M\).

\begin{lemma}\label{lem:substitution}
  Let \(q\colon E\to E^{(n)}\) and \(t_{1}\), \dots, \(t_{n}\colon E\to E^{(m)}\) be terms, and let \(u\DefEq\langle t_{1},\dots,t_{n}\rangle\colon E^{(n)}\to E^{(m)}\), so that \(u\circ q\colon E\to E^{(m)}\) is the term obtained by substituting the \(t_{k}\) into \(q\). Then, for all \(\mathbf{m}\in M^{m}\),
  \[
    (u\circ q)^{\Psi M}(\mathbf{m})=q^{\Psi M}\bigl(t_{1}^{\Psi M}(\mathbf{m}),\dots,t_{n}^{\Psi M}(\mathbf{m})\bigr)\text{.}
  \]
\end{lemma}
\begin{proof}
  By Proposition~\ref{prop:Theta}, \(u=\Theta(\overline{u})\) for a unique \(\R\)-module map \(\overline{u}\colon \Re^{(n)}\to\Re^{(m)}\), namely the one with \(\overline{u}(e_{k})=\theta_{m}^{-1}(t_{k})=a_{t_{k}}\). By~\eqref{eq:Theta-theta}, \(\theta_{m}(\overline{u}(a_{q}))=\Theta(\overline{u})\circ\theta_{n}(a_{q})=u\circ q\), so \(a_{u\circ q}=\overline{u}(a_{q})\): the representative of the substituted term is obtained by substituting the representatives. Now \(\ev_{\mathbf{m}}\circ\overline{u}\colon \Re^{(n)}\to M\) is a module morphism sending \(e_{k}\) to \(\ev_{\mathbf{m}}(a_{t_{k}})=t_{k}^{\Psi M}(\mathbf{m})\), so it equals \(\ev_{\mathbf{t}}\) with \(\mathbf{t}\DefEq(t_{k}^{\Psi M}(\mathbf{m}))_{k}\). Hence \((u\circ q)^{\Psi M}(\mathbf{m})=\ev_{\mathbf{m}}(\overline{u}(a_{q}))=\ev_{\mathbf{t}}(a_{q})=q^{\Psi M}(\mathbf{t})\).
\end{proof}

Fix a presentation of \(\V\) by a signature and identities. Each operation symbol \(\omega\) of arity \(n\) is a term \(\omega\colon E\to E^{(n)}\), and we interpret it on the underlying set of \(M\) by~\(\omega^{\Psi M}\). By Lemma~\ref{lem:substitution} and induction on the construction of terms, the interpretation of any term \(t\) in this algebra is \(t^{\Psi M}\). If \(t=t'\) is an identity of \(\V\), then \(t\) and \(t'\) are the same morphism \(E\to E^{(n)}\), so \(t^{\Psi M}=t'^{\Psi M}\). Hence \(\Psi M\) is an algebra of \(\V\), and its structure does not depend on the chosen presentation: every term of \(\V\) is interpreted by~\eqref{eq:interpretation}.

\begin{example}\label{ex:interpretation}
  For a unary term \(r\in\Re=\V(E,E)\), we have \(\theta_{1}(e_{1}\cdot r)=r\), so \(a_{r}=e_{1}\cdot r\) and \(r^{\Psi M}(m)=m\cdot r\). By Proposition~\ref{prop:normal-form-V}, a binary term is uniquely of the form \(q=i_{1}\circ r_{1}+i_{2}\circ r_{2}+\langle i_{1},i_{2}\rangle\circ\iota_{E,E}\circ x\), and then \(a_{q}=e_{1}\cdot r_{1}+e_{2}\cdot r_{2}+\brk{e_{1}}{e_{2}}\cdot x\), so \(q^{\Psi M}(m,n)=m\cdot r_{1}+n\cdot r_{2}+\brk{m}{n}\cdot x\). In particular, the comultiplication \(\mu_{E}=i_{1}+i_{2}\) is interpreted as the addition of \(M\), the coinverse \(v_{E}=-1_{E}\) as \(m\mapsto-m\), and the term \(\iota_{E,E}\circ x\), for \(x\in\Ree\), as \((m,n)\mapsto\brk{m}{n}\cdot x\). So the group structure on \(\Psi M\) induced by the cogroup \(E\) is the group structure of \(M\).
\end{example}

Let \(f\colon M\to N\) be a morphism of \(\R\)-modules. For \(\mathbf{m}\in M^{n}\), the module morphisms \(f\circ\ev_{\mathbf{m}}\) and \(\ev_{f(\mathbf{m})}\colon \Re^{(n)}\to N\) agree on the generators, hence are equal; evaluating at \(a_{q}\) gives \(f(q^{\Psi M}(\mathbf{m}))=q^{\Psi N}(f(\mathbf{m}))\). Thus \(f\) is a homomorphism \(\Psi f\colon \Psi M\to\Psi N\) of algebras of \(\V\), and since composition and identities are those of the underlying maps we obtain a functor \(\Psi\colon \Mod_{\R}\to\V\).

\subsection{Proof of the Main Theorem}\label{sec:proof-main}

\begin{lemma}\label{lem:ev-theta}
  Let \(X\) be an object of \(\V\), \(\alpha_{1}\), \dots, \(\alpha_{n}\in\V(E,X)\), \(\boldsymbol{\alpha}=(\alpha_{k})_{k}\) and \(u\DefEq\langle\alpha_{1},\dots,\alpha_{n}\rangle\colon E^{(n)}\to X\). Then \(\ev_{\boldsymbol{\alpha}}=u_{*}\circ\theta_{n}\colon \Re^{(n)}\to\Phi(X)\), and consequently \(q^{\Psi\Phi(X)}(\boldsymbol{\alpha})=u\circ q\) for every \(n\)-ary term \(q\).
\end{lemma}

\begin{proof}
  Both are module morphisms \(\Re^{(n)}\to\Phi(X)\) (Proposition~\ref{prop:hom-module}) sending \(e_{k}\) to \(\alpha_{k}\), since \(u\circ\theta_{n}(e_{k})=u\circ i_{k}=\alpha_{k}\). Hence \(q^{\Psi\Phi(X)}(\boldsymbol{\alpha})=\ev_{\boldsymbol{\alpha}}(a_{q})=u\circ\theta_{n}(a_{q})=u\circ q\).
\end{proof}

\begin{theorem}\label{thm:main}
  Let \(\V\) be a variety of algebras. The following are equivalent:
  \begin{tfae}
    \item\label{it:main-nil} \(\V\) is two-step nilpotent;
    \item\label{it:main-mod} there exists a square ring \(\R\) such that \(\V\simeq\Mod_{\R}\).
  \end{tfae}
  More precisely, if \(\V\) is two-step nilpotent, then the functors \(\Phi\colon \V\to\Mod_{\R(\V)}\) and \(\Psi\colon \Mod_{\R(\V)}\to\V\) are mutually quasi-inverse equivalences.
\end{theorem}

\begin{proof}
  \ref{it:main-mod} implies~\ref{it:main-nil} by Theorem~\ref{thm:mod-nilpotent}, since a category equivalent to a two-step nilpotent variety is one itself. Conversely, let \(\V\) be two-step nilpotent and \(\R=\R(\V)\). We show that \(\Phi\) and \(\Psi\) are quasi-inverse equivalences.

  \emph{The isomorphism \(\eta\colon \Phi\Psi\Rightarrow 1_{\Mod_{\R}}\).} For an \(\R\)-module \(M\), define \(\eta_{M}\colon \Phi\Psi(M)=\V(E,\Psi M)\to M\colon \alpha\mapsto\alpha(e)\). It is a bijection because \(E\) is free on \(e\). It is a module morphism: for \(\alpha\), \(\beta\in\V(E,\Psi M)\), any \(n\)-ary term \(q\) and \(\alpha_{1}\), \dots, \(\alpha_{n}\in\V(E,\Psi M)\), the homomorphism \(\langle\alpha_{1},\dots,\alpha_{n}\rangle\colon E^{(n)}\to\Psi M\) preserves the term \(q\), so \((\langle\alpha_{1},\dots,\alpha_{n}\rangle\circ q)(e)=q^{\Psi M}(\alpha_{1}(e),\dots,\alpha_{n}(e))\). Applying this to \(q=\mu_{E}\), to the unary term \(r\in\Re\) and to \(q=\iota_{E,E}\circ x\), Example~\ref{ex:interpretation} gives
  \[
    \eta_{M}(\alpha+\beta)=\alpha(e)+\beta(e)\text{,}\qquad\eta_{M}(\alpha\cdot r)=\alpha(e)\cdot r\text{,}\qquad\eta_{M}(\brk{\alpha}{\beta}\cdot x)=\brk{\alpha(e)}{\beta(e)}\cdot x\text{.}
  \]
  Naturality is clear: for a module morphism \(f\colon M\to N\), \(f(\alpha(e))=(\Psi f\circ\alpha)(e)\).

  \emph{The isomorphism \(\varepsilon\colon \Psi\Phi\Rightarrow 1_{\V}\).} For an object \(X\) of \(\V\), define
  \[
    \varepsilon_{X}\colon \Psi\Phi(X)\to X\colon \alpha\mapsto\alpha(e)\text{,}
  \]
  a bijection on underlying sets. For an \(n\)-ary term \(q\) and \(\boldsymbol{\alpha}\in\V(E,X)^{n}\), with \(u=\langle\alpha_{1},\dots,\alpha_{n}\rangle\), Lemma~\ref{lem:ev-theta} gives
  \[
    \varepsilon_{X}\bigl(q^{\Psi\Phi(X)}(\boldsymbol{\alpha})\bigr)=(u\circ q)(e)=q^{X}(\alpha_{1}(e),\dots,\alpha_{n}(e))=q^{X}\bigl(\varepsilon_{X}(\alpha_{1}),\dots,\varepsilon_{X}(\alpha_{n})\bigr)\text{,}
  \]
  so \(\varepsilon_{X}\) is a bijective homomorphism, hence an isomorphism in \(\V\). Naturality in \(X\) is again clear.
\end{proof}

\section{Syntactic characterisations}\label{sec:syntactic}

As for abelian varieties~\cite{Freyd} and semi-abelian varieties~\cite{BoJa}, the Main Theorem yields descriptions of two-step nilpotent varieties in terms of operations and identities.

\subsection{Signature and identities}\label{sec:signature}
Let \(\R\) be a square ring. By the proof of Proposition~\ref{prop:mod-semi-abelian}, the variety \(\Mod_{\R}\) is presented by the signature \(\Sigma_{\R}\) consisting of a constant~\(0\), a unary operation \(-\), a binary operation \(+\), a unary operation \(\rho_{r}\) for each \(r\in\Re\) and a binary operation \(\phi_{x}\) for each \(x\in\Ree\), subject to the group identities for \((0,+,-)\) and to the identities, for all \(r\), \(s\in\Re\) and \(x\), \(y\in\Ree\),
\begin{modaxioms}
  \item \(\rho_{1}(a)=a\), \(\rho_{s}(\rho_{r}(a))=\rho_{r\cdot s}(a)\), \(\rho_{r+s}(a)=\rho_{r}(a)+\rho_{s}(a)\);
  \item \(\rho_{r}(a+b)=\rho_{r}(a)+\rho_{r}(b)+\phi_{H(r)}(a,b)\);
  \item \(\rho_{P(x)}(a)=\phi_{x}(a,a)\);
  \item \(\phi_{T(x)}(a,b)=\phi_{x}(b,a)\);
  \item \(\phi_{x}(\rho_{r}(a),\rho_{s}(b))=\phi_{(r,s)\cdot x}(a,b)\), \(\rho_{r}(\phi_{x}(a,b))=\phi_{x\cdot r}(a,b)\);
  \item \(\phi_{x}(a+b,c)=\phi_{x}(a,c)+\phi_{x}(b,c)\), \(\phi_{x}(a,b+c)=\phi_{x}(a,b)+\phi_{x}(a,c)\) and \(\phi_{x+y}(a,b)=\phi_{x}(a,b)+\phi_{y}(a,b)\);
  \item \(\phi_{y}(\phi_{x}(a,b),c)=0\).
\end{modaxioms}
These are the axioms of Definition~\ref{def:module} with \(\rho_{r}(a)=a\cdot r\) and \(\phi_{x}(a,b)=\brk{a}{b}\cdot x\). The normalisation identities \(\rho_{0}(a)=0\), \(\rho_{r}(0)=0\), \(\phi_{0}(a,b)=0\) and \(\phi_{x}(0,b)=0=\phi_{x}(a,0)\) follow from ~\ref{M:action}, ~\ref{M:trilinear} and Lemma~\ref{lem:module-basic}\ref{it:zero}. Note that \(\phi_{x}(a,b)\) is not the group commutator of \(a\) and \(b\); by Lemma~\ref{lem:module-basic}\ref{it:commutator}, the commutator \(b+a-b-a\) is \(\phi_{H(2)}(a,b)\). By Corollary~\ref{cor:free-modules}, every \(n\)-ary term operation of \(\Mod_{\R}\) has a unique normal form
\begin{equation}\label{eq:normal-form-term}
  u(a_{1},\dots,a_{n})=\sum_{k=1}^{n}\rho_{r_{k}}(a_{k})+\sum_{1\leq k<l\leq n}\phi_{x_{kl}}(a_{k},a_{l})\text{,}
\end{equation}
where uniqueness is understood modulo the identities, i.e., as term operations. If \(\Ree=0\), the operations \(\phi_{x}\) disappear, \(\Re\) is a ring and the presentation is the usual one of right \(\Re\)-modules.

Recall that two varieties are \defn{term equivalent} if they have the same algebras up to a bijective correspondence between their term operations that preserves arities, projections and substitutions---equivalently, if their categories of algebras are isomorphic over \(\Set\). This is the appropriate notion here, since a signature is never unique: redundant operation symbols may always be added.

\begin{corollary}\label{cor:term-equivalence}
  A variety \(\V\) is two-step nilpotent if and only if there is a square ring \(\R\) such that \(\V\) is term equivalent to the variety presented by the signature \(\Sigma_{\R}\), the group identities and the identities~\ref{M:action}--\ref{M:nil} above.
\end{corollary}

\begin{proof}
  If \(\V\) is term equivalent to \(\Mod_{\R}\), then \(\V\cong\Mod_{\R}\) as categories and \(\V\) is two-step nilpotent by Theorem~\ref{thm:mod-nilpotent}. Conversely, let \(\V\) be two-step nilpotent, and let \(\R=\R(\V)\). By Proposition~\ref{prop:normal-form-V} and Lemma~\ref{lem:ev-theta}, for any algebra \(X\) of \(\V\) and any \(n\)-ary term \(q=\sum_{k}i_{k}\circ r_{k}+\sum_{k<l}\langle i_{k},i_{l}\rangle\circ\iota_{E,E}\circ x_{kl}\) of \(\V\), the term operation \(q^{X}\) is given by~\eqref{eq:normal-form-term} with \(\rho_{r}=r^{X}\), \(\phi_{x}=(\iota_{E,E}\circ x)^{X}\), and \(+\) the term operation \(\mu_{E}^{X}\). Hence every term operation of \(\V\) is a term operation in the derived signature \(\Sigma_{\R}\), whose members are themselves terms of \(\V\). By Theorem~\ref{thm:main} and Example~\ref{ex:interpretation}, the algebras of \(\V\) equipped with these derived operations are precisely the \(\R\)-modules, on the same underlying sets, and the \(\V\)-homomorphisms are the module morphisms. So \(\V\) is term equivalent to \(\Mod_{\R}\).
\end{proof}

\subsection{Quadratic term operations}\label{sec:quadratic-terms}
Following~\cite{BHP}, a map \(u\colon G^{n}\to G\), where \(G\) is a group written additively, is \defn{quadratic} if its \defn{deviation}
\[
  (\mathbf{a}\mid\mathbf{b})_{u}\DefEq u(\mathbf{a}+\mathbf{b})-u(\mathbf{b})-u(\mathbf{a})\qquad(\mathbf{a},\mathbf{b}\in G^{n})
\]
is biadditive and takes values in the centre of \(G\). A quadratic map is reduced, \(u(\mathbf{0})=0\): indeed \((\mathbf{0}\mid\mathbf{0})_{u}=-u(\mathbf{0})\) vanishes by biadditivity. In a two-step nilpotent variety, every term reduced in each of three disjoint blocks of variables vanishes:

\begin{lemma}\label{lem:trireduced}
  Let \(\V\) be a two-step nilpotent variety and \(t(\mathbf{x};\mathbf{y};\mathbf{z})\) a term whose variables are divided into three disjoint blocks. If the identities \(t(\mathbf{0};\mathbf{y};\mathbf{z})=0\), \(t(\mathbf{x};\mathbf{0};\mathbf{z})=0\) and \(t(\mathbf{x};\mathbf{y};\mathbf{0})=0\) hold in \(\V\), then so does \(t=0\).
\end{lemma}

\begin{proof}
  Let \(X\), \(Y\), \(Z\) be the free algebras on the variables in \(\mathbf{x}\), \(\mathbf{y}\), \(\mathbf{z}\), so that \(X+Y+Z\) is free on all of them, and consider the element \(t(\mathbf{x};\mathbf{y};\mathbf{z})\) of \(X+Y+Z\). The hypotheses say that it is sent to \(0\) by the three retractions \(r^{(k)}\), since homomorphisms preserve terms; as kernels in a variety are computed on underlying sets, it lies in \(X\cosmash Y\cosmash Z\) (see also~\cite[Example~2.6]{HVdL}), which vanishes. So \(t=0\) in the free algebra, i.e., \(t=0\) is an identity of \(\V\).
\end{proof}

\begin{corollary}\label{cor:quadratic-terms}
  For a variety \(\V\), the following are equivalent:
  \begin{tfae}
    \item\label{it:qt-nil} \(\V\) is two-step nilpotent;
    \item\label{it:qt-quad} there are terms \(0\), \(-x\) and \(x+y\) which make the underlying set of every algebra of \(\V\) a group, and such that with respect to this group structure, for every \(n\geq 1\), every \(n\)-ary term \(u\) and every algebra \(X\), the term operation \(u_{X}\colon \U(X)^{n}\to\U(X)\) is quadratic.
  \end{tfae}
\end{corollary}

\begin{proof}
  \ref{it:qt-nil} implies \ref{it:qt-quad}: take the group structure given by \(\mu_{E}\) (Proposition~\ref{prop:hom-group}, Example~\ref{ex:interpretation}). Since \(\V\) is pointed, \(u_{X}(\mathbf{0})=0\) for every term \(u\), hence \((\mathbf{0}\mid\mathbf{b})_{u}=0=(\mathbf{a}\mid\mathbf{0})_{u}\). For biadditivity in the first variable, consider the term \(A(\mathbf{a};\mathbf{a}';\mathbf{b})\DefEq(\mathbf{a}+\mathbf{a}'\mid\mathbf{b})_{u}-(\mathbf{a}'\mid\mathbf{b})_{u}-(\mathbf{a}\mid\mathbf{b})_{u}\), which vanishes when any of its three blocks of variables is zero; by Lemma~\ref{lem:trireduced}, \(A=0\) is an identity. Biadditivity in the second variable is symmetric. For centrality, consider \(C(\mathbf{a};\mathbf{b};d)\DefEq(\mathbf{a}\mid\mathbf{b})_{u}+d-(\mathbf{a}\mid\mathbf{b})_{u}-d\), which again vanishes when any block is zero, so \(C=0\) by Lemma~\ref{lem:trireduced}.

  \ref{it:qt-quad} implies \ref{it:qt-nil}: \(\V\) is protomodular by~\cite{BoJa}, since it admits group operations as terms. It is pointed: for a nullary term \(c\), the unary term \(x\mapsto c+x\) is quadratic, hence reduced, so \(c=c+0=0\). Thus \(\V\) is semi-abelian. It remains to see that \(X\cosmash Y\cosmash Z=0\) for all \(X\), \(Y\), \(Z\) in \(\V\). Put \(W\DefEq X+Y+Z\) and let \(h\in\U(W)\) be in the kernel of the three retractions. Since \(W\) is generated by the images of \(X\), \(Y\) and~\(Z\), we may write \(h=u_{W}(\mathbf{w})\) for some \(n\)-ary term \(u\) and some \(\mathbf{w}\in\U(W)^{n}\) each of whose components comes from one of the three factors. Decompose \(\mathbf{w}=\mathbf{x}+\mathbf{y}+\mathbf{z}\) componentwise, where \(\mathbf{x}\) keeps the components coming from \(X\) and replaces the others by \(0\), and similarly for \(\mathbf{y}\) and \(\mathbf{z}\). Writing \(j\colon X+Y\to W\) for the coproduct injection and \(r^{(3)}\colon W\to X+Y\) for the retraction, we have \(\mathbf{x}+\mathbf{y}=(j\circ r^{(3)})^{n}(\mathbf{w})\), hence \(u_{W}(\mathbf{x}+\mathbf{y})=j(r^{(3)}(u_{W}(\mathbf{w})))=j(r^{(3)}(h))=0\), and likewise \(u_{W}(\mathbf{x}+\mathbf{z})=u_{W}(\mathbf{y}+\mathbf{z})=0\) and \(u_{W}(\mathbf{x})=u_{W}(\mathbf{y})=u_{W}(\mathbf{z})=0\). By definition of the deviation, \((\mathbf{x}\mid\mathbf{z})_{u}=u_{W}(\mathbf{x}+\mathbf{z})-u_{W}(\mathbf{z})-u_{W}(\mathbf{x})=0\), and similarly \((\mathbf{y}\mid\mathbf{z})_{u}=0\). By biadditivity, \(h=u_{W}(\mathbf{x}+\mathbf{y}+\mathbf{z})=(\mathbf{x}+\mathbf{y}\mid\mathbf{z})_{u}+u_{W}(\mathbf{x}+\mathbf{y})+u_{W}(\mathbf{z})=(\mathbf{x}\mid\mathbf{z})_{u}+(\mathbf{y}\mid\mathbf{z})_{u}=0\).
\end{proof}

\begin{remark}
  The choice of the group structure in~\ref{it:qt-quad} need not be unique, as the comultiplication \(\mu_{E}\) need not be (Remark~\ref{rem:choice-of-mu}); some choice is needed for the deviation to make sense, much as in the syntactic characterisation of abelian varieties.
\end{remark}

\begin{remark}
  Being two-step nilpotent, every category of modules over a square ring is algebraically coherent by~\cite[Theorem~5.2]{acc}; being moreover a variety, it is even (fibre-wise) algebraically cartesian closed~\cite{Gray2014}.
\end{remark}

\section{Compatibility with the abelian core; role of the bilinear product}\label{sec:compatibility}
Let \(\V\) be an arbitrary semi-abelian variety with free algebra \(E\) on one generator. The abelian core \(\Ab(\V)\) is an abelian variety, whose free algebra on one generator is \(\ab(E)\), so that \(\Ab(\V)\simeq\Mod_{R}\) for the ring \(R\DefEq\Ab(\V)(\ab(E),\ab(E))\)~\cite{Freyd}. Likewise, the two-nilpotent core \(\Nil_{2}(\V)\) is a two-step nilpotent variety with free algebra \(E'\DefEq\nil_{2}(E)\) on one generator, so that \(\Nil_{2}(\V)\simeq\Mod_{\R}\) for the square ring
\[
  \R\DefEq\R(\Nil_{2}(\V))=\sqring{\V(E',E')}{\V(E',E'\cosmash_{2}E')}{H}{P}\text{,}
\]
where \(\cosmash_{2}\) denotes the cosmash product in \(\Nil_{2}(\V)\). The ring \(R\) is encoded in the square ring \(\R\):

\begin{proposition}\label{prop:compatibility}
  In the situation above, \(R\cong\Rbar=\Coker(P)\) as rings.
\end{proposition}

\begin{proof}
  Since \(E'\) is projective in \(\Nil_{2}(\V)\) and \(\ab\cong\ab\circ\nil_{2}\) (both \(\ab\) and \(\nil_{2}\) are left adjoints to inclusions, and \(\Ab(\V)\subseteq\Nil_{2}(\V)\)), we have \(\ab(E)\cong\ab(E')\), and the adjunction gives a bijection
  \[
    R=\Ab(\V)(\ab(E'),\ab(E'))\cong\V(E',\ab(E'))
  \]
  sending \(\phi\) to \(\phi\circ\ab_{E'}\), under which the addition of \(R\) corresponds to the addition of Subsection~\ref{sec:hom-groups} (Lemma~\ref{lem:abelian-target}). Consider the map \(\V(E',E')\to R\colon f\mapsto\ab(f)\), which corresponds to \((\ab_{E'})_{*}\colon \V(E',E')\to\V(E',\ab(E'))\) since \(\ab(f)\circ\ab_{E'}=\ab_{E'}\circ f\). It is multiplicative by functoriality of \(\ab\), additive by~\eqref{eq:post-additive}, and surjective because \(E'\) is projective: given \(\phi\), the morphism \(\phi\circ\ab_{E'}\) lifts along the regular epimorphism \(\ab_{E'}\) to some \(f\), and then \(\ab(f)=\phi\). Its kernel is the image of \(P=(\cc{E'})_{*}\), by the argument in the proof of Proposition~\ref{prop:hom-group}. Hence it induces a ring isomorphism \(\Rbar=\Coker(P)\to R\).
\end{proof}

This agrees with Proposition~\ref{prop:cosmash-mod}: \(\Ab(\Nil_{2}(\V))=\Ab(\V)\simeq\Ab(\Mod_{\R})\simeq\Mod_{\Rbar}\).

In~\cite{DHVdL} an intrinsic tensor product, the \defn{bilinear product}, is introduced on any semi-abelian category \(\C\): the bilinear product of two objects \(X\) and \(Y\) is the abelian object \(X\otimes Y\DefEq\ab(X)\cosmash_{2}\ab(Y)\) of \(\C\), where \(\cosmash_{2}\) denotes the cosmash product in \(\Nil_{2}(\C)\). Its restriction to \(\Ab(\C)\) is a \defn{symmetric bi-cocontinuous bifunctor} \(\otimes\colon \Ab(\C)\times\Ab(\C)\to\Ab(\C)\). Now let \(\C=\V\). The equivalence \(\Ab(\V)\simeq\Mod_{R}\) of~\cite{Freyd} sends an abelian object \(A\) to the right \(R\)-module \(\Ab(\V)(\ab(E),A)\), and under it \(X\otimes Y\) corresponds to the \(R\)-module
\[
  \bigl(\V(\ab(E),\ab(X))\otimes\V(\ab(E),\ab(Y))\bigr)\otimes_{R\otimes R}\V(\ab(E),\ab(E)\cosmash_{2}\ab(E))\text{,}
\]
naturally in \(X\) and \(Y\). Here \(\V(\ab(E),\ab(E)\cosmash_{2}\ab(E))\cong\V(E',E'\cosmash_{2}E')=\Ree\), since \(\ab(E)\cong\ab(E')\) and \(E'\cosmash_{2}E'\) only depends on \(\overline{E'}=\ab(E')\) by Proposition~\ref{prop:cosmash-mod}. Moreover, the symmetry isomorphism \(\tau_{X,Y}\colon X\otimes Y\to Y\otimes X\) corresponds to the map \(f\otimes g\otimes h\mapsto g\otimes f\otimes T(h)\). It follows that the pair \((\Rbar,\Ree)\), where \(\Ree\) carries its structure of an \((\Rbar\wr\mathfrak{S}_{2})\)-\(\Rbar\)-bimodule from Subsection~\ref{sec:square-rings-def}, encodes the abelian core of \(\V\) together with its bilinear product. Note that in~\cite{DHVdL} the ring \(\Rbar\wr\mathfrak{S}_{2}\) is written \((\Rbar\otimes\Rbar)\wr\mathfrak{S}_{2}\).

In~\cite{DHVdL} the converse problem is raised.

\begin{question}\label{q:realisation}
  Given an abelian variety \(\A\) and a symmetric bi-cocontinuous bifunctor \(B\colon \A\times\A\to\A\), is there a two-step nilpotent variety \(\V\) whose abelian core is equivalent to \(\A\) and whose bilinear product is isomorphic to \(B\) under this equivalence? If so, determine all such varieties \(\V\) up to equivalence.
\end{question}

The existence of \(\V\) is settled in the affirmative in~\cite[Theorem~7.6]{DHVdL}, as follows. For a free algebra \(F\) on a single generator of \(\A\), the abelian group \(\A(F,B(F,F))\) is an \((R\wr\mathfrak{S}_{2})\)-\(R\)-bimodule, where \(R\) is the endomorphism ring of \(F\) in \(\A\)~\cite[proof of Theorem~7.6]{DHVdL}; so the pair \(\P=(R,\A(F,B(F,F)))\) is a \emph{two-step nilpotent reduced symmetric operad in the monoidal category of abelian groups}~\cite[Remark~6.6]{DHVdL}, and the category of \(\P\)-algebras is a two-step nilpotent variety with the required properties. Thanks to Theorem~\ref{thm:main} and the above observations, we can now settle the second half of Question~\ref{q:realisation}.

\begin{theorem}\label{thm:classification}
  Any two-step nilpotent variety \(\V\) such that \(\Ab(\V)\) is equivalent to \(\A\) and the bilinear product of \(\V\) is isomorphic to \(B\) under this equivalence is the category of modules over a square ring \(\R\) such that the ring \(\Rbar\) is isomorphic to \(R\) and \(\Ree\) is isomorphic to \(\A(F,B(F,F))\) as an \((\Rbar\wr\mathfrak{S}_{2})\)-\(\Rbar\)-bimodule.\noproof
\end{theorem}

In subsequent work we develop a means to construct such square rings, with prescribed ring \(\Rbar\) and bimodule \(\Ree\), from classical (i.e., linear and abelian) data, using the abelian models of nilpotent groups introduced in the second author's PhD thesis \cite{Hartl-Diss} and recalled in a more general form in his habilitation thesis~\cite{Hartl-Habil}.

Finally, given an arbitrary square ring \(\R\), one may ask what the role of the two-step nilpotent reduced symmetric operad \(\P=(\Rbar,\Ree)\) in abelian groups is with respect to the category of \(\R\)-modules. The answer is that it encodes the multilinear structure of the graded abelian group formed by the quotients of the lower central series of an \(\R\)-module: the multilinear operations acting on it and the multilinear identities they satisfy. This goes as follows.

For an \(\R\)-module \(M\), let \(\Gr(M)\) be the graded abelian group
\[
  \Gr(M)\DefEq\gamma_{1}(M)/\gamma_{2}(M)\oplus\gamma_{2}(M)/\gamma_{3}(M)=\overline{M}\oplus[M]\text{,}
\]
where \(\gamma_{2}(M)=[M,M]=[M]\) by Proposition~\ref{prop:cosmash-mod} and \(\gamma_{3}(M)=0\) by Theorem~\ref{thm:mod-nilpotent}. To see that \(\Gr(M)\) is naturally a \(\P\)-algebra, we again use \cite[Remark~6.6]{DHVdL}. Both summands of \(\Gr(M)\) are right \(\Rbar\)-modules: \(\overline{M}\) by Lemma~\ref{lem:module-basic}\ref{it:quotient}, and \([M]\) by Proposition~\ref{prop:cosmash-mod}, being the image of \(\cc{M}\colon M\cosmash M\to M\). Let \(D(\Gr(M))\DefEq 0\oplus[M]\subseteq\Gr(M)\), so that \(\Gr(M)/D(\Gr(M))\cong\overline{M}\). A surjective homomorphism of right \(\Rbar\)-modules
\[
  \overline{\mu}_{2}\colon (\overline{M}\otimes\overline{M})\otimes_{\Rbar\wr\mathfrak{S}_{2}}\Ree\to[M]\colon \overline{m}\otimes\overline{n}\otimes x\mapsto\brk{m}{n}\cdot x\text{,}
\]
where \(\tau\) acts on \(\overline{M}\otimes\overline{M}\) by the swap, is well defined by Proposition~\ref{prop:cosmash-mod} and Definition~\ref{def:module}\ref{M:T}: the identity \(\brk{n}{m}\cdot x=\brk{m}{n}\cdot T(x)\) says precisely that the map of Proposition~\ref{prop:cosmash-mod} factors through the tensor product over \(\Rbar\wr\mathfrak{S}_{2}\). It is easy to check that \(\Gr\) extends to a functor from \(\Mod_{\R}\) to the category of graded \(\P\)-algebras.

This observation is a special case of the global multilinear algebra structure on the quotients of the lower central series in semi-abelian varieties exhibited in~\cite{CCC}.

\section*{Acknowledgements}
This work originates in an idea due to Teimuraz Pirashvili~\cite{Pirashvili-Square-Rings}, whom we thank warmly. Generative AI tools (OpenAI's ChatGPT 6.1 Sol and Anthropic's Claude Fable 5.1) were used to assist in revising and restructuring a draft written by the authors, with the aim of improving its presentation; they were not used to generate new mathematical results. The authors reviewed and edited all output and take full responsibility for the content of this article.

\bibliographystyle{amsplain-nodash}
\bibliography{tim}

\providecommand{\noopsort}[1]{}
\providecommand{\bysame}{\leavevmode\hbox to3em{\hrulefill}\thinspace}
\providecommand{\MR}{\relax\ifhmode\unskip\space\fi MR }
\providecommand{\MRhref}[2]{%
  \href{http://www.ams.org/mathscinet-getitem?mr=#1}{#2}
}
\providecommand{\href}[2]{#2}
\begin{thebibliography}{10}

\bibitem{Amaglobeli-Miasnikov}
M.~Amaglobeli and A.~Miasnikov, \emph{Tensor completions of {$2$}-nilpotent
  finitely generated torsion-free groups}, Preprint,
  \href{https://arxiv.org/abs/2601.17600}{arXiv:2601.17600}, 2026.

\bibitem{Barratt-TrackGroups}
M.~G. Barratt, \emph{Track groups. {II}}, Proc. Lond. Math. Soc. (3) \textbf{5}
  (1955), 285--329.

\bibitem{Baues:Metastable}
H.-J. Baues, \emph{Quadratic functors and metastable homotopy}, J. Pure Appl.
  Algebra \textbf{91} (1994), no.~1-3, 49--107.

\bibitem{Baues:HomotopyType}
H.-J. Baues, \emph{Homotopy type and homology}, Oxford Math. Monogr., Oxford
  Univ. Press, 1996.

\bibitem{BHP}
H.-J. Baues, M.~Hartl, and T.~Pirashvili, \emph{Quadratic categories and square
  rings}, J. Pure Appl. Algebra \textbf{122} (1997), no.~1-2, 1--40.

\bibitem{Baues-Iwase}
H.-J. Baues and N.~Iwase, \emph{Square rings associated to elements in homotopy
  groups of spheres}, Groups of homotopy self-equivalences and related topics,
  Contemp. Math., vol. 274, 2001, pp.~57--78.

\bibitem{BJP:SquareGroups}
H.-J. Baues, M.~Jibladze, and T.~Pirashvili, \emph{Quadratic algebra of square
  groups}, Adv. Math. \textbf{217} (2008), no.~3, 1236--1300.

\bibitem{Baues-Muro}
H.-J. Baues and F.~Muro, \emph{Secondary homotopy groups}, Forum Math.
  \textbf{20} (2008), no.~4, 631--677.

\bibitem{Baues-Muro:Smash}
H.-J. Baues and F.~Muro, \emph{Smash products for secondary homotopy groups},
  Appl. Categ. Structures \textbf{16} (2008), no.~5, 551--616.

\bibitem{Baues-Muro:RingSpectra}
H.-J. Baues and F.~Muro, \emph{The algebra of secondary homotopy operations in
  ring spectra}, Proc. Lond. Math. Soc. (3) \textbf{102} (2011), no.~4,
  637--696.

\bibitem{Baues-Pirashvili}
H.-J. Baues and T.~Pirashvili, \emph{Quadratic endofunctors of the category of
  groups}, Adv.~Math. \textbf{141} (1999), no.~1, 167--206.

\bibitem{Baues-Pirashvili:Spaces}
H.-J. Baues and T.~Pirashvili, \emph{Spaces associated to quadratic
  endofunctors of the category of groups}, Extracta Math. \textbf{20} (2005),
  no.~2, 99--136.

\bibitem{Baumslag:UniqueRoots}
G.~Baumslag, \emph{Some aspects of groups with unique roots}, Acta Math.
  \textbf{104} (1960), 217--303.

\bibitem{BeBou}
C.~Berger and D.~Bourn, \emph{Central reflections and nilpotency in exact
  {M}al'tsev categories}, J. Homotopy Relat. Struct. \textbf{12} (2017),
  765--835.

\bibitem{Borceux-Bourn}
F.~Borceux and D.~Bourn, \emph{Mal'cev, protomodular, homological and
  semi-abelian categories}, Math. Appl., vol. 566, Kluwer Acad. Publ., 2004.

\bibitem{Bourn-Gran}
D.~Bourn and M.~Gran, \emph{Central extensions in semi-abelian categories},
  J.~Pure Appl. Algebra \textbf{175} (2002), 31--44.

\bibitem{BoJa}
D.~Bourn and G.~Janelidze, \emph{Characterization of protomodular varieties of
  universal algebras}, Theory Appl. Categ. \textbf{11} (2003), No. 6, 143--147.

\bibitem{CJ03}
A.~Carboni and G.~Janelidze, \emph{Smash product of pointed objects in
  lextensive categories}, J. Pure Appl. Algebra \textbf{183} (2003), no.~1--3,
  27--43.

\bibitem{acc}
A.~S. Cigoli, J.~R.~A. Gray, and T.~Van~der Linden, \emph{Algebraically
  coherent categories}, Theory Appl. Categ. \textbf{30} (2015), no.~54,
  1864--1905.

\bibitem{Defourneau}
T.~Defourneau, \emph{Linéarisation de structures algébriques à l'aide
  d'opérades et de foncteurs polynomiaux}, Ph.D. thesis, Université de
  Valenciennes et du Hainaut-Cambrésis, 2017.

\bibitem{DHVdL}
B.~S. Deval, M.~Hartl, and T.~Van~der Linden, \emph{Intrinsic tensor products
  and a {G}anea-type extension of the five-term exact sequence}, J.\ Algebra
  \textbf{710} (2027), 645--715.

\bibitem{Eilenberg-MacLane}
S.~Eilenberg and S.~{Mac\,Lane}, \emph{On the groups {{\(H(\pi,n)\)}} {II}},
  Ann. of Math. \textbf{60} (1954), 49--139.

\bibitem{Freese-McKenzie}
R.~Freese and R.~McKenzie, \emph{Commutator theory for congruence modular
  varieties}, London Math. Soc. Lecture Note Ser., vol. 125, Cambridge Univ.
  Press, 1987.

\bibitem{Freyd}
P.~Freyd, \emph{Abelian categories}, Harper and Row, New York, 1964, Repr.
  Theory Appl. Categ. \textbf{3} (2003).

\bibitem{P2RM}
H.~Gaudier and M.~Hartl, \emph{Quadratic maps between modules}, J. Algebra
  \textbf{322} (2009), 1667--1688.

\bibitem{Gaudier-Hartl:Morphismes}
H.~Gaudier and M.~Hartl, \emph{Morphismes quadratiques entre modules sur un
  anneau carr{\'e}}, Preprint,
  \href{https://arxiv.org/abs/1001.2849}{arXiv:1001.2849}, 2010.

\bibitem{Hartl-WIP}
F.~Goichot, M.~Hartl, and B.~Loiseau, \emph{Foundations of quadratic algebra},
  unpublished notes, LAMAV, Universit{\'e} de Valenciennes, 2010.

\bibitem{Gran-Rosicky:Monadic}
M.~Gran and J.~Rosick{\'y}, \emph{Semi-abelian monadic categories}, Theory
  Appl. Categ. \textbf{13} (2004), no.~6, 106--113.

\bibitem{Gray2014}
J.~R.~A. Gray, \emph{A note on the distributivity of the {H}uq commutator over
  finite joins}, Appl.\ Categ.\ Structures \textbf{22} (2014), 305--310.

\bibitem{Hall:Edmonton}
P.~Hall, \emph{The {E}dmonton notes on nilpotent groups}, Queen Mary College
  Mathematics Notes, Queen Mary College, London, 1969, Lectures given at the
  University of Alberta, 1957.

\bibitem{Hartl-Linearization}
M.~Hartl, \emph{Linearization of algebraic structures via functor calculus}, in
  preparation.

\bibitem{Hartl-Diss}
M.~Hartl, \emph{Abelsche {M}odelle nilpotenter {G}ruppen},
  Inauguraldissertation, Universit\"at Bonn, 1991.

\bibitem{Hartl-Habil}
M.~Hartl, \emph{Structures polynomiales en th{\'e}orie des groupes nilpotents},
  M{\'e}moire d'habilitation, Universit{\'e} Louis Pasteur, Strasbourg, 1998,
  \url{https://homepages.abdn.ac.uk/d.j.benson/pages/html/archive/hartl.html}.

\bibitem{Hartl:QuadraticMaps}
M.~Hartl, \emph{Quadratic maps between groups}, Georgian Math. J. \textbf{16}
  (2009), no.~1, 55--74.

\bibitem{CCC}
M.~Hartl, \emph{Algebraic functor calculus and commutator theory}, in
  preparation, 2019.

\bibitem{WorkingMathematician}
M.~Hartl and B.~Loiseau, \emph{A characterization of semi-abelian categories},
  Cah. Topol. G{\'e}om. Diff{\'e}r. Cat{\'e}g. \textbf{52} (2011), no.~1,
  77--80.

\bibitem{Actions}
M.~Hartl and B.~Loiseau, \emph{On actions and strict actions in homological
  categories}, Theory Appl.~Categ. \textbf{27} (2013), no.~15, 347--392.

\bibitem{HPV}
M.~Hartl, T.~Pirashvili, and C.~Vespa, \emph{Polynomial functors from algebras
  over a set-operad and nonlinear {M}ackey functors}, Int. Math. Res. Not. IMRN
  \textbf{2015} (2015), no.~6, 1461--1554.

\bibitem{HVdL}
M.~Hartl and T.~Van~der Linden, \emph{The ternary commutator obstruction for
  internal crossed modules}, Adv.~Math. \textbf{232} (2013), no.~1, 571--607.

\bibitem{Hartl-Vespa}
M.~Hartl and C.~Vespa, \emph{Quadratic functors on pointed categories},
  Adv.~Math. \textbf{226} (2011), no.~5, 3927--4010.

\bibitem{Higgins}
P.~J. Higgins, \emph{Groups with multiple operators}, Proc. Lond. Math. Soc.
  (3) \textbf{6} (1956), 366--416.

\bibitem{Huq}
S.~A. Huq, \emph{Commutator, nilpotency and solvability in categories}, Q. J.
  Math. \textbf{19} (1968), no.~2, 363--389.

\bibitem{IKK}
S.~O. Ivanov, G.~Kadantsev, and A.~Krasilnikov, \emph{Ring-induced
  localizations of nilpotent groups}, 2026, Preprint,
  \href{https://arxiv.org/abs/2606.15392}{arXiv:2606.15392}.

\bibitem{Janelidze-Kelly}
G.~Janelidze and G.~M. Kelly, \emph{Galois theory and a general notion of
  central extension}, J.~Pure Appl. Algebra \textbf{97} (1994), no.~2,
  135--161.

\bibitem{Janelidze-Marki-Tholen}
G.~Janelidze, L.~M{\'a}rki, and W.~Tholen, \emph{Semi-abelian categories},
  J.~Pure Appl. Algebra \textbf{168} (2002), no.~2--3, 367--386.

\bibitem{Jibladze-Pirashvili:Linear}
M.~Jibladze and T.~Pirashvili, \emph{Linear extensions and nilpotence of
  {M}altsev theories}, Beitr. Algebra Geom. \textbf{46} (2005), no.~1, 71--102.

\bibitem{Jibladze-Pirashvili}
M.~Jibladze and T.~Pirashvili, \emph{Quadratic envelope of the category of
  class two nilpotent groups}, Georgian Math. J. \textbf{13} (2006), no.~4,
  693--722.

\bibitem{Lyndon1960}
R.~C. Lyndon, \emph{Groups with parametric exponents}, Trans. Amer. Math. Soc.
  \textbf{96} (1960), 518--533.

\bibitem{MM}
S.~Mantovani and G.~Metere, \emph{Internal crossed modules and {P}eiffer
  condition}, Theory Appl. Categ. \textbf{23} (2010), no.~6, 113--135.

\bibitem{MR94}
A.~G. Myasnikov and V.~N. Remeslennikov, \emph{Groups with exponents {I}.
  {F}undamentals of the theory and tensor completions}, Siberian Math. J.
  \textbf{35} (1994), no.~5, 986--996.

\bibitem{Passi:GroupRings}
I.~B.~S. Passi, \emph{Group rings and their augmentation ideals}, Lecture Notes
  in Math., vol. 715, Springer, 1979.

\bibitem{Pirashvili-Square-Rings}
T.~Pirashvili, \emph{Report on square rings}, unpublished manuscript.

\bibitem{Sm}
J.~D.~H. Smith, \emph{{M}al'cev varieties}, Lecture Notes in Math., vol. 554,
  Springer, 1976.

\bibitem{Warfield:Nilpotent}
R.~B. Warfield, Jr., \emph{Nilpotent groups}, Lecture Notes in Math., vol. 513,
  Springer, 1976.

\end{thebibliography}

\end{document}